\documentclass[a4paper,11pt]{article}
\usepackage[a4paper,left=1.0in,right=1.0in,top=1.0in,bottom=1.00in,nohead]{geometry}

\usepackage[utf8]{inputenc}
\usepackage[T1]{fontenc}

\usepackage{hyperref}

\usepackage{amsmath}
\usepackage{amssymb}
\usepackage{amsthm}
\usepackage{amsfonts}
\usepackage{mathtools}
\usepackage{bigints}
\usepackage{xcolor}
\usepackage{enumitem}
\usepackage{afterpage}
\usepackage{setspace}
\usepackage{cleveref}
\usepackage[numbers]{natbib}

\newtheorem{theorem}{Theorem}[section]
\newtheorem{proposition}{Proposition}[section]
\newtheorem{corollary}{Corollary}[section]
\newtheorem{lemma}{Lemma}[section]

\theoremstyle{definition}
\newtheorem{definition}{Definition}[section]

\newtheorem{example}{Example}[section]
\newtheorem{rmk}{Remark}[section]

\numberwithin{equation}{section}

\newcommand{\B}{\mathcal{B}}
\newcommand{\C}{\mathcal{C}}
\newcommand{\E}{\mathbb{E}}
\newcommand{\e}{\text{e}}
\newcommand{\F}{\mathcal{F}}
\newcommand{\TF}{\tilde{\mathcal{F}}}
\newcommand{\G}{\mathcal{G}}
\newcommand{\M}{\mathcal{M}}
\newcommand{\N}{\mathbb{N}}
\newcommand{\OO}{\mathcal{O}}
\newcommand{\TO}{\tilde{\Omega}}
\newcommand{\PP}{\mathcal{P}}
\newcommand{\PB}{\mathbb{P}}
\newcommand{\Q}{\mathcal{Q}}
\newcommand{\QB}{\mathbb{Q}}
\newcommand{\R}{\mathbb{R}}
\newcommand{\Rd}{\mathbb{R}^d}

\newcommand{\UNN}{\mathcal{U}^N_N}
\newcommand{\UNone}{\mathcal{U}^N_1}

\newcommand{\UTNN}{\tilde{\mathcal{U}}^N_N}
\newcommand{\UTNone}{\tilde{\mathcal{U}}^N_1}
\newcommand{\V}{\mathcal{V}}
\newcommand{\wiener}{\mathcal{W}^{\nu}}
\newcommand{\X}{\mathcal{X}}

\newcommand{\bu}{\textbf{\textit{u}}}
\newcommand{\blambda}{\boldsymbol{\lambda}}
\newcommand{\Ni}{^{N,i}}
\newcommand{\None}{^{N,1}}

\newcommand{\rbr}[1]{\left(#1\right)}

\newcommand{\cbr}[1]{\left\{#1\right\}}
\newcommand{\nor}[1]{\left\|#1\right\|}
\newcommand{\modl}[1]{\left| #1\right|}
\newcommand{\infnorm}[1]{\| #1\|_{\infty}}
\newcommand{\infnormt}[1]{\| #1\|_{\infty,t}}
\newcommand{\infnorms}[1]{\| #1\|_{\infty,s}}

\newcommand{\round}[1]{\left(#1\right)}
\newcommand{\brackets}[1]{\left[#1\right]}
\newcommand{\braces}[1]{\left\{#1\right\}}
\newcommand{\conv}{\underset{n\rightarrow\infty}{\longrightarrow}}

\newcommand{\weakconv}{\overset{w}{\rightharpoonup}}

\newcommand{\simdistr}{\overset{d}{\sim}}
\newcommand{\ba}{\begin{eqnarray}}
\newcommand{\ea}{\end{eqnarray}}
\newcommand{\bi}{\begin{itemize}}
\newcommand{\ei}{\end{itemize}}
\newcommand{\ben}{\begin{enumerate}}
\newcommand{\een}{\end{enumerate}}
\newcommand{\blem}{\begin{lemma}}
\newcommand{\elem}{\end{lemma}}
\newcommand{\bteo}{\begin{theorem}}
\newcommand{\eteo}{\end{theorem}}
\newcommand{\bdefi}{\begin{definition}}
\newcommand{\edefi}{\end{definition}}
\newcommand{\bcor}{\begin{corollary}}
\newcommand{\ecor}{\end{corollary}}
\newcommand{\bprop}{\begin{proposition}}
\newcommand{\eprop}{\end{proposition}}
\newcommand{\brem}{\begin{rmk}}
\newcommand{\erem}{\end{rmk}}

\hypersetup{
    colorlinks=true,
    linkcolor=blue,
    filecolor=magenta,      
    urlcolor=black,
    pdftitle={Sharelatex Example},
    bookmarks=true,
}

\begin{document}
\title{\textbf{$N$-player games and mean-field games with \\ smooth dependence on past absorptions}\\
}
\author{Luciano Campi\thanks{%
London School of Economics, Department of Statistics, Columbia House, Houghton Street, London, WC2A 2AE.
Universit\`{a} degli Studi di Milano, Dipartimento di Matematica  ``Federigo Enriques'', Via Saldini 50, 20133, Milano, Italy.
Email:  L.Campi@lse.ac.uk.
} \and Maddalena Ghio\thanks{%
Scuola Normale Superiore, Piazza dei Cavalieri 7, 56126, Pisa. Email: maddalena.ghio@sns.it.} \and Giulia Livieri\thanks{%
Scuola Normale Superiore, Piazza dei Cavalieri 7, 56126, Pisa. Email: giulia.livieri@sns.it.
}}
\date{\today}
\maketitle
\begin{abstract}
\noindent Mean-field games with absorption is a class of games that has been introduced in \cite{campi2018n} and that can be viewed as natural limits of symmetric stochastic differential games with a large number of players who, interacting through a mean-field, leave the game as soon as their private states hit some given boundary.

In this paper, we push the study of such games further, extending their scope along two main directions. First, we allow the state dynamics and the costs to have a very general, possibly infinite-dimensional, dependence on the (non-normalized) empirical sub-probability measure of the survivors' states. This includes the particularly relevant case where the mean-field interaction among the players is done through the empirical measure of the survivors together with the fraction of absorbed players over time. Second, the boundedness of coefficients and costs has been considerably relaxed including drift and costs with linear growth in the state variables, hence allowing for more realistic dynamics for players' private states. We prove the existence of solutions of the MFG in strict as well as relaxed feedback form, and we establish uniqueness of the MFG solutions under monotonicity conditions of Lasry-Lions type. Finally, we show in a setting with finite-dimensional interaction  that such solutions induce approximate Nash equilibria for the $N$-player game with vanishing error as $N \to \infty$.\medskip

\noindent{\small \textbf{Key words and phrases:} Nash equilibrium, mean-field game, absorbing boundary, McKean-Vlasov limit, controlled martingale problem, relaxed control.}\medskip 

\noindent{\small \textbf{2000 AMS subject classifications:} 60B10, 60K35, 91A06, 93E20.}
\end{abstract}
\newpage
{
  \hypersetup{hidelinks}
  \tableofcontents
}
\section{Introduction}\label{sec:introduction}
\noindent Mean-field games (MFGs for short) are, loosely speaking, limits of symmetric stochastic differential games with a large number of players, where each of them interacts with the average behaviour of his/her competitors. They were introduced in the seminal papers by \citet{lasry-lions2006a,lasry-lions2006b,lasry-lions2007} and, simultaneously, by \citet{huang2006}. An increasing stream of research has been flourishing since then, producing theoretical results as well as a wide range of applications in many fields such as economics, finance, crowd dynamics and social sciences in general. For an excellent presentation of the theory we refer to the lecture notes of \citet{cardaliaguet} and the two-volume monograph by \citet{carmonadelarue}.\medskip
	
\noindent\textit{Motivation.}\quad In most of the literature on MFGs, all players stay in the game until the end of the period, while in many applications, especially in economics and finance, it is natural to have a mechanism deciding when some player has to leave. Such a mechanism can be modelled by introducing an absorbing boundary for the state space as in \citet{campi2018n}, which is the starting point of our study (other related references will be discussed later in detail). Therein, existence of solutions of the MFG and construction of approximate Nash equilibria for the $N$-player games were provided under some boundedness assumptions on the coefficients and without including the effect of past absorption on the survivors' behaviour. The present paper continues the investigation of this kind of games, with the following main extensions.
\bi
\item[(i)] We recast MFGs with absorption in a more general setting, most common to the MFG literature, where the dependence of the dynamics and costs on the empirical measure is infinite-dimensional.
\item[(ii)] We introduce a direct dependence on past absorptions in the drift of the Stochastic Differential Equations (SDEs) describing the evolution of the players' states by letting the initial distribution of players lose mass over time. Such a loss of mass corresponds to the exit of the absorbed players from the game, so that the proportion of the absorbed players has an effect on the future evolution of the survivors. This feature was not present in \cite{campi2018n}, where the empirical measure of the survivors was re-normalized at each time. Such a dependence on past absorptions is also included in the costs.
\item[(iii)] We allow both the drift and the cost functional of the players to grow at most linearly with the state, hence they are not necessarily bounded unlike in \cite{campi2018n}. Moreover, the set of non-absorbing states $\mathcal O$ can also be unbounded. Dropping the boundedness of the game data increases the flexibility of our setting, which can include more realistic dynamics from the viewpoint of applications (for more details, see later in this introduction). 
\ei
\noindent To be more precise, the purpose of this paper is to study $N$-player games and related MFGs in the presence of an absorbing set (i.e. a player is eliminated from the game once his/her private state leaves a given open set $\mathcal{O} \in \mathbb{R}^{d}$), and where the vector of private states $\textbf{X}^{N} \doteq (X^{N,1}, \ldots, X^{N,N})$ evolves according to
\ba 
X\Ni_t = X\Ni_0 &+& \int_{0}^{t} \bar{b}\round{s, X\Ni_s,\mu_s^N, u\Ni\round{s,\textbf{X}^N}}\,ds+ \sigma W\Ni_t, \,\,\, t \in \brackets{0,T},
\label{eq:U:SDE}
\ea
\noindent for $i \in \braces{1,\ldots ,N}$, where $\bu^N \doteq (u^{N,1},\ldots ,u^{N,N})$ is a vector of feedback strategies, $W^{N,1},\ldots,$ $W^{N,N}$ are independent $d$-dimensional Wiener processes defined on some filtered probability space, $\sigma$ is the (non-degenerate) diffusion matrix and $\bar{b}$ is a given  drift functional.  Finally, $\mu^{N}$ is the random flow of empirical sub-probability measures representing the empirical distribution of the survivors
\ba 
\mu_t^N\round{\cdot} \doteq \frac{1}{N}\sum_{i = 1}^{N} \delta_{X_t^{N,i}}\round{\cdot}\mathbf{1}_{[0,\tau^{X^{N,i}})}\round{t}.
\nonumber
\ea
\noindent Each player evaluates a strategy vector $\bu^N$ according to his/her expected costs
\ba
J^{N,i}\round{\bu^N} \doteq \E\Biggl[\int_{0}^{\tau\Ni}\bar{f}\round{s, X_s\Ni,\mu_s^N, u\Ni\round{s,\textbf{X}^N}}ds
+F\round{\tau\Ni,X\Ni_{\tau\Ni}}\Biggr]
\label{eq:U:Cost}
\ea
\noindent over a random time horizon. In Eq.\eqref{eq:U:Cost}, $\textbf{X}^N$ is the $N$-player dynamics under $\bu^N$ and $\tau^{N,i} \doteq \tau^{X^{N,i}} \wedge T$.  In the present work, we are interested in drifts $\bar{b}$ and costs $\bar{f}$ with sub-linear growth, hence possibly unbounded. 
Further details on the setting with all the technical assumptions will be given in Section \ref{sec:preliminaries}.\\
\indent The dynamics above is also motivated by economic models for corporate finance, systemic risk, and asset allocation. For instance, we can interpret players as firms whose values are represented by the state variables $X^{N,i}$ for $i\in \{1,\ldots,N\}$. Each company is affected by the fraction of both defaulted and non-defaulted firms and takes strategic decisions accordingly. Moreover, sub-linearity of the drift allows to include a mean-reversion term representing some herding behaviour. A possible application is the pricing of portfolio credit derivatives where the pricing depends upon the so called distance-to-default of the assets in the portfolio (\citet{hambly2017stochastic}). Alternatively, each player can be interpreted as a bank, whose monetary reserve evolves according to the stochastic dynamics in Eq.\eqref{eq:U:SDE} where the drift depends on both the rate of interbank borrowing/lending and on a controlled borrowing/lending rate to a central bank, as in \cite{carmona2015mean}. However, in \cite{carmona2015mean} no absorbing boundary conditions are considered. The latter features could be incorporated in the model by introducing absorbing boundary conditions at the default level, similarly to \cite{hambly2017stochastic}. This would enable to study the impact of defaults on systemic risk and stability of the financial system described by the game. Last but not least, the proposed set-up allows for a Brownian motion with an Ornstein–Uhlenbeck type drift modelling for the private state, a model that has been used (for instance) for the notion of flocking to default in the financial literature (\citet{fouque2013systemic}). However, in the present paper we focus on the mathematical properties of the proposed family of games and we leave the applications for future research.\medskip

\noindent\textit{Main results.}\quad The main contributions of the paper can be summarized as follows:
\begin{itemize}
\item We introduce the MFG with smooth dependence on past absorptions, i.e. the limit model corresponding to the above $N$-player games as $N$ tends to infinity. For a solution of the MFG, the empirical sub-probability measures  $(\mu_t^N)_{t \in [0,T]}$  are replaced by flows of sub-probability measures on $\R^d$; see Definition \ref{def:MFGsol}.
\item We prove existence of a relaxed feedback MFG solution and, under an additional convexity assumption, we show that there are optimal feedback strategies in strict form; see Theorem \ref{teo:existenceRelFeedSolMFG}, Proposition \ref{prop:existenceRelFeedlSolMFG} and Proposition \ref{prop:existenceStrictFeedlSolMFG}. Additionally, we show  that there exist  relaxed and strict feedback solutions that are Markovian up to the exit time; see Proposition \ref{lem:existencemarkovSolMFG}.
\item We prove uniqueness of the MFG solution under standard monotonicity conditions of the Lasry-Lions type formulated for sub-probability measures; see Theorem \ref{teo:uniqueness}.
\item We study approximate Nash equilibria for the $N$-player game in a setting where the dependence on the measure variable is finite-dimensional. Precisely, we show that if we have a feedback solution of the MFG (either relaxed or strict), we can construct a sequence of approximate Nash equilibria for the corresponding $N$-player games with a vanishing approximation error as $N \rightarrow \infty$; see Theorem \ref{teo:approxNashRel} and Corollary \ref{cor:approxNash}.
It is worth stressing that the construction produces approximate $N$-player equilibria in feedback strategies (instead of the more common open-loop strategies).
\end{itemize}

\noindent The proof of the existence of feedback solutions of the MFG is inspired by the truncation procedure introduced by \cite{lacker}. We construct a sequence of approximating MFGs, each one with bounded drift and cost functional, to which we can apply the results of \cite{campi2018n}. Then, we prove convergence of the solutions of these approximating MFGs to a solution of the original one. Nonetheless, the procedure in \cite{lacker} cannot be applied directly to our case mainly due to the history dependency and the discontinuities induced by past absorptions. In particular, a different instance of the mimicking result of \cite{brunick} applies to our framework.\\ 
\indent To establish the uniqueness result we follow standard monotonicity arguments, with some adjustments due to the dependence of the coefficients on a flow of sub-probability measures instead of probability measures. In particular, the uniqueness result relies on an additional (standard) monotonicity assumption on the running cost of the Lasry-Lions type.\\
\indent The proof of the construction of approximate Nash equilibria for the $N$-player game is based on weak convergence arguments and controlled martingale problems. The use of martingale problems in proving convergence to the  McKean-Vlasov limit and propagation of chaos for weakly interacting systems goes back to \cite{funaki1984certain}, \cite{oelschlager1984martingale} and \cite{meleard1996asymptotic}. We observe that, whereas standard results prove convergence in law of the empirical measures, in the present paper we follow the approach of \cite{lacker2018strong} to obtain a strong form of propagation of chaos with possibly unbounded and path-dependent drift. We show that the empirical measures converge in a stronger topology (the $\tau$-topology), a result that enables us to take the limit as $N \rightarrow \infty$ without assuming any regularity of the feedback strategies with respect to the state process.
In our framework, unlike  \cite{campi2018n}, the continuity of the MFG optimal control for almost every path of the state variable with respect of the Wiener measure is no longer feasible. Indeed, the PDE-based estimates that were used in  \cite{campi2018n} to get such a regularity are not available anymore due to the possible unboundedness of the drift and the running cost.
\medskip

\noindent\textit{Related literature.}\quad We have already discussed the paper \cite{campi2018n}, so here we focus on some other contributions in the literature of mean-field models and games related to our study. First, we cite the works of \cite{giesecke2013default} and \cite{giesecke2015large} where a model based on point processes for correlated defaults timing in a portfolio of firms is introduced and analysed.
\cite{giesecke2013default} prove a LLN for the default rate as the number $N$ of firms goes to infinity.\\
\indent Motivated by modelling the contagion effect are the works of \cite{ hambly2017stochastic}, \cite{hambly2019mckean} and \cite{hambly2019spde} too. The first work provides a LLN for the empirical measure of a system of finitely many (uncontrolled) diffusions on the half-line, absorbed when they hit zero and correlated through the proportion of absorbed processes. In \cite{hambly2019mckean} the model is extended to include a positive feedback mechanism when the particles hit the barrier, thus modelling contagious blow-ups. A mathematical complement to the previous work is provided in \cite{ledger2020uniqueness}. More recently,  \cite{hambly2019spde} have proposed a general model for systemic (or macroscopic) events. By working on a set-up similar to \cite{hambly2017stochastic}, they interpret the diffusions as distances-to-default of financial institutions and model the correlation effect through a common source of noise and a form of mean-reversion in the drift.  A form of endogenous contagion mechanism is also considered.\\
\indent On the side of applications to economics, \cite{chan2015bertrand} and \cite{chan2017fracking} study oligopolistic models with exhaustible resources formulated as MFGs with absorption at zero. Their model keeps track of the fraction of active players at each time. However, this fraction appears in the objective functions but not in the state variable.\\  
\indent Two more papers are those by \cite{delarue2015global} and  \cite{delarue2015particle}, where a particle system approach is used to study the mathematical properties of an integrate-and-fire model from neurology. The particles' dynamics have some resetting mechanism which activates as soon as some particle hits a given boundary. Besides, we cite two recent papers by \citet{nadtochiy2019particle,nadtochiy2020mean}. The first one focuses on the cascade effect in an interbank mean-field model with defaults and a contagion effect modelled via a singular interaction through hitting times. The second one investigates the associated mean-field game also including more general dynamics and connection structures.\\
\indent Finally, we mention a class of MFGs that has been considered quite recently especially in relation to bank run models, that is MFGs of optimal stopping or timing; see, for instance, \cite{bertucci2018optimal}, \cite{bouveret2020mean}, \cite{carmona2017mean} and \cite{nutz2018mean}. Therein, the agents solve an optimal stopping problem so that the terminal time is directly chosen by them instead of being determined by the evolution of the controlled state as in our setting. In both settings the terminal time is in fact a random time and the state evolution might be affected by the fraction of leavers and the empirical measure of the remainers. 
\medskip

\noindent\textit{Structure of the paper.}\quad In Section \ref{sec:preliminaries} we introduce the notation and present both the $N$-player and the MFGs along with the main assumptions. Section \ref{sec:mfg} contains the results on the existence of feedback MFG solutions. In Section \ref{sec:uniqueness} we prove the uniqueness of MFG solutions under some monotonicity condition of the Lasry-Lions type. In Section \ref{sec:Nplayer} we specialize to a finite dimensional setting and construct approximate Nash equilibria in feedback form for the $N$-player game using the MFG solutions. The technical results used in the paper can be found in the Appendix \ref{app:Appendix}.

\section{Preliminaries and assumptions}\label{sec:preliminaries}
\noindent In this section, we provide the definitions of the different spaces of trajectories and measures used in the paper along with the corresponding topologies, distances and notions of convergence.  
In addition, we describe the MFG with smooth dependence on past absorptions and give the definition of solution of the MFG.
We conclude the section by introducing the MFGs with truncated coefficients, which will be used in the proof of existence of MFG solutions.\medskip

\noindent\textit{Spaces of trajectories.}\quad Let $d\in\mathbb{N}$. We denote by $\OO\subset\Rd$ an open subset of $\Rd$ representing the space of the players' private states and by $\X\doteq C([0,T];\Rd)$ the space of $\Rd$-valued continuous trajectories on the time interval $[0,T]$, $T<\infty$. The space $\Rd$ is equipped with the standard Euclidean norm, always indicated by $|\cdot|$,  while $\X$ with the sup-norm, denoted by $\infnorm{\cdot}$, which makes $\X$ separable and complete. We use the notation $\| \cdot\|_{\infty,t}$ whenever the sup-norm is computed over the time interval $[0,t]$, $t < T$. Besides, we denote with $\X^N\doteq C([0,T];\R^{d\times N})$ the space of $N$-dimensional vectors of continuous trajectories and identify it with $\X^{\times N}$.\medskip

\noindent\textit{Spaces of measures.}\quad We use flows of probability and sub-probability measures to describe the distribution of players and its time evolution in $\OO$. For $E$ a Polish space, let $\M_f(E)$ denote  the space of finite Borel measures on $E$, $\PP(E)$ the space of Borel probability measures on $E$ and $\M_{\leq 1}(E)$ the space of Borel sub-probability measures on $E$, i.e. measures $\mu \in \M_f(E)$ such that $\mu(E)\leq 1$. These spaces are endowed with the weak convergence of measures (\citet{billingsley}). We will often write $\mu^n \weakconv \mu$ to indicate weak convergence of $\mu^n$ towards $\mu$ as $n \to \infty$ and $\xi_n\overset{\mathcal{L}}{\longrightarrow}\xi$ to denote convergence in law of a sequence of random variables $(\xi_n)_{n\in\N}$ (defined on possibly different probability spaces) to a limit random variable $\xi$.

We define by $\Upsilon_{\PP}^{T}(E)$ (resp. by $\Upsilon_{\leq 1}^{T}(E)$) the spaces of measurable flows of probability (resp. sub-probability) measures on $E$, i.e. the space of Borel measurable maps $\pi$ (resp. $\mu$) from the time interval $[0,T]$ to $\PP(E)$ (resp. $\M_{\leq 1}(E)$). Wherever possible without confusion, we use $\Upsilon_{\PP}^{T}$ (resp. $\Upsilon_{\leq 1}^{T}$) when $E=\Rd$. We denote by $\PP_1(E)$ and by $\M_{\leq 1,1}(E)$ the following subsets of $\PP(E)$ and $\mathcal{M}_{\le 1}(E)$:
\ba 
\PP_1\rbr{E}&\doteq &\cbr{\pi\in\PP\rbr{E}:\int_E d_E(x,x_0)\pi(dx) <\infty\text{ for some }x_0\in E},\nonumber\\
\M_{\leq 1,1}\rbr{E}&\doteq &\cbr{\mu\in\M_{\leq 1}\rbr{E}:\int_E d_E(x,x_0)\mu(dx) <\infty\text{ for some }x_0\in E}.
\nonumber
\label{eq:firstMom}
\ea
\noindent We endow $\PP_1(E)$ with the 1-Wasserstein distance $W_1$
\ba 
W_1(\mu,\nu)\doteq\inf_{\pi\in\Pi(\mu,\nu)}\int_{E\times E}d_E\rbr{x,y}d\pi(x,y)=\sup_{f\in\text{Lip}_1\rbr{E;\R}}\int_{E}f(x)d(\mu-\nu)(x)
\label{eq:was1Dist}
\ea
\noindent where $\Pi(\mu,\nu) \subset \PP_1(E\times E)$ represents the set of probability measures with given marginals $\mu$ and $\nu$, and $\text{Lip}_1(E;\R)$ the set of Lipschitz functions on $E$ with unitary Lipschitz constant. The second equality in Eq.\eqref{eq:was1Dist} is due to the Kantorovich-Rubinstein Theorem  (see, for instance, Theorem 6.1.1 in \citet{ambrosiogigli}). Notice that $(\PP_1(E),W_1)$  is a separable and complete metric space whenever $(E,d_E)$ is separable and complete. Finally, let $\Upsilon_{\PP,1}^{T}(E)$ (resp. $\Upsilon_{\leq 1,1}^{T}(E)$) denote the space of measurable flows of probability measures in $\PP_1(E)$ (resp. in $\M_{\leq 1,1}(E)$). Again,  wherever possible without confusion, we use  $\Upsilon_{\PP,1}^{T}$ and $\Upsilon_{\leq 1,1}^{T}$ when $E=\Rd$.\medskip

\noindent\textit{The canonical space.}\quad We will often work on the canonical filtered probability space, denoted by $(\Omega,\F,(\F_t)_{t\in[0,T]},\PB)$ and defined as follows. Set $\Omega\doteq\X$, let $\xi$ be an $\Rd$-valued random variable with law $\nu\in\PP(\Rd)$ and let $W$ be a $d$-dimensional Wiener process on $\X$ independent of $\xi$. Define $\wiener\in\PP(\X)$ as the law of $\xi+\sigma W$. Set $\F$ as the $\wiener$-completion of the Borel $\sigma$-algebra $\B(\X)$ and $(\F)_{t\in[0,T]}$ as the $\wiener$-augmentation of the filtration generated by the canonical process $\hat{X}$ on $\X$, i.e. $\hat{X}_t(\varphi)\doteq \varphi(t)$ for all $(t,\varphi)\in[0,T]\times\X$. In particular, $(\F)_{t\in[0,T]}$ satisfies the usual conditions. Finally set $\PB\doteq\wiener$ and $W\doteq\sigma^{-1}(\xi-\hat{X})$, which is a Wiener process on $\X$. Where no confusion is possible, we will write $X$ for $\hat{X}$.\medskip

\noindent Now, let $\mathcal{O} \subset \R^{d}$ be a non-empty open set, the set of non-absorbing states, and let $\Gamma\subset\R^{d}$ be the set of control actions.
For each $\varphi\in\X$ we set $\tau^{\varphi}\doteq\inf\lbrace t\in[0,T]:\,\varphi(t)\not\in\OO\rbrace$, with the convention $\inf\emptyset = \infty$, and $\tau(\varphi)\doteq \tau^{\varphi}\wedge T$. In order to set up the dynamics of the players' states, we need to introduce the following functions: 
\ba 
&&\bar{b}:\brackets{0,T} \times \Rd \times \mathcal{M}_{\leq 1,1}(\mathbb{R}^d)\times\Gamma \rightarrow \R^{d},\quad\quad\quad\quad\quad\,\,\, \sigma \in \R^{d \times d},\nonumber\\
&& \bar{f}:\brackets{0,T} \times \Rd \times  \mathcal{M}_{\leq 1,1}(\mathbb{R}^d)\times\Gamma \rightarrow [0,\infty),\quad\quad F:\brackets{0,T} \times  \Rd \rightarrow [0,\infty).
\nonumber
\ea
\noindent Since we will have to impose some joint continuity property for the functions above, in particular with respect to the $\mu$-variable, and there is no natural metrizable topology over the set of sub-probability measures $\mathcal M_{\le 1,1}(\R^d)$, it will be convenient to work with the following reparameterization of a suitable restriction of $\bar b$ and $\bar f$:
\begin{eqnarray}
b(t,\varphi,\theta,u)&\doteq &\bar{b}(t,\varphi(t),g(t,\theta),u), \quad \nonumber\\
f(t,\varphi,\theta,u)&\doteq &\bar{f}(t,\varphi(t),g(t,\theta),u)
\nonumber
\end{eqnarray}
where $b$ and $f$ are progressively measurable functionals such that
\ba 
&&b:\brackets{0,T} \times \X \times \PP_{ 1}(\X)\times\Gamma \rightarrow \R^{d} \nonumber,\\
&& f:\brackets{0,T} \times \X \times  \PP_{ 1}(\X)\times\Gamma \rightarrow [0,\infty)
\nonumber
\ea
while $g:[0,T]\times\mathcal{P}_1(\X)\rightarrow \mathcal{M}_{\leq 1,1}(\mathbb{R}^d)$ is defined by its action on the test functions of the 1-Wasserstein convergence, i.e., on the functions $\psi\in C(\mathbb{R}^d)$ with sub-linear growth, as
\begin{equation}
\int_{\mathbb{R}^d}\psi(x)g(t,\theta)(dx)\doteq \int_{\X}\psi(\varphi(t))\mathbf{1}_{[0,\tau^{\varphi})}(t)\theta(d\varphi).
\label{eq:defGtheta}
\end{equation}
In words, the functions $b$ and $f$ above are reparameterizatons of the \emph{restrictions} of $\bar b$ and $\bar f$, respectively, to the range of the map
\[ (t,\varphi , \theta , u) \mapsto (t, \varphi(t), g(t,\theta),u).\]
Moreover, for each $\mu\in\mathcal{M}_{\leq 1,1}(\mathbb{R}^d)$ and $\theta\in\mathcal{P}_1(\X)$ we introduce the notation
\begin{equation}
m(\mu)\doteq\int_{\mathbb{R}^d}|x|\mu(dx)\quad \text{and} \quad m(t;\theta)\doteq \int_{\X}|\varphi(t)|\mathbf{1}_{[0,\tau^{\varphi})}(t)\theta(d\varphi).
\nonumber
\end{equation}
\noindent Now, we collect the necessary assumptions on all initial data in order to state our main results. Some further assumptions will be given later in the paper when necessary.

\begin{itemize}
\item[\hypertarget{H1}{(\textrm{H1})}] The drift $\bar{b}$ satisfies the following uniform Lipschitz continuity:
\begin{eqnarray}
\left|\bar{b}(t,x,\mu,u)-\bar{b}(t,x',\mu,u)\right|\leq L|x-x'|,\quad x,x'\in\Rd
\nonumber
\label{eq:driftLip}
\end{eqnarray}
\noindent for any $(t,\mu,u)\in[0,T]\times\mathcal{M}_{\leq 1,1}(\mathbb{R}^d)\times\Gamma$. Moreover it has sub-linear growth, i.e.
\begin{eqnarray}
\left|\bar{b}(t,x,\mu,u)\right|\leq C\left(1+|x|+m(\mu)\right)
\nonumber
\label{eq:driftLin}
\end{eqnarray} 
\noindent for all $(t,x,\mu,u)\in[0,T] \times \Rd \times\mathcal{M}_{\leq 1,1}(\mathbb{R}^d)\times\Gamma$  and for a positive constant $C>0$.
\item[\hypertarget{H2}{(\textrm{H2})}] The running costs $\bar{f}$ and the terminal cost $F$ have sub-linear growth, i.e.
\begin{eqnarray}
&\bar{f}(t,x,\mu,u)&\leq C(1+|x|+m(\mu)),\nonumber \\
&F(t,x)&\leq  C(1+|x|),
\nonumber
\end{eqnarray}
\noindent for all $(t,x,\mu,u)\in[0,T]\times\Rd\times\mathcal{M}_{\leq 1,1}(\mathbb{R}^d)\times\Gamma$, $(t,x)\in[0,T]\times\mathbb{R}^d$ and for a positive constant $C>0$.
\item[\hypertarget{H3}{(\textrm{H3})}] $\bar{b}$ and $\bar{f}$ are such that their reparametrizations $b$ and $f$ are jointly continuous at points $(t,\varphi,\theta, u)\in\brackets{0,T} \times \X \times \PP_{ 1}(\X)\times\Gamma$ such that 
$\theta\ll\wiener$. Moreover, $F$ is jointly continuous on $[0,T]\times\Rd$.
\item[\hypertarget{H4}{(\textrm{H4})}] The set $\OO$ is open, convex and strictly included in $\Rd$ with $\C^2$-boundary, i.e. $\partial\OO$ is the graph of a $\C^2$ function. Alternatively, $\OO=(0, \infty)^{\times d}$ is also allowed.
\item[\hypertarget{H5}{(\textrm{H5})}] The set $\Gamma\subset\Rd$ is compact.
\item[\hypertarget{H6}{(\textrm{H6})}] The diffusion matrix $\sigma\in\R^{d\times d}$ has full rank.
\item[\hypertarget{H7}{(\textrm{H7})}] The initial distribution $\nu\in\PP(\Rd)$ has support in $\OO$ and satisfies $\int_{\OO}\e^{\lambda|x|^2}\nu(dx)<\infty$ for some $\lambda >0$.
\item[\hypertarget{H8}{(\textrm{H8})}] The initial conditions of the $N$-player game $X^{N,i}_0$, $i\in\{1,\ldots,N\}$, are i.i.d. and with the initial condition of the MFG $X_0$, they are all distributed as $\nu\in\PP(\Rd)$.
\end{itemize}
\noindent Before turning to the MFG dynamics, some remarks on the assumptions above are in order.

\brem The growth assumptions in \hyperlink{H1}{(\textrm{H1})} and  \hyperlink{H2}{(\textrm{H2})} could be further refined. For instance, one could assume sub-linear and sub-polynomial growth of the drift and diffusion matrix with suitable exponents as, e.g., in \cite{lacker}. Moreover, the running cost $f$ could certainly take real values; however, without loss of generality and given the interpretation as a cost term, we have assumed $f\geq 0$.\erem

\brem\label{rem:contReparam} The continuity properties in \hyperlink{H3}{(\textrm{H3})} are crucial in the passage to the limit performed in Proposition \ref{prop:approxMart}. Since the laws of the processes that we consider are absolutely continuous with respect to the Wiener measure $\wiener$ (they belong to the set $\Q\subset\PP(\X)$ of laws of Brownian-driven processes with sub-linear drift that we introduce and characterize in the Appendix \ref{app:Appendix}, cfr. Lemma \ref{lem:regularitySubLinear}), it is sufficient to require continuity at points $\theta\ll\wiener$. The passage to the limit in the measure argument can then be performed by Lemma \ref{lem:regularity} together with Lemma \ref{lem:regularityConvergence}. \erem

\brem Admittedly, compactness of $\Gamma$ is a strong assumption, but it will play an important role in order to obtain existence and uniqueness of weak solutions of the SDEs for the player state's dynamics in both the MFG and the $N$-player games. In particular, it enables a line of arguments based on Ben\v{e}s' condition -- ensured by the boundedness of the coefficient in the control variable -- and Girsanov's theorem (see Remark \ref{rem:existenceUniqWeak} for more precise references), which is one of the main tools of our approach.
\erem

\brem The nondegeneracy of $\sigma$  as in \hyperlink{H6}{(\textrm{H6})} is justified by the counter-example in \cite{campi2018n}, Section 7, where it was shown that a feedback MFG solution does not necessarily induce a sequence of approximate Nash equilibria with vanishing error. A careful inspection of such a counter-example reveals that it can be easily adapted to our setting since, in that particular context, dividing by the initial number of players $N$ (as in our setting) or renormalizing each time by the current number of players (as in the counter-example) turn out to be equivalent for $N$ large.  Finally, even though state dependency of the diffusion matrix can be handled using very similar techniques, we have decided to leave it out and focus on other more interesting aspects of the model. For the same reason we leave aside a possible dependence of $\sigma$ on the control, as it would just increase the level of technicality of the proofs due to the use of martingale measures (see \cite{lacker}).\erem
\medskip

\noindent\textit{The mean-field dynamics.}\quad Given a flow of sub-probability measures $\mu\in\Upsilon^T_{\leq 1,1}$ and a feedback progressively measurable control $u :\brackets{0,T}\times\X \rightarrow \Gamma$, the representative player's state evolves according to the equation
\ba 
X_t = X_0 + \int_{0}^{t}\bar{b}\round{s, X_s,\mu_s, u\round{s,X} }\,ds + \sigma W_t,\quad t \in \brackets{0,T},
\label{eq:MFG}
\ea
\noindent where $X$ is a $d$-dimensional stochastic process starting at $X_0 \overset{d}{\sim} \nu \in \PP(\Rd)$ and $W$ is a $d$-dimensional Wiener process on some filtered probability space $(\Omega, \F,(\F_t)_{t\in[0,T]},\PB)$. Solutions of Eq.\eqref{eq:MFG} are understood to be in the weak sense (see Remark \ref{rem:existenceUniqWeak} below).\\
\indent Let $\mathcal{U}_{fb}$ denote the set of all \textit{feedback controls} defined as
\ba
\mathcal{U}_{fb} \doteq \lbrace u:\brackets{0,T}\times\X \rightarrow \Gamma : \text{\,$u$ is progressively measurable}\rbrace.
\nonumber 
\ea
\noindent The cost associated with a strategy $u \in \mathcal{U}_{fb}$, a flow of sub-probability measures $\mu\in\Upsilon^T_{\leq 1,1}$ and an initial distribution $\nu \in\PP(\Rd)$ is given by (we omit, for the sake of simplicity, the explicit dependence on $\nu$)
\ba
J^{\mu}\round{u} \doteq \E\Biggl[\int_{0}^{\tau}\bar{f}\round{s, X_s,\mu_s, u\round{s,X} }ds+F\round{\tau,X_{\tau}}\Biggr]
\label{eq:MFGcost}
\ea
\noindent where $(\Omega, \F, (\F_t)_{t\in[0,T]},\PB, W, X)$ is a solution of Eq.\eqref{eq:MFG} under $u$ with initial distribution $\nu$, and $\tau \doteq \tau^{X} \wedge T$ the random time horizon. Finally we set
\[V^{\mu}\doteq \inf_{u \in \mathcal{U}_{fb}} J^{\mu}(u).\]
\brem\label{rem:existenceUniqWeak}
\noindent For a given flow of sub-probability measures $\mu$, thanks to the linear growth of $\bar{b}$ in the state variable $\varphi$ and to the boundedness of the action space $\Gamma$, we have that both existence and uniqueness in law of a weak solution of Eq.\eqref{eq:MFG} is guaranteed by Lemma \ref{lem:benevs}, and by Proposition 5.3.6, Remark 5.3.8 and Proposition 5.3.10 in \cite{karaztas} (see our Lemma \ref{lem:existenceUniquenessWeak}). Precisely, this can be proved by means of Girsanov's theorem and Ben\v{e}s' condition \citep{benevs}.
\erem
\noindent The notion of solution we consider for the MFG is the following.
\bdefi[\textit{Feedback MFG solution}]\label{def:MFGsol} A feedback solution of the MFG is a pair $(u,\mu) \in \mathcal{U}_{fb} \times \Upsilon_{\leq 1,1}^{T}$ such that:
\bi
\item[(i)] Strategy $u$ is optimal for $\mu$, i.e. $V^{\mu} = J^{\mu}(u)$.
\item[(ii)] Let $(\Omega,\F,(\F_t)_{t\in[0,T]},\PB,X,W)$ is a weak solution of Eq.\eqref{eq:MFG} with flow of sub-probability measures $\mu$, strategy $u$ and initial condition $\nu$. Then 
\[\mu_t(\cdot)=\PB(\{X_t \in \cdot\}\cap\{ \tau^X > t\}), \quad t\in[0,T].\]
\ei
\edefi

\noindent\textit{Relaxed controls.}\quad It will be very convenient to use relaxed controls (see \cite{karoui} for a precise definition), which allow us to view progressively measurable controls with values on a compact set $\Gamma$ as elements of the space of probability measures on $\Gamma$. The latter space is compact when endowed with the weak convergence of measures. The space $\V$ of relaxed controls is given by
\ba 
\V\doteq \cbr{q\in \M_f\round{[0,T]\times \Gamma}:q(dt,d\gamma)=dt q_t (d \gamma),\,t\mapsto q_t\in \PP(\Gamma)\,\text{Borel measurable}}
\nonumber
\ea
i.e. it is the set of all finite positive measures on $[0,T] \times \Gamma$ with Lebesgue time marginal. With a slight abuse of notation, we denote with $\hat{\Lambda}$ both the identity map and the canonical process on $\V$ (where no confusion is possible, we drop the hat and write $\Lambda$ in place of $\hat{\Lambda}$). Precisely, a single-player relaxed control is a $\mathcal{V}$-valued random variable $\Lambda$ such that $(\Lambda_t)_{t\in[0,T]}$ is a progressively measurable $\PP(\Gamma)$-valued stochastic process. We say that $\Lambda$ is a feedback control if there exists a progressively measurable functional $\lambda:[0,T]\times\X\rightarrow\PP(\X)$ such that $\Lambda_t = \lambda(t,X)$ for all $t\in[0,T]$, with $X$ denoting the player's dynamics. Moreover, we say that $\Lambda$ is a strict and feedback control if there exists $u\in\mathcal{U}_{fb}$ such that $\lambda(t,X) = \delta_{u(t,X)}$ for all $t\in[0,T]$.\\
\indent Let $\widetilde{\mathcal{U}}_{fb}$ be the set of relaxed feedback controls for the MFG.
We rewrite the dynamics and the cost functional of the MFG (Eq.\eqref{eq:MFG}) and Eq.\eqref{eq:MFGcost}) using relaxed controls:
\ba 
X_t &=& X_0 + \int_{[0,t]\times\Gamma}\bar{b}\round{s,X_s,\mu_s,u}\lambda\round{s,X}(du)ds+ \sigma W_t, \label{eq:MFGRel}\\
J^{\mu}\round{\lambda}&=&\E\Biggl[\int_{[0,\tau]\times\Gamma}\bar{f}\round{s,X_s,\mu_s,u}\lambda\round{s,X}(du)ds+F\round{\tau,X_{\tau}}\Biggr]
\label{eq:costMFGRel}
\nonumber 
\ea
\noindent where $t\in[0,T]$ and $\lambda\in\widetilde{\mathcal{U}}_{fb}$. Moreover, we extend accordingly the notion of feedback solutions of the MFG.
\bdefi[\textit{Relaxed feedback MFG solution}]\label{def:MFGsolRel} A relaxed feedback solution of the MFG is a pair $(\lambda,\mu) \in \tilde{\mathcal{U}}_{fb} \times \Upsilon_{\leq 1,1}^{T}$ such that:
\bi
\item[(i)] $\lambda$ is optimal, i.e. $V^{\mu} = J^{\mu}(\lambda)$.
\item[(ii)] Let $(\TO,\TF,(\TF_t)_{t\in[0,T]},\mathbb{Q},X,W)$ be a weak solution of Eq.\eqref{eq:MFGRel} with flow of sub-probability measures $\mu$, control $\lambda$ and initial condition $\nu$. Then 
\[\mu_t(\cdot)=\mathbb{Q}(\{X_t \in \cdot\}\cap\{ \tau^X > t\}), \quad t\in[0,T].\]
\ei
\edefi

\noindent\textit{Feedback and open-loop controls.}\quad
Feedback controls induce stochastic open-loop controls, i.e. tuples $(\Omega,\F,(\F_t)_{t\in[0,T]},\PB,X,u,W)$ that are weak solutions of
\ba 
X_t = X_0 + \int_{0}^{t}\bar{b}\round{s, X_s, \mu_s, u_s}\,ds + \sigma W_t,\quad t \in \brackets{0,T}
\label{eq:MFG2}
\ea
\noindent where $u$ is a progressively measurable $\Gamma$-valued stochastic process. As a consequence, the computation of the infimum of $J^{\mu}(\cdot)$ over the class of stochastic open-loop controls would imply a lower value for $V^{\mu}$. However, thanks to Proposition 2.6 in \cite{karoui}, the two minimization problems are equivalent from the point of view of the value function.\\
\noindent A similar argument holds also in the case of feedback relaxed controls, that induce relaxed stochastic open-loop controls, tuples  $(\TO,\TF,(\TF_t)_{t\in[0,T]},\mathbb{Q},X,\Lambda,W)$ that are weak solutions of
\ba 
X_t &=& X_0 + \int_{[0,t]\times\Gamma}\bar{b}\round{s,X_s,\mu_s,u}\Lambda_s(du)ds+ \sigma W_t,\quad t\in[0,T]
\label{eq:SDEopen}
\ea
\noindent where $\Lambda$ is a progressively measurable $\PP(\Gamma)$-valued stochastic process.\\
\noindent In the rest of the paper we will call $\mathbb{U}$ the set of open-loop controls and, for the sake of brevity and where no confusion is possible, denote with $u$ an element of $\mathbb{U}$ implying the whole tuple $(\Omega,\F,(\F_t)_{t\in[0,T]},\PB,X,u,W)$. Similarly,  we will call $\tilde{\mathbb{U}}$ the set of open-loop relaxed controls and denote with $\Lambda$ an element of $\tilde{\mathbb{U}}$ implying the whole tuple $(\TO,\TF,(\TF_t)_{t\in[0,T]},\mathbb{Q},X,\Lambda,W)$.\\
\medskip

\noindent\textit{The extended canonical probability space.}\quad When dealing with relaxed controls we will work on the following extension of the canonical probability space $\X$. Set $\TO\doteq \X\times\V$, let $\F$ and $(\F_t)_{t\in[0,T]}$ be the canonical $\sigma$-algebra and the canonical filtration on $\X$, respectively, whereas $\G$ and $(\G_t)_{t\in[0,T]}$ denote the Borel $\sigma$-algebra and the filtration generated by the canonical process $\hat{\Lambda}$ on $\V$, respectively. Finally, we set $\TF_t\doteq \F_t\otimes\G_t$ for all $t\in[0,T]$, and $\TF\doteq \F\otimes \G$.\\

\noindent\textit{Approximating MFGs.}\quad We conclude this preliminary section by introducing a suitable sequence of approximating MFGs, which is obtained by truncation of the coefficients of the original MFG similarly as in \cite{lacker}. Such a sequence will be useful in the proof of existence of a MFG solution along the following lines: we will prove existence of feedback MFG solutions of the approximating MFGs in the sequence by extending the existence result of \cite{campi2018n}. Then, by letting the truncation threshold go to infinity, we will obtain a solution of the original MFG. This approach relies on two additional assumptions (Assumptions \hyperlink{C1}{(\textrm{C1})} and \hyperlink{C2}{(\textrm{C2})} below) that will be introduced later in this part.\\
\indent Let $(K_n)_{n\in\mathbb{N}}\subset\R_+$ be an increasing sequence such that $K_n\nearrow+\infty$. The $n^{\rm th}$ approximating MFG model, denoted by  MFG($n$), is obtained as follows.
\bi
\item[$\round{\mathbf{T}_n}$] $\bar{b}^n(x)=\bar{b}(x)$ when $|\bar{b}(x)|\leq K_n$, while it is continuously truncated at level $K_n$, i.e. $|\bar{b}^n(x)|= K_n$, otherwise. Similarly for the costs $\bar{f}^n$ and $F^n$ and for the associated functions $b^n$ and $f^n$.
\ei
\noindent Notice that we do not truncate the possibly unbounded set $\mathcal O$ of non-absorbing states. In each MFG($n$) the representative player's state evolves as in Eq.\eqref{eq:MFG} with $\bar b$ replaced by $\bar b^n$, i.e.
\ba 
X_t = X_0 + \int_{0}^{t}\bar{b}^n\round{s, X_s, \mu_s,u(s,X)}\,ds + \sigma W_t , \quad t\in [0,T]
\label{eq:MFGn} 
\ea
when the player is using the strict control $u$, and similarly when he/she is using a relaxed control. Moreover, in the cost functional $\bar{f}$ and $F$ are replaced by their truncated counterpart $\bar f^n$ and $F^n$. The associated cost functional is denoted by $J^{n,\mu}\round{u}$ or $J^{n,\mu}\round{\lambda}$ depending on whether the player is implementing a strict strategy $u$ or a relaxed one $\lambda$. The optimal values are defined, accordingly, by 
\[V^{n,\mu}\doteq \inf_{u \in \mathcal{U}_{fb}} J^{n,\mu}(u).\] 
The definitions of strict and relaxed MFG solutions given above for the (un-truncated) MFG can clearly be applied to the approximating MFG($n$)s with the obvious modifications. We associate to the MFG($n$)s the following Hamiltonians: 
\ba
h^n(t,x,\theta ,z,u)&\doteq& f^n(t,x,\theta ,u) + z\,\sigma^{-1}\,b^n(t,x,\theta ,u),\nonumber\\
H^n(t,x,\theta ,z)&\doteq&\inf_{\substack{u\in\Gamma}} h^n(t,x,\theta ,z,u)
\label{eq:hamiltonian}
\nonumber 
\ea
\noindent and the set of minimizers
\ba 
A^n(t,x,\theta ,z) &\doteq &\left\{u\in\Gamma\,:\,h^n(t,x,\theta ,z,u)=H^n(t,x,\theta ,z)\right\}
\label{eq:minimizers} 
\nonumber
\ea
\noindent for $(t,x,\theta ,z)\in[0,T]\times\Rd\times\PP_1(\X)\times\R^{d}$. In the next section on existence of MFG solutions we will rely on the following additional convexity assumptions:
\bi 
\item[\hypertarget{C1}{(\textrm{C1})}] For each $n\in\N$, $A^n(t,x,\theta,z)$ is convex for all $(t,x,\theta,z)\in[0,T]\times\Rd\times \PP_1(\X)\times\R^{d}$.
\item[\hypertarget{C2}{(\textrm{C2})}] The running cost $f$ is convex in the control variable $u\in\Gamma$.
\ei
\brem
\noindent Assumption \hyperlink{C1}{(\textrm{C1})}  is common in control theory and it is crucial in order to apply fixed point theorems. In our case it is satisfied if, for instance, the running cost $f$ is bounded and convex in the control variable $u\in\Gamma$. Indeed in this case, due to the flexibility in the choice of the truncation thresholds, choosing $K^n\geq\infnorm{f}$ for all $n\in\mathbb{N}$ we have $f^n=f$ for all $n\in\mathbb{N}$. Then convexity is preserved by adding any sub-linear term. Finally, we observe that Assumption \hyperlink{C2}{(\textrm{C2})} will be used in Section \ref{sec:existenceSol} for obtaining the existence of strict MFG solutions. 
\erem

\section{Existence of solutions of the mean-field game}\label{sec:mfg}
\noindent Throughout this section Assumptions \hyperlink{H1}{(\textrm{H1})}-\hyperlink{H8}{(\textrm{H8})} are in force. Under these and the additional convexity Assumptions \hyperlink{C1}{(\textrm{C1})} and \hyperlink{C2}{(\textrm{C2})} we show that both a relaxed and a strict feedback solution of the MFG exist; see Theorem \ref{teo:existenceRelFeedSolMFG} below together with Proposition \ref{prop:existenceRelFeedlSolMFG} and Proposition \ref{prop:existenceStrictFeedlSolMFG}. In addition, we guarantee the existence of a feedback solution of the MFG with Markovian feedback strategy up to the exit time; see Proposition \ref{lem:existencemarkovSolMFG}. Our main existence result can be stated as follows. 

\bteo[\textit{Existence of relaxed and strict feedback MFG solutions}]\label{teo:existenceRelFeedSolMFG} Under Assumptions  \hyperlink{H1}{(\textrm{H1})}-\hyperlink{H8}{(\textrm{H8})} and \hyperlink{C1}{(\textrm{C1})}, there exists a relaxed feedback MFG solution $(\lambda,\mu)$. Moreover, under the additional Assumption \hyperlink{C2}{(\textrm{C2})} , there exists a strict feedback MFG solution $(u,\mu)$.
\eteo
\noindent To prove Theorem \ref{teo:existenceRelFeedSolMFG}, we proceed by approximation in the sense that, first, we prove that each MFG($n$) introduced in the previous section has a feedback (strict) solution by extending the results in \cite{campi2018n}; see Subsection \ref{subsec:approximatingMFG}. Then, we prove the convergence of such approximating solutions to a feedback (relaxed) solution of the original MFG by passing to the limit with the truncation thresholds; see Subsection \ref{subsec:convergence}. \\

\noindent Before proceeding, we ensure the well-posedness of the game in the sense that we show that the private state $X$ of the representative agent remains in $\OO$ up to time $T$ with some positive probability. This is the content of the following lemma. 

\blem \label{lem:exitO} Grant Assumptions  \hyperlink{H1}{(\textrm{H1})}-\hyperlink{H8}{(\textrm{H8})}. Let $(\Omega,\F,(\F_t)_{t\in[0,T]},\PB,X,W)$ be a weak solution of Eq.\eqref{eq:MFG}. Then $\PB(\tau^X>t)>0$ for all $t\in[0,T]$. 
\elem
\begin{proof}\noindent Set $b_t\doteq \bar{b}(t,X_t,\mu_t,u(t,X))$ for $t\in[0,T]$, and define $Z \doteq (Z_t)_{t\in[0,T]}$ as
\ba
Z_t\doteq\mathcal{E}_t\round{-\int_0^{\cdot}\sigma^{-1}b_sdW_s},\quad t\in[0,T],\nonumber
\ea
where $\mathcal E_t (\cdot)$ denotes the Dol\'eans-Dade stochastic exponential. By Lemma \ref{lem:benevs}, $Z$ is a true martingale. Define $\QB$ by $\frac{d\QB}{d\PB}\doteq Z_T$. By Girsanov's theorem $\widetilde{W}_t\doteq W_t+\int_0^t\sigma^{-1}b_sds$, $t\in[0,T]$, is a $\QB$-Wiener process, and under $\QB$ the process $X$ has law $\wiener$. As a consequence of the law of iterated logarithms, any Wiener process remains in an open set, hence in $\OO\subset \Rd$, for a finite time with strictly positive probability. Therefore $\QB(\tau^X>T)>0$ and thus $\PB(\tau^X>T)>0$.
\end{proof}

\subsection{Approximating MFGs}\label{subsec:approximatingMFG}
\noindent In this subsection we prove existence of solutions of the approximating MFG($n$)s.

\bteo[\textit{Existence of solutions of MFG($n$)}]\label{teo:existenceSolMfgn} Let $n\in\mathbb{N}$. Under Assumptions  \hyperlink{H1}{(\textrm{H1})}-\hyperlink{H8}{(\textrm{H8})} and  \hyperlink{C1}{(\textrm{C1})} there exists a feedback solution $(u^n,\mu^n)$ of MFG($n$).
\eteo
\begin{proof} The proof follows similar steps to those in Section 6 of \cite{campi2018n}: we only sketch here the main steps. The main difference with \cite{campi2018n} is that, due to Assumption \hyperlink{C1}{(\textrm{C1})}, we have to deal with set-valued maps, hence to apply a version of Kakutani's fixed point theorem instead of Brouwer's. We use the version proposed by \cite{carmonaweak}, Proposition 7.4, which is in turn based on the results of \cite{cellina}. Other adjustments are due to  the fact that $\mu$ is a flow of sub-probability measures (instead of probability measures) and that $\OO$ can be unbounded.
\\
\indent Fix $n\in\mathbb{N}$. The proof is based on the construction of a suitable map $\Psi:\PP(\X)\times\mathbb{U}\rightarrow \PP(\X)$ on an appropriate compact and convex subset of $\PP(\X)$, where $\mathbb{U}$ is the space of progressively measurable $\Gamma$-valued stochastic processes. The fixed points of $\Psi$ will provide MFG($n$) solutions. More in detail, define $\mathcal{Q}_{\nu,K}$ as the set of laws $\theta\in\PP(\X)$ of any process of the type
\ba 
\xi+\int_0^tb_sds+\sigma W_t,\quad t\in[0,T]
\nonumber
\ea
defined on some filtered probability space with a Wiener process $W$, $\xi\simdistr \nu$, drift $(b_t)_{t\in[0,T]}$ adapted and bounded by $K>0$. Let us consider
\[ \Psi : \Q_{\nu,K_n}\times\mathbb{U} \ni (\theta,u) \mapsto \mathbb{P}^{\theta,u} \circ X^{-1} \in \Q_{\nu,K_n},\]
where $X$ is the canonical process on $\X$ and the probability measure $\mathbb{P}^{\theta,u}$ is defined as follows. Let $(\theta,u)\in\Q_{\nu,K_n}\times\mathbb{U}$ and let $\mu^{\theta}\in\Upsilon_{\leq 1}^{T}$ be defined as $\mu_t^{\theta}(\cdot)\doteq\theta(\{X_t\in\cdot \}\cap\{ \tau^{X}>t\})$ for all $t\in[0,T]$. Let $(\Omega,\F^{u},(\F^{u}_t)_{t\in[0,T]},\mathbb{P}^{\theta,u},X,W^{u})$ be the weak solution of
\ba 
X_t = X_0 + \int_{0}^{t}\bar{b}^n(s, X_s, \mu_s^{\theta},u_s)\,ds + \sigma W^{u}_t,\quad t\in[0,T]
\nonumber 
\ea
\noindent on the canonical space $(\Omega\doteq \X,\F,(\F_t)_{t\in[0,T]},\PB)$. Moreover, for $\theta\in \Q_{\nu,K_n}$ we call $u^{\theta}$ an optimal control for the cost
\ba
J^{n,\mu^{\theta}}\round{u} \doteq \E^{\mathbb{P}^{\theta,u}}\Biggl[\int_{0}^{\tau}\bar{f}^n(s, X_s, \mu_s^{\theta},u_s)ds+F^n\round{\tau,X_{\tau}}\Biggr].
\nonumber
\ea
\noindent Such optimal controls $u^{\theta}$ can be constructed by standard BSDE techniques as in \cite{campi2018n}, Section 6.1, by means of \cite{darling}, Theorem 3.4, due to the random terminal times. Under Assumption \hyperlink{C1}{(\textrm{C1})} optimal controls $u^{\theta}$ are in general not unique. Indeed
\ba 
A^{n}(\theta)\doteq \cbr{u^\theta \in\mathbb{U}: u^\theta \in A^{n}(\cdot,X_{\cdot},\theta,Z^{\theta}_{\cdot}),\,\mathcal{L}_T\otimes\PB-a.e.}
\nonumber
\ea
\noindent provides an entire set of optimal controls, where $Z^{\theta}$ is part of the the solution of the associated adjoint BSDE and $\mathcal{L}_T$ denotes the Lebesgue measure on $[0,T]$. Moreover, by measurable selection there exists a measurable function $\hat{u}^{n,\theta}:[0,T]\times\R^d \times\mathcal{Q}_{\nu,K_n}\times\R^d\rightarrow \Gamma$ such that
\ba 
\hat{u}^{n,\theta}(\cdot ,X_\cdot ,\theta,Z^\theta _\cdot )\in A^n(\theta), \quad \mathcal{L}_T\otimes\PB-\textrm{a.e.}
\nonumber
\ea
\noindent Additionally, $\hat{u}^{n,\theta}(t,X_t,\theta,Z^{\theta}_t)$, for $t\in[0,T]$, is a progressively measurable control process that can be written in feedback form. Indeed, since $Z^{\theta}$ is progressively measurable for the canonical filtration, it can expressed as $Z^{\theta}_t=\zeta^{\theta}(t,X)$ for some progressively measurable functional $\zeta^{\theta}:[0,T]\times\X\rightarrow\R^d$ and for any $t\in[0,T]$.\\
\indent Now, a fixed point for the map $\Psi$ is a probability measure $\theta\in \mathcal{Q}_{\nu,K_n}$ such that $\theta\in\Psi(\theta,A(\theta))$. Existence is provided by Proposition 7.4 in \cite{carmonaweak}, so to conclude the proof it suffices to check that all the required assumptions are satisfied in our case.
The set $\mathcal{Q}_{\nu,K_n}\subset\PP(\X)$ is a (weakly) compact, convex and metrizable subset of $C_b^*(\X)$, the dual of the space of bounded and continuous functions on $\X$, which is a locally convex topological vector space with the weak* topology (that induces the weak convergence of measures on $\PP(\X)$). We endow the vector space $\mathbb{U}$ with the norm $\nor{\cdot}_{\mathbb{U}}$ defined as $\nor{u}_{\mathbb{U}}\doteq\E[\int_0^T|u_t|dt]$. As a consequence of Berge's maximum theorem \citep[Theorem 17.31]{aliprantis} and of Assumption \hyperlink{C1}{(\textrm{C1})} the set-valued map $A^n:\mathcal{Q}_{\nu,K_n}\rightarrow\mathbb{U}$ is upper hemicontinuous and has non-empty convex and closed values (see the proof of Lemma 7.11 in \cite{carmonaweak}). Therefore, Proposition 7.4 in \cite{carmonaweak} applies, yielding the existence of a feedback solution of MFG($n$).
\end{proof}

\noindent\textit{A-priori estimates.}\quad Here, we show that the moments up to any order $\alpha\geq 1$ of the state process remain bounded uniformly in $n$. Such estimates will be very useful when we will relax the truncation in the next section.
\blem[\textit{A-priori estimates}]\label{lem:alphaEstimate} Grant Assumptions \hyperlink{H1}{(\textrm{H1})}-\hyperlink{H8}{(\textrm{H8})} and \hyperlink{C1}{(\textrm{C1})}. Consider feedback solutions $(u^n,\mu^n)_{n\in\N}$ and $(u,\mu)$ of the MFG(n)'s and of the MFG, respectively. 
Let $(\Omega^n,\F^n,(\F^n_t)_{t\in[0,T]},\PB^n,X^n,W^n)_{n\in\N}$ be a sequence of weak solutions of the SDEs in Eq.\eqref{eq:MFGn} and $(\Omega,\F,(\F_t)_{t\in[0,T]},\PB,X,W)$ a weak solution of the SDE in Eq.\eqref{eq:MFG}. Then for any $\alpha\geq 1$
\ba 
\sup_{n\in\N}\E^{\PB^n}\brackets{\infnorm{X^n}^{\alpha}}\leq K(\alpha)\quad\text{and}\quad\E^{\PB}\brackets{\infnorm{X}^{\alpha}}\leq K(\alpha)
\nonumber
\label{eq:alphaEstimate}
\ea 
where $K(\alpha)<\infty$ is a positive constant independent of $n$.
\elem
\begin{proof}
\noindent This follows from standard estimates that rely on the drift's sub-linear growth and on Gr\"onwall's lemma.
\end{proof}

\subsection{Convergence of the approximating MFGs}\label{subsec:convergence}
\noindent Let $(u^n,\mu^n)_{n\in\N}$ be a sequence of feedback solutions of the approximating MFGs introduced in the previous Subsection \ref{subsec:approximatingMFG}, whose existence is guaranteed by Theorem \ref{teo:existenceSolMfgn}.  In addition, let $(\Omega^n,\F^n,(\F^n_t)_{t\in[0,T]},\PB^n,X^n,W^n)_{n\in\N}$ be a sequence of weak solutions of the SDEs in Eq.\eqref{eq:MFGn} associated to $(u^n,\mu^n)_{n\in\N}$. Let $\theta^n$ be defined as $\theta^n \doteq \PB^n\circ (X^n)^{-1}$ for each $n\in\N$. 

To prove the convergence of the approximating MFGs we proceed in the following way. First, we show that there exists a subsequence of $(\theta^n)_{n \in \mathbb{N}}$, say $(\theta^{n_k})_{n_k \in \mathbb{N}}$, that converges in  $\PP_1(\X)$ to some limit $\theta\in\PP_1(\X)$. To prove this, we interpret $(u^n,\mu^n)_{n\in\N}$ as relaxed feedback solutions, $(\lambda^n,\mu^n)_{n\in\N}$. Second, we show that also the sequence of the corresponding extended laws $(\Theta^n)_{n\in\N}\subset\PP(\X\times\V)$  converges in $\PP_1(\X\times\V)$ to some limit $\Theta\in\PP_1(\X\times\V)$.
Finally, we characterize the limit points by means of the martingale problem of Stroock and Varadhan (see \citet{stroock1969,stroock}).
\blem[\textit{Relative compactness}]\label{lem:tightness}$(\theta^n)_{n\in\N}$ is relatively compact in $\PP(\X)$.
\elem
\begin{proof}First, we prove tightness by applying Aldous' criterion (see, e.g., \cite{jacod}, Condition VI.4.4), that is
\ba 
\lim_{\delta\rightarrow 0}\underset{n\rightarrow\infty}{\lim\sup}\sup_{\tau\leq\sigma\leq\tau+\delta}\PB^n\round{\modl{X^n_{\sigma}-X^n_{\tau}}\geq r}=0
\nonumber
\ea
\noindent for all $r>0$ and where $\tau$ and $\sigma$ are stopping times bounded by $T$. Indeed, we have
\ba 
\PB^n\round{\modl{X^n_{\sigma}-X^n_{\tau}}\geq r}\leq \frac{\E^{\PB^n}\brackets{\modl{X^n_{\sigma}-X^n_{\tau}}}}{r}
\nonumber
\ea
\noindent and
\ba 
\E^{\PB^n}\brackets{\modl{X^n_{\sigma}-X^n_{\tau}}}&\leq &\E^{\PB^n}\brackets{\int_{\tau}^{(\tau+\delta)\wedge T}\modl{\bar{b}^n(t,X^n_t,\mu_t^n,u^n(t,X^n))}dt}+ |\sigma|((\tau+\delta)\wedge T-\tau)^{\frac{1}{2}}C^W_T\nonumber\\
&\leq &\E^{\PB^n}\brackets{C\int_{\tau}^{(\tau+\delta)\wedge T}(1+\infnormt{X^n}+\sup_{n\in\N}\E^{\PB^n}{\infnormt{X^n}}+|u^n(t,X^n)|)dt}\nonumber\\
&&+ |\sigma|((\tau+\delta)\wedge T-\tau)^{\frac{1}{2}}C^W_T\nonumber\\
&\leq &\E^{\PB^n}\brackets{C\int_{\tau}^{(\tau+\delta)\wedge T}(1+\infnorm{X^n}+K+|u^n(t,X^n)|)dt}\nonumber\\
&&+ |\sigma|((\tau+\delta)\wedge T-\tau)^{\frac{1}{2}}C^W_T
\nonumber
\ea
\noindent for some constants $C^W_T,K>0$ independent of $n\in\N$. Then we conclude by Lemma \ref{lem:alphaEstimate}. Relative compactness then follows from Prohorov's Theorem.
\end{proof}
\noindent Now, let $\theta\in\PP(\X)$ be a limit point for $(\theta^n)_{n\in\N}$ and let $(\theta^{n_k})_{n_k\in\N}$ be a subsequence of $(\theta^n)_{n\in\N}$ such that 
$\theta^{n_k}\weakconv\theta$ as $n_k\rightarrow\infty$. With a slight abuse of notation, in what follows we identify $(\theta^{n_k})_{n_k\in\N}$ with $(\theta^n)_{n\in\N}$. We now show that the latter convergence is actually stronger by proving that $(\theta^n)_{n\in\N}$ converges to $\theta$ in the 1-Wasserstein distance. 
 
\blem[\textit{Convergence in the 1-Wasserstein distance}]\label{lem:conv1Wass}Let $(\theta^n)_{n\in\N}$  be as above. Then
$W_1(\theta^n,\theta)\rightarrow 0$ and $\theta\in\PP_1(\X)$.  
\elem
\begin{proof}\noindent Notice that by Lemma \ref{lem:alphaEstimate} we have $(\theta^n)_{n\in\N}\subset\PP_1(\X)$. To prove convergence in the 1-Wasserstein distance, we have to show that (see, for instance, Theorem 7.12.ii in \citet{villani})
\ba 
\lim_{R\rightarrow\infty}\sup_{n\in\mathbb{N}}\E^{\PB^n}\brackets{\infnorm{X^{n}}\mathbf{1}_{\cbr{\infnorm{X^{n}}\geq R}}}=0.
\nonumber
\ea
\noindent Set $\alpha,\beta>1$ such that $\frac{1}{\alpha}+\frac{1}{\beta}=1$. Then, for any  $\epsilon>0$ by Young's and Markov's inequalities, and by Lemma  \ref{lem:alphaEstimate} we have
\ba 
\E^{\PB^n}\brackets{\infnorm{X^{n}}\mathbf{1}_{\cbr{\infnorm{X^{n}}\geq R}}}
&\leq &\epsilon^{\alpha}\frac{\E^{\PB^n}\brackets{\infnorm{X^{n}}^{\alpha}}}{\alpha}+\dfrac{\PB^n(\infnorm{X^{n}}\geq R)}{\epsilon^{\beta}\beta}\nonumber\\
&\leq &\epsilon^{\alpha}\frac{K(\alpha)}{\alpha}+\dfrac{K}{\epsilon^{\beta}\beta R}
\nonumber
\ea
\noindent for some positive constants $K(\alpha)$ and $K$ independent of $n\in\N$.  The conclusion immediately follows thanks to the fact that convergence in the 1-Wasserstein distance preserves the finiteness of the first moment.
\end{proof}

\bprop[\textit{Absolute continuity of limit measures}]\label{prop:absContLimit} Let $\theta,(\theta^n)_{n\in\N}\subset\PP_1(\X)$ be as in Lemma \ref{lem:conv1Wass}. Then $\theta\ll\wiener$, i.e. $\theta$ is absolutely continuous with respect to $\wiener$.
\eprop
\begin{proof} By construction $\theta^n\ll\wiener$ for all $n\in\N$, hence we have to make sure that the absolute continuity is also preserved in the limit. For doing so, we apply Theorem X.3.3 in \cite{jacod}. In particular, we have to verify that all assumptions therein are fulfilled, which in our setting are reduced to the following properties:
\bi
\item[(i)] The contiguity of the sequence of $\theta^n$ with respect to the Wiener measure $\wiener$, i.e. for any sequence of measurable sets $B_n$ with $\wiener(B_n) \to 0$ we have $\theta^n (B_n) \to 0$ as $n \to \infty$ (see, e.g., Definition V.1.1 in \citet{jacod}).
\item[(ii)] The tightness of the sequence of $\wiener$-martingales $(M^n)_{n\in\N}$, where each $M^n=(M^n_t)_{t\in[0,T]}$ is defined as
\ba 
M^n_t\doteq\mathcal{E}_t\round{\int_0^{\cdot}\sigma^{-1}\bar{b}^n(s,X_s,\mu^n_s,u^n(s,X))dW_s},\quad t\in[0,T].
\nonumber
\ea
\ei
\noindent In order to check property (i), we first show that the sequence of Radon-Nikodym derivatives $(\frac{d\theta^n}{d\wiener})_{n\in\mathbb{N}}$ is uniformly integrable under $\wiener$. This is a consequence of the following bound:
\ba
\sup_{n\in\N}\E^{\wiener}\left[\left(\frac{d\theta^n}{d\wiener}\right)^p\right]<\infty, \quad p\in[1,\infty) \label{supRNder}
\ea
which follows from Corollary \ref{cor:expmartMoments} and by fact that, by inspection of the proofs of Lemma \ref{lem:benevs} and Corollary \ref{cor:expmartMoments}, all bounds are uniform in $n\in\N$.\\
\indent Now, property (i) can be obtained as follows: for all sequences of measurable sets $B_n$ with $\wiener(B_n) \to 0$, we have
\[ \theta^n (B_n) = \mathbb E^{\wiener} \left[  \frac{d\theta^n}{d\wiener} \mathbf 1_{B_n}\right] \to 0, \quad n \to \infty,\]
by an application of dominated convergence theorem due to the bound in Eq.\eqref{supRNder}. Hence the sequence of measures $\theta^n$ is contiguous to $\wiener$.\\
\indent Property (ii) follows from Aldous criterion \cite[Condition VI.4.4]{jacod}, that is
\ba 
\lim_{\delta\rightarrow 0}\underset{n\rightarrow\infty}{\lim\sup}\sup_{\tau\leq\sigma\leq\tau+\delta}\wiener\round{\modl{M^n_{\sigma}-M^n_{\tau}}\geq r}=0
\label{eq:jacod}
\ea
\noindent for all $r>0$ and where $\tau$ and $\sigma$ are stopping times bounded by $T$. As a consequence, we will also have the tightness property for the pair $(X,M^n)_{n\in\mathbb{N}}$ under the measure $\wiener$. By Theorem VI.4.13 in \cite{jacod} it is sufficient to check the tightness property for the corresponding quadratic variation processes
\ba 
\langle M^n\rangle_t=\int_0^t \modl{\sigma^{-1}\bar{b}^n(s,X_s,\mu^n_s,u^n(s,X))M^n_s}^2ds,\quad t\in[0,T].
\nonumber
\ea
\noindent First, by Markov's inequality $\wiener(|\langle M^n\rangle_{\sigma}-\langle M^n\rangle_{\tau}|\geq r)\leq\frac{1}{r}\E^{\wiener}[|M^n_{\sigma}-M^n_{\tau}|]$. Then,  by Young's inequality for all $p,q>1$ such that $\frac{1}{p}+\frac{1}{q}=1$ we have
\ba 
\mathbb{E}^{\wiener}\brackets{\modl{\langle M^n\rangle_{\sigma}-\langle M^n\rangle_{\tau}}} &\leq & \mathbb{E}^{\wiener}\brackets{\int_{\tau}^{(\tau+\delta)\wedge T} \modl{\sigma^{-1}}^2\modl{\bar{b}^n(s,X_s,\mu^n_s,u^n(s,X))}^2\modl{M^n_s}^2ds}\nonumber\\
&\leq &\frac{1}{p}\modl{\sigma^{-1}}^2\int_{\tau}^{(\tau+\delta)\wedge T} \mathbb{E}^{\wiener}\brackets{\modl{\bar{b}^n(s,X_s,\mu^n_s,u^n(s,X))}^{2p}}ds\nonumber\\
&& +\frac{1}{q}\modl{\sigma^{-1}}^2\int_{\tau}^{(\tau+\delta)\wedge T} \mathbb{E}^{\wiener}\brackets{\modl{M^n_s}^{2q}}ds\nonumber\\
&\leq & \round{\frac{K(p)}{p}+\frac{K(q)}{q}}\modl{\sigma^{-1}}^2\round{(\tau+\delta)\wedge T-\tau}
\nonumber
\ea
\noindent for some positive constants $K(p)$ and $K(q)>0$ independent of $n \in \mathbb{N}$. Notice that the last inequality is a consequence of Lemma \ref{lem:alphaEstimate} and Property (i). Therefore, Aldous' criterion in Eq.\eqref{eq:jacod} is satisfied.\medskip

\noindent After checking properties (i) and (ii) above, we can at last apply Theorem X.3.3 in \cite{jacod}, yielding that the tightness of $(\wiener\circ(X,M^n)^{-1})_{n\in\mathbb{N}}$ implies the tightness of $(\theta^n\circ(X,M^n)^{-1})_{n\in\mathbb{N}}$. In particular, if $(\wiener\circ(X, M^n)^{-1})_{n\in\mathbb{N}}$ weakly converges to some $\Theta'$ in $\PP(\mathcal{X}\times \mathcal X)$ then $(\theta^n\circ(X,M^n)^{-1})_{n\in\mathbb{N}}$ weakly converges to some other $\Theta''\ll\Theta'$ in $\PP(\mathcal{X}\times \mathcal X)$, and the same holds true for their first marginals on $\mathcal{X}$. Therefore, we can conclude that $\theta\ll\wiener$.
\end{proof}

\noindent\textit{Compactification method.}\quad So far we have established the convergence of the laws $(\theta^n)_{n\in\N}$ to some limit law $\theta$ in the 1-Wasserstein distance. Now, in order to prove the convergence of the approximating feedback solutions $(u^n,\mu^n)_{n\in\N}$ to some feedback MFG solution $(u,\mu)$, we need to show that the sequence of optimal controls $(u^n)_{n\in\N}$ converges to a control $u$, which is optimal for the limit game.\\
\indent To do this, we interpret the sequence of strict        
feedback solutions $(u^n,\mu^n)_{n\in\N}$ as a sequence of relaxed feedback solutions  $(\lambda^n,\mu^n)_{n\in\N}$, by defining $\lambda^n:[0,T]\times\X\rightarrow\PP(\Gamma)$ as $\lambda^n(t,\varphi)\doteq\delta_{u^n(t,\varphi)}$ for all $(t,\varphi)\in[0,T]\times\X$ and for all $n\in\N$. Furthermore, we identify each $\lambda^n$ with a stochastic relaxed control $\Lambda^n$.  We then fix a sequence of associated weak solutions $(\TO^{n},\TF^{n},(\TF^{n}_t)_{t\in[0,T]},\QB^n,X^n,W^n)$ of Eq.\eqref{eq:MFGRel} and set $\Theta^n\doteq\QB^n\circ (X^n,\Lambda^n)^{-1}\in\PP(\X\times\V)$ for all $n\in\N$. Finally, we associate to each MFG($n$) and to the limit MFG a martingale problem (\citet{stroock1969,stroock}) and show that the limit points $\Theta\in\PP(\X\times\V)$ of $(\Theta^n)_{n\in\N}$ solve the limit relaxed martingale problem. We start with the following lemma. 
\blem[\textit{Tightness in the 1-Wasserstein distance and absolute continuity}]\label{lem:conv1WassRel} Let $(\Theta^n)_{n\in\N}$ be as above. Then the following two properties hold:\begin{enumerate}
\item[(i)] $(\Theta^n)_{n\in\N}$ is tight in $\mathcal P_1 (\X\times \V)$; 
\item[(ii)] Any limit point $\Theta$ of the sequence $(\Theta^n)_{n\in\N}$ in $\mathcal P_1(\X \times \V)$ satisfies $\Theta\circ X^{-1}\ll\wiener$.\end{enumerate}
\elem
\begin{proof}\noindent\textit{(i).} It follows from Lemma \ref{lem:conv1Wass} and the compactness of $\Gamma$.\\
\indent\textit{(ii).} This is a consequence of Proposition \ref{prop:absContLimit}, the fact that by construction $\theta^n=\Theta^n\circ X^{-1}$ for all $n\in\N$, and the fact that weak convergence of the joint laws implies weak convergence of the marginals.
\end{proof}

\noindent By the previous lemma, we can assume without loss of generality that the original sequence $(\Theta^n)_{n \in \N}$ converges to some limit measure $\Theta$ in $\mathcal P_1 (\X \times \V)$. In order to characterize the limit point $\Theta$, we associate to each approximating MFG($n$) and to the limit MFG a (relaxed) martingale problem, henceforth RM($n$) and RM, respectively. Then, we show that $\Theta$ is also a solution of RM. We will use the notation $Dg$ and $D^2 g$ for the gradient and the Hessian of a smooth function $g: \R^d \to \R$, while $\textrm{Tr}[A]$ denote the trace of a square matrix $A$. Notice that in the following definition we have used the repameterization $b$ of the drift $\bar b$.

\bdefi{\textit{The approximating martingale problems (RM($n$))}}\label{def:relMartProbln} We say that $\widehat{\Theta}\in\PP(\X\times\V)$ is a solution of RM($n$) if for all $g\in \C^{2}_c(\Rd)$ the process
\ba 
M^{n,g}_t(\varphi,q;\widehat{\Theta})\doteq g(\varphi(t))-g(\varphi(0))-\int_{[0,t]\times\Gamma}b^n(s,\varphi,\hat{\theta},u)^\top D g(\varphi(s))q(ds,du)&&\nonumber\\
-\frac{1}{2}\int_0^t \text{Tr}\left[{\sigma \sigma^\top} D^2 g(\varphi(s))\right]ds,\quad t\in[0,T]&&
\nonumber 
\ea
\noindent is a $\widehat{\Theta}$-martingale, where $\hat{\theta}\doteq \widehat{\Theta}\circ X^{-1}$ and $X$ is the canonical process on $\X$. \edefi
\noindent Observe that, by construction, each $\Theta^n$ solves RM($n$). In Proposition \ref{prop:approxMart} below we will characterize the limit points as solutions of the following (relaxed) martingale problem.
\bdefi{\textit{The limit martingale problem (RM)}}\label{def:relMartProbl} We say that $\widehat{\Theta}\in\PP(\X\times\V)$ is a solution of RM if for all $g\in \C^{2}_c(\Rd)$ the process
\ba 
M^{g}_t(\varphi,q;\widehat{\Theta})\doteq g(\varphi(t))-g(\varphi(0))-\int_{[0,t]\times\Gamma}b(s,\varphi,\hat{\theta},u)^\top  D g(\varphi(s))q(ds,du)&&\nonumber\\
-\frac{1}{2}\int_0^t\text{Tr}\left[{\sigma \sigma^\top} D^2 g(\varphi(s))\right] ds,\quad t\in[0,T]&&
\nonumber 
\ea
\noindent is a $\widehat{\Theta}$-martingale, where $\hat{\theta}\doteq \widehat{\Theta}\circ X^{-1}$.
\edefi
\brem \label{rem:martProble}
\noindent The martingale property in both RM($n$) and in RM is understood to hold on $(\X\times\V,\B(\X\times\V))$ with respect to the $\Theta$-augmentation of the canonical filtration made right continuous by a standard procedure. Nonetheless, to conclude it is sufficient to check that the martingale property holds with respect to the canonical filtration on $\X\times\V$ (see, for instance, Problem 5.4.13 in \citet{karaztas}).
\erem
\noindent Now, we can characterize the limit points via the martingale problems.
\bprop[\textit{Characterization of limit points via martingale problems}]\label{prop:approxMart} $\Theta$ solves RM as in Definition \ref{def:relMartProbl}.
\eprop
\begin{proof}
\noindent Fix $t_1,t_2\in[0,T]$, $t_1<t_2$, $g\in \C^2_c(\Rd)$ and $\psi\in\C_b(\X\times\V)$ measurable with respect to $\B_{t_1}(\X\times\V)$. Define $\Psi,\Psi^n:\PP(\X\times\V)\rightarrow\R$ as
\ba 
\Psi\round{\Theta';\Theta}&\doteq &\E^{\Theta'}\brackets{\psi\round{M^{g}_{t_2}(\cdot\,;\Theta)-M^{g}_{t_1}(\cdot\,;\Theta)}},\nonumber\\
\Psi^n\round{\Theta';\Theta}&\doteq &\E^{\Theta'}\brackets{\psi\round{M^{n,g}_{t_2}(\cdot\,;\Theta)-M^{n,g}_{t_1}(\cdot\,;\Theta)}}
\nonumber
\label{eq:psi}
\ea
\noindent for $\Theta',\Theta\in\PP(\X\times\V)$ and for all $n\in\N$. Since $\Psi^n(\Theta^n;\Theta^n)=0$ for all $n\in\N$,
it suffices to prove that $\Psi^n(\Theta^n;\Theta^n)\rightarrow\Psi(\Theta;\Theta)$ as $n\rightarrow\infty$.\\
\indent First, we observe that $\Psi^n(\Theta^n;\Theta^n)$ and $\Psi(\Theta;\Theta)$ can be written as
\ba 
\Psi^n(\Theta^n;\Theta^n)&=&\int_{\X\times\V}\psi(\varphi,q)\int_{[t_1,t_2]\times\Gamma}b^n(s,\varphi,\theta^n,u)^\top D g(\varphi(s))q(ds,du)\Theta^n(d\varphi,dq)\nonumber\\
&&+ \int_{\X\times\V}\psi(\varphi,q)\int_{t_1}^{t_2} \frac{1}{2} \text{Tr}\left[{\sigma \sigma^\top} D^2 g(\varphi(s))\right] ds\,\Theta^n(d\varphi,dq)
\nonumber
\ea
\noindent and
\ba 
\Psi(\Theta;\Theta)&=&\int_{\X\times\V}\psi(\varphi,q)\int_{[t_1,t_2]\times\Gamma}b(s,\varphi,\theta,u)^\top D g(\varphi(s))q(ds,du)\Theta(d\varphi,dq)\nonumber\\
&&+\int_{\X\times\V}\psi(\varphi,q)\int_{t_1}^{t_2} \frac{1}{2}  \text{Tr}\left[{\sigma \sigma^\top} D^2 g(\varphi(s))\right] ds\,\Theta(d\varphi,dq).
\nonumber
\ea 
\noindent The convergence of the diffusion terms is a straightforward consequence of the weak convergence $\Theta^n\weakconv\Theta$ and the fact that the map
\ba 
(\varphi,q)\mapsto \psi(\varphi,q)\,\int_{t_1}^{t_2}\frac{1}{2}  \text{Tr}\left[{\sigma \sigma^\top} D^2 g(\varphi(s))\right] ds \nonumber
\ea
\noindent is in $C_b(\X\times\V)$, leading to
\ba 
\int_{\X\times\V}\psi(\varphi,q)\int_{t_1}^{t_2} \frac{1}{2}  \text{Tr}\left[{\sigma \sigma^\top} D^2 g(\varphi(s))\right] ds\,\Theta^n(d\varphi,dq)\nonumber\\
\conv\int_{\X\times\V}\psi(\varphi,q)\int_{t_1}^{t_2} \frac{1}{2}  \text{Tr}\left[{\sigma \sigma^\top} D^2 g(\varphi(s))\right] ds\,\Theta(d\varphi,dq).
\nonumber 
\ea
\noindent Hence, we only need to study the convergence of the drift terms. We split the rest of the proof in two steps.\\
\indent\textit{Step 1.}\quad We prove that
\ba 
\int_{\X\times\V}\psi(\varphi,q)\int_{[t_1,t_2]\times\Gamma}\round{b^n(s,\varphi,\theta^n,u)-b(s,\varphi,\theta^n,u)}^\top  Dg(\varphi(s))q(ds,du)\Theta^n(d\varphi,dq)\conv 0.
\nonumber
\ea
\noindent Indeed,
\begin{equation*}
\begin{split}
& \modl{\int_{\X\times\V} \psi(\varphi,q)\int_{[t_1,t_2]\times\Gamma}\round{b^n(s,\varphi,\theta^n,u)-b(s,\varphi,\theta^n,u)}^\top D g(\varphi(s))q(ds,du)\Theta^n(d\varphi,dq)} \\
 & \leq  C_{Dg}C_{\psi}\int_{\X\times\V}\int_{[t_1,t_2]\times\Gamma}\modl{b^n(s,\varphi,\theta^n,u)-b(s,\varphi,\theta^n,u)}q(ds,du)\Theta^n(d\varphi,dq) \\
& \leq C_{Dg}C_{\psi}\int_{\X\times\V}\int_{[t_1,t_2]\times\Gamma}\modl{b(s,\varphi,\theta^n,u)}\mathbf{1}_{\{|b|\geq K_n\}}q(ds,du)\Theta^n(d\varphi,dq)\\
& \leq C_{Dg}C_{\psi} \frac{\epsilon^{\alpha}\int_{\X\times\V}\int_{[t_1,t_2]\times\Gamma} \modl{b(s,\varphi,\theta^n,u)}^{\alpha}  q(ds,du)\Theta^n(d\varphi,dq)}{2\alpha}\\
& +C_{Dg}C_{\psi} \frac{\int_{\X\times\V}\int_{[t_1,t_2]\times\Gamma} \mathbf{1}_{\{|b|\geq K_n\}} q(ds,du)\Theta^n(d\varphi,dq)}{2\beta\epsilon^{\beta}}\\
& \leq C_{Dg}C_{\psi} \frac{\epsilon^{\alpha}\sup_{n\in\N}\int_{\X\times\V}\int_{[t_1,t_2]\times\Gamma} \modl{b(s,\varphi,\theta^n,u)}^{\alpha}  q(ds,du)\Theta^n(d\varphi,dq)}{2\alpha}\\
& +C_{Dg}C_{\psi} \frac{\sup_{n\in\N}\int_{\X\times\V}\int_{[t_1,t_2]\times\Gamma} \modl{b(s,\varphi,\theta^n,u)}  q(ds,du)\Theta^n(d\varphi,dq)}{2K_n\beta\epsilon^{\beta}}
\end{split}
\end{equation*}
\noindent for all $\epsilon>0$, where $C_{D g}$ and $C_{\psi}$ are uniform bounds on $D g$ and $\psi$, respectively. We applied Young's inequality with exponents $\alpha, \beta>1$, $\frac{1}{\alpha}+\frac{1}{\beta}=1$ for the third inequality, while for the last one we used the Markov's inequality with respect to the measure $\pi(ds,du,d\varphi,dq)\doteq q(ds,du)\Theta^n(d\varphi,dq)$ on $\X\times\V\times [0,T]\times\Gamma$:
\ba 
\int_{\X\times\V}\int_{[t_1,t_2]\times\Gamma} \mathbf{1}_{\{|b|\geq K_n\}} q(ds,du)\Theta^n(d\varphi,dq)\leq\frac{\int_{\X\times\V}\int_{[t_1,t_2]\times\Gamma} \modl{b(s,\varphi,\theta^n,u)}  q(ds,du)\Theta^n(d\varphi,dq)}{K_n}.
\nonumber
\ea
\noindent The suprema over $n\in\N$ are bounded due to Lemma \ref{lem:alphaEstimate}. We conclude this step by letting first $n\rightarrow\infty$ (so that $K_n\nearrow\infty$) then $\epsilon\rightarrow 0$.\\
\indent\textit{Step 2.}\quad We prove that
\ba 
\int_{\X\times\V}\psi(\varphi,q)\int_{[t_1,t_2]\times\Gamma}b(s,\varphi,\theta^n,u)^\top D g(\varphi(s))q(ds,du)\Theta^n(d\varphi,dq)\nonumber\\
\conv \int_{\X\times\V}\psi(\varphi,q)\int_{[t_1,t_2]\times\Gamma}b(s,\varphi,\theta,u)^\top D g(\varphi(s))q(ds,du)\Theta(d\varphi,dq).
\nonumber
\ea
\noindent To this aim we show that:
\ba
(\theta,\varphi,q)\mapsto \psi(\varphi,q)\int_{[t_1,t_2]\times\Gamma}b(s,\varphi,\theta,u)^\top D g(\varphi(s))q(ds,du)
\nonumber
\ea
\noindent is continuous on $\PP_1(\X)\times\X\times\V$ at points such that $\theta\ll\wiener$ and that it has sub-linear growth in $(\varphi,q)\in \X\times\V$ so that we can conclude by using the property $W_1(\Theta^n,\Theta)\rightarrow 0$ together with Theorem 7.12.iv in \cite{villani}. Since $\psi\in \C(\X\times\V)$, we only need to show the continuity of the second (integral) term. Let $(\theta^n,\varphi^n,q^n,u^n)_{n\in\N}\subset \PP_1(\X)\times\X\times\V\times\Gamma$ converge to some point $(\theta,\varphi,q,u)\in\PP_1(\X)\times\X\times\V\times \Gamma$ where $\theta\ll\wiener$. Then
\ba 
b(t,\varphi^n,\theta^n,u^n)^\top Dg(\varphi^n(t))\conv b(t,\varphi,\theta,u)^\top Dg(\varphi(t))
\nonumber
\ea
\noindent for all $t\in[t_1,t_2]$ by the continuity assumptions on $b$ and $D g$, i.e. $b(t,\cdot)^\top D g(\cdot)$ is jointly continuous for each $t\in[t_1,t_2]$ at points $(\theta,\varphi,q,u)$ with $\theta\ll\wiener$. Moreover
\ba 
\modl{b(t,\varphi,\theta,u)^\top D g(\varphi(t))}&\leq & C_{D g}C\round{1+\infnormt{\varphi}+m(t;\theta)+|u|}\nonumber\\
&\leq & C_{D g}C\round{1+K+\infnormt{\varphi}+|u|}
\nonumber
\ea
\noindent for some constants $C_{D g},C,K>0$ (this replaces Assumption (2) of Corollary A.5 in \cite{lacker}). We conclude by means of Corollary A.5 in \cite{lacker}.

\end{proof}
\noindent We conclude this subsection by characterizing any limit measure $\Theta$ as the joint law of state and (relaxed) control for a weak solution of the limit SDE in Eq.\eqref{eq:SDEopen} with drift $\bar b$. The next corollary is a fairly standard result establishing a well-known connection between solutions of RM and weak solutions of SDEs:
\bcor[\textit{Representation of limit points}]\label{cor:limitPointsRep} Let $\Theta$ be a solution of RM, as in Definition \ref{def:relMartProbl}. Then there exists a weak solution $(\TO,\TF,(\TF_t)_{t\in[0,T]},\QB,X,\Lambda,W)$ of
\ba
X_t &=& X_0 + \int_{[0,t]\times\Gamma}\bar{b}\round{s,X_s,\mu_s,u}\Lambda_s(du)ds+ \sigma W_t,\quad t\in[0,T]
\nonumber
\ea
\noindent such that $\Theta = \QB \circ (X,\Lambda)^{-1}$, $\theta=\Theta\circ X^{-1}$ and $\mu_t=g(t,\theta)$ with $g:[0,T]\times\PP_1(\X)\rightarrow\mathcal{M}_{\leq 1,1}(\R^d)$ as in Eq.\eqref{eq:defGtheta}.
\ecor
\begin{proof}
Arguing analogously as in the proofs of Proposition 5.4.6 and Corollary 5.4.8 in \cite{karaztas} gives the existence of a weak solution $(\TO,\TF,(\TF_t)_{t\in[0,T]},\QB,X,\Lambda,W)$ of the SDE
\ba
X_t &=& X_0 + \int_{[0,t]\times\Gamma}b\round{s,X,\theta,u}\Lambda_s(du)ds+ \sigma W_t,\quad t\in[0,T]
\label{eq:MFGRel2}
\ea
\noindent such that $\Theta$ is the law of $(X,\Lambda)$ under $\QB$ and $\theta=\Theta\circ X^{-1}$. The conclusion is obtained by going back to the original drift $\bar b$, that we recall is given by
\[ \bar{b}(t,\varphi(t),g(t,\theta),u) = b(t,\varphi,\theta,u), \quad (t,\varphi, \theta, u) \in [0,T] \times \X \times \mathcal P_1 (\X) \times \Gamma,\]
and $g(t, \theta) = \mu_t$ as in Eq.\eqref{eq:defGtheta}.
\end{proof}

\subsection{Optimality of the limit points}\label{subsec:optimality}
\noindent In this subsection, we show that any limit point $\Theta\in\PP(\X\times\V)$ of $(\Theta^n)_{n\in\N}$ is optimal according to the cost functional of the MFG. In order to do that, we will extend the notion of relaxed MFG solution to controls that are not necessarily in feedback form. In this case we evaluate optimality according to the following cost functional:
\ba 
J^{\mu}\round{\Lambda}&\doteq &\E\Biggl[\int_{[0,\tau]\times\Gamma}\bar{f}\round{s,X_s,\mu_s,u}\Lambda_s(du)ds+F\round{\tau,X_{\tau}}\Biggr],
\nonumber
\ea
\noindent where $\Lambda$ is any relaxed stochastic control and $\tau \doteq \tau^{X} \wedge T$, subject to the dynamics
\ba
X_t &=& X_0 + \int_{[0,t]\times\Gamma}\bar{b}\round{s,X_s,\mu_s,u}\Lambda_s(du)ds+ \sigma W_t,\quad t\in[0,T].
\label{eq:MFGRel3}
\ea
\noindent We set $V^{\mu} =\inf_{\Lambda} J^{\mu}(\Lambda)$, where the minimization is actually performed over the set of relaxed stochastic open-loop controls, i.e. over the tuples $(\TO,\TF,(\TF_t)_{t\in[0,T]},\QB,X,\Lambda,W)$ that are weak solutions of Eq.\eqref{eq:MFGRel3} and where $\Lambda$ is a progressively measurable $\PP(\Gamma)$-valued stochastic process. To simplify the notation, we will just write $\Lambda$ to refer to the whole tuple. Moreover, when working on the canonical space $\X\times\V$, where the canonical process $(X,\Lambda)$ is completely characterized by its law $\Theta$, we will simply write $J^{\mu}(\Theta)$ in place of $J^{\mu}(\Lambda)$.
\bdefi[\textit{Relaxed MFG solution}]\label{def:MFGsolRel2} A relaxed solution of the MFG is a pair $(\Lambda,\mu)$, where $\Lambda$ is a relaxed stochastic control and $\mu\in\Upsilon_{\leq 1,1}^{T}$, such that:
\bi
\item[(i)] $\Lambda$ is optimal, i.e. $V^{\mu} = J^{\mu}(\Lambda)$.
\item[(ii)] Let $(\TO,\TF,(\TF_t)_{t\in[0,T]},\QB,X,\Lambda,W)$ be a weak solution of Eq.\eqref{eq:MFGRel3} with flow of sub-probability measures $\mu$, stochastic control $\Lambda$ and initial condition $\nu$. Then 
\[ \mu_t(\cdot)=\QB(\{X_t \in \cdot\}\cap\{ \tau^X > t\}), \quad t\in[0,T].\]
\ei
\edefi
\bprop[\textit{Existence of relaxed MFG solutions}]\label{prop:existenceRelSolMFG} Grant Assumptions (H1)-(H8) and (C1). Let $\Theta$ be a limit point of $(\Theta^n)_{n\in\mathbb{N}}$ in $\PP_1(\X\times \V)$. Set $\mu\in\Upsilon^T_{\leq 1,1}$ as
\ba 
\mu_t\round{\cdot}\doteq\Theta\round{\cbr{X_t\in \cdot}\cap\cbr{\tau^{X}>t}}\quad t\in[0,T].
\nonumber
\ea
\noindent Then $(\Theta,\mu)$ is a relaxed MFG solution according to Definition \ref{def:MFGsolRel2}.
\eprop
\begin{proof}\noindent By construction we immediately have that $\Lambda$ is a relaxed stochastic control and $\mu\in\Upsilon_{\leq 1,1}^{T}$. Moreover, property (ii) is a consequence of the fact that $\Theta$ is a solution of RM as in Definition \ref{def:relMartProbl}. To prove property (i), we proceed through the following steps:
\bi 
\item[(j)]\quad Let $\tilde{\Theta}\in\PP(\X\times\V)$ be a solution of RM. Then there exists a sequence of solutions $(\tilde{\Theta}^n)_{n\in\N}$ of RM($n$) such that $\lim_{n\rightarrow\infty}J^{n,\mu^n}(\tilde{\Theta}^n)=J^{\mu}(\tilde{\Theta})$.
\item[(jj)]\quad $\lim_{n\rightarrow\infty}J^{n,\mu^n}(\Theta^n)=J^{\mu}(\Theta)$.
\item[(jjj)]\quad $J^{\mu}(\Theta)\leq J^{\mu}(\tilde{\Theta})$ for any $\tilde{\Theta}\in\PP(\X\times\V)$ solution of RM.
\ei
\noindent The proof of (j)-(jjj) largely follows that of Theorem 3.6 in \cite{lacker}. Therefore, we highlight only the main differences with respect to our setting, which are due to the sub-linear growth of the drift and the cost functional and to the path dependency induced by the exit time from $\mathcal{O}$.\\
\indent\textit{Proof of (j).}\quad Let $\tilde{\Theta}\in\PP(\X\times\V)$ be a solution of RM 
and let $(\tilde{\Omega},\tilde{\F},(\tilde{\F}_t)_{t\in[0,T]},\tilde{\Theta},X,\Lambda,W)$ be a weak solution of Eq.\eqref{eq:MFGRel3} on the canonical space $\TO=\X\times\V$. The existence of this solution is guaranteed by Corollary \ref{cor:limitPointsRep}. Now fix $\Lambda$ and let $X^n$ be a sequence of strong solutions of:
\ba 
X^n_t=\xi+\int_{[0,t]\times\Gamma}\bar{b}^n\round{s,X_s^n,\mu_s^n,u}\Lambda_s(du)ds+\sigma W_t,\quad t\in[0,T]
\nonumber
\ea
on the filtered probability space $(\tilde{\Omega},\tilde{\F},(\tilde{\F}_t)_{t\in[0,T]},\tilde{\Theta})$. Set $\tilde{\Theta}^n\doteq\tilde{\Theta}\circ (X^n,\Lambda)^{-1}$ for each $n\in\mathbb{N}$. Notice that $(\tilde{\Theta}^n)_{n\in\mathbb{N}}\subset\PP_1(\mathcal{X}\times\mathcal{V})$. Moreover each $\tilde{\Theta}^n$ solves RM($n$) as in Definition \ref{def:relMartProbln}. We now show that:
\ba 
\mathbb{E}^{\tilde{\Theta}}\brackets{\infnorm{X^n-X}}\conv 0\quad\text{and}\quad W_1(\tilde{\Theta}^n,\tilde{\Theta})\conv 0.
\label{eq:convL1}
\ea
Regarding the first limit, it is sufficient to note that:
\ba 
\E^{\tilde{\Theta}}\brackets{\infnormt{X^n-X}}&\leq & L\int_0^t\E^{\tilde{\Theta}}\brackets{\infnorms{X^n-X}}ds+\E^{\tilde{\Theta}}\brackets{\int_{[0,t]\times\Gamma}\Delta b^n(s,u)\Lambda_s(du)ds}
\nonumber
\ea
\noindent where we set
\ba 
\Delta b^n(t,u)\doteq |\bar{b}^n(t,X_t,\mu_t,u)-\bar{b}(t,X_t,\mu_t,u)|.
\nonumber
\ea 
\noindent The first term can be handled with Gr\"onwall's Lemma, whereas the second one by applying a similar argument as in the first step of the proof of Proposition \ref{prop:approxMart}. Regarding the second limit in Eq.\eqref{eq:convL1} we can proceed as follows. First, notice that the first limit in Eq.\eqref{eq:convL1} implies convergence in probability, hence in law, of $X^n$ to $X$. 
Thus, by an argument similar to that of Lemma \ref{lem:conv1WassRel}, we can prove the convergence in the 1-Wasserstein distance. At this point, the convergence of the costs is a consequence of the convergence in the 1-Wasserstein distance and the sub-linear growth of the running cost (combined with Theorem 7.12.iv in \cite{villani}), as in the second step of the proof of Proposition \ref{prop:approxMart}.\\
\indent\textit{Proof of (jj).}\quad This follows from an argument similar to the second part of (j).\\
\indent\textit{Proof of (jjj).}\quad Let $\tilde{\Theta}\in\PP(\X\times\V)$ be a solution of RM and let $(\tilde{\Theta}^n)_{n\in\N}\subset\PP(\X\times\V)$ be an approximating sequence as in (j). By the optimality of $\Theta^n$ we have
\ba 
J^{n,\mu^n}\round{\Theta^n}\leq J^{n,\mu^n}\round{\tilde{\Theta}^n}
\nonumber
\ea
\noindent for all $n\in\N$. The optimality of $\Theta$ follows by taking the limit for $n\rightarrow\infty$ on both sides of the inequality above and using the previous properties (j) and (jj).
\end{proof}
\subsection{Existence of solutions}\label{sec:existenceSol}
In this subsection we finally conclude the proof of Theorem \ref{teo:existenceRelFeedSolMFG} by proving the existence of a relaxed feedback MFG solution and, under additional convexity assumptions, the existence of a strict feedback MFG solution. In addition, we also prove existence of solutions that are Markovian up to the exit time.\\

\noindent\textit{Relaxed feedback MFG solutions.}\quad The main mathematical tool here is the mimicking result of \cite{brunick}. We follow the procedure in \cite{lacker} but with modifications due to the peculiarities of our model induced mainly by the presence of absorptions. We give more details in the proof below.
\bprop[\textit{Existence of relaxed feedback MFG solutions}]\label{prop:existenceRelFeedlSolMFG} Grant Assumptions  \hyperlink{H1}{(\textrm{H1})}- \hyperlink{H8}{(\textrm{H8})} and  \hyperlink{C1}{(\textrm{C1})}. Let $(\Theta,\mu)$ be a relaxed MFG solution as in Definition \ref{def:MFGsolRel2}. 

Then there exists another relaxed MFG solution $(\Theta',\mu)$ and a progressively measurable functional $\lambda:[0,T]\times\X\rightarrow\PP(\Gamma)$ such that $\Theta'((\varphi,q)\in\X\times\V:q_t=\lambda(t,\varphi))=1$ for $\mathcal L_T$-a.e. $t\in[0,T]$ and $J^{\mu}(\Theta')=J^{\mu}(\Theta)=V^{\mu}$, i.e. $(\lambda,\mu)$ is a relaxed feedback solution of the MFG as in Definition \ref{def:MFGsolRel}.
\eprop
\begin{proof}
\noindent We adapt the proof of Theorem 3.7 in \cite{lacker} to our setting, by exploiting the mimicking result in Corollary 3.11 of \cite{brunick} instead of Corollary 3.7 as in \cite{lacker}. As a consequence, the mimicking process that we get is not Markovian as in Lacker. However, it has the same law as the original process and not only the same marginals. This is important in our setting due to the path dependency induced by the exit time $\tau$.

We start with the construction of $\lambda$ by disintegration. Precisely, define $\eta\in\PP([0,T]\times\mathcal{X}\times\Gamma)$ as:
\ba 
\eta\round{I\times B\times G}&\doteq &\frac{1}{T}\mathbb{E}^{\Theta}\brackets{\int_{[0,T]\times\Gamma}\mathbf{1}_{\round{I\times B\times G}}\round{t,X,u}\Lambda\round{dt,du}}
\nonumber\ea
\noindent and disintegrate it as $\eta(dt,d\varphi,du)=\tilde{\eta}(dt,d\varphi)\lambda_{t,\varphi}(du)$. Then:
\ba
\eta\round{I\times B\times G}
&=&\int_{[0,T]\times\mathcal{X}}\int_{\Gamma}\mathbf{1}_{\round{I\times B\times G}}\round{t,\varphi,u}\lambda_{t,\varphi}\round{du}\tilde{\eta}\round{dt,d\varphi}
\nonumber
\ea
\noindent for all $I\in\mathcal{B}([0,T])$, $B\in\mathcal{B}(\mathcal{X})$ and $G\in\mathcal{B}(\Gamma)$. By the disintegration theorem, $(t,\varphi)\mapsto \lambda_{t,\varphi}(\cdot)\in\PP(\Gamma)$ is Borel-measurable. Now set $\tilde{\F}^X_t\doteq \sigma(X_s,s\in[0,t])$ for each $t\in[0,T]$. We claim that:
\ba 
\lambda_{t,X}\round{\cdot}&=&\mathbb{E}^{\Theta}\brackets{\Lambda_t\round{\cdot}\big| \tilde{\F}^X_t}\quad\Theta\text{-a.s. and for $\mathcal L_T$-a.e.}\,t\in[0,T]
\label{eq:lambdaProgmeas}
\ea
\noindent which is measurable and adapted, hence it has a progressively measurable modification $\lambda$. We show that for any bounded measurable functional $g:[0,T]\times\mathcal{X}\times\Gamma\rightarrow\R$ such that $g(t,\cdot,u)$ is $\tilde{\F}^X_t$-measurable for all $t\in[0,T]$ and $u\in\Gamma$
\ba 
\int_{\Gamma}g\round{t,X,u}\lambda_{t,X}\round{du}
&=&\int_{\Gamma}g\round{t,X,u}\mathbb{E}^{\Theta}\brackets{\Lambda_t\round{du}\big| \tilde{\F}^X_t}
\nonumber
\ea
\noindent $\Theta\text{-a.s. and for $\mathcal L_T$-a.e.}\,t\in[0,T]$. Indeed, for any other bounded measurable functional $h:[0,T]\times\mathcal{X}\rightarrow\R$ such that $h(t,\cdot)$ is $\tilde{\F}^X_t$-measurable for all $t\in[0,T]$, we have
\ba 
&&\frac{1}{T}\mathbb{E}^{\Theta}\brackets{\int_0^Th\round{t,X}\int_{\Gamma}g\round{t,X,u}\lambda_{t,X}\round{du}dt}\\
&=&\int_{[0,T]\times\mathcal{X}}h\round{t,\varphi}\int_{\Gamma}g\round{t,\varphi,u}\lambda_{t,\varphi}\round{du}\tilde{\eta}\round{dt,d\varphi}\nonumber\\
&=&\int_{[0,T]\times\mathcal{X}\times\Gamma}h\round{t,\varphi}g\round{t,\varphi,u}\eta\round{dt,d\varphi,du}\nonumber\\
&=&\frac{1}{T}\mathbb{E}^{\Theta}\brackets{\int_0^Th\round{t,X}\int_{\Gamma}g\round{t,X,u}\Lambda_t\round{du}dt}
\nonumber
\ea
\noindent where the first equality comes from the definition of $\tilde{\eta}$, the second one is due to the disintegration of $\eta$ and the third one holds by definition of $\eta$.\\
\noindent Now, let $(\tilde{\Omega},\tilde{\F},(\tilde{\F}_t)_{t\in[0,T]},\QB,W,X,\Lambda)$ be a weak solution of Eq.\eqref{eq:MFGRel3} with relaxed control $\Theta=\QB\circ(X,\Lambda)^{-1}$. By Corollary 3.11 in \cite{brunick} there exists a weak solution $(\tilde{\Omega}',\tilde{\F}',(\tilde{\F}'_t)_{t\in[0,T]},\QB',W',X')$ of
\ba 
X'_t &=&\xi+\int_0^{t}\int_{\Gamma}\bar{b}\round{s,X'_s,\mu_s,u}\lambda_{s,X'}(du)ds+\sigma W'_t,\quad t\in[0,T]
\nonumber
\ea
\noindent such that $\QB'\circ (X')^{-1}=\QB\circ X^{-1}$. Define $\Theta'\doteq \QB'\circ (X',\Lambda')^{-1}$ where $\Lambda'(dt,du)\doteq dt\lambda_{t,X'}(du)$. Notice that if $\mu'$ is the flow of sub-probability measures associated to $\Theta'$ then $\mu'=\mu$. Finally, $\Theta'$ solves the same relaxed martingale problem as $\Theta$, and it has the same cost as $\Theta$ as required:
\ba 
J^{\mu}\round{\Theta'}&=& 
\mathbb{E}^{\QB'}\brackets{\int_0^{\tau'}\int_{\Gamma}\bar{f}\round{t,X'_t,\mu_t,u}\lambda_{t,X'}\round{du}dt+F\round{\tau',X'_{\tau'}}}\nonumber\\
&=&\mathbb{E}^{\QB}\brackets{\int_0^{\tau}\int_{\Gamma}\bar{f}\round{t,X_t,\mu_t,u}\lambda_{t,X}\round{du}dt+F\round{\tau,X_{\tau}}}\nonumber\\
&=&\mathbb{E}^{\QB}\brackets{\int_0^{\tau}\int_{\Gamma}\bar{f}\round{t,X_t,\mu_t,u}\mathbb{E}^{\QB}\brackets{\Lambda_t\round{du}\big| \tilde{\F}^X_t}dt+F\round{\tau,X_{\tau}}}\nonumber\\
&=&\mathbb{E}^{\QB}\brackets{\int_0^{\tau}\int_{\Gamma}\mathbb{E}^{\QB}\brackets{\bar{f}\round{t,X_t,\mu_t,u}\Lambda_t\round{du}\big| \tilde{\F}^X_t}dt+F\round{\tau,X_{\tau}}}\nonumber\\
&=&\mathbb{E}^{\QB}\brackets{\int_{[0,\tau]\times\Gamma}\bar{f}\round{t,X_t,\mu_t,u}\Lambda\round{dt,du}+F\round{\tau,X_{\tau}}}\nonumber\\
&=&J^{\mu}\round{\Theta}.
\nonumber 
\ea
\end{proof}
\brem \label{rem:brunick-shreve}
\noindent We observe that, due to the discontinuity induced by the exit time $\tau$, it is not possible in general to apply Theorem 3.6 of \cite{brunick} to $Z_t=(X_t,\mathbb{I}_{[0,\tau)}(t))$, $t\in[0,T]$, to obtain a control which is Markovian in $Z$. Moreover the few mimicking results available in the literature for discontinuous processes hold under very restrictive or hardly verifiable assumptions. Nonetheless, Theorem 3.6 of \cite{brunick} could still be applied in some particular cases when, for instance, $\mathcal{O}=(0,\infty)$ and $Z_t=(X_t,\inf_{s\in[0,t]} X_s)$.
\erem

\noindent\textit{Strict feedback MFG solutions.}\quad Under additional convexity assumptions (\citet{filippov,haussmann}), we prove existence of feedback MFG solutions in strict form. Let $(\Theta,\mu)$ be a relaxed MFG solution according to Definition \ref{def:MFGsolRel2} and
for each $(t,\varphi)\in[0,T]\times\X$ define $K(t,\varphi,\mu)$ as:
\ba 
K\round{t,\varphi,\mu}\doteq\cbr{\round{\bar{b}\round{t,\varphi(t),\mu_t,u},z}\,:\,z\geq \bar{f}\round{t,\varphi(t),\mu_t,u}\quad\text{and}\quad u\in\Gamma}.
\label{eq:setKConvex}
\nonumber
\ea
\noindent Existence of strict MFG solutions is established under the additional Assumption  \hyperlink{C2}{(\textrm{C2})}.
\brem
\noindent Assumption \hyperlink{C2}{(\textrm{C2})} is equivalent to requiring that the set $K(t,\varphi,\mu)$ is convex. This assumption is crucial to apply the measurable selection arguments in \cite{haussmann,dufour}.
\erem
\bprop[\textit{Existence of strict feedback MFG solutions}]\label{prop:existenceStrictFeedlSolMFG} Grant Assumptions \hyperlink{H1}{(\textrm{H1})}- \hyperlink{H8}{(\textrm{H8})},  \hyperlink{C1}{(\textrm{C1})} and Assumption \hyperlink{C2}{(\textrm{C2})}. Let $(\Theta,\mu)$ be a relaxed MFG solution as in Definition \ref{def:MFGsolRel2}. 

Then there exists another relaxed MFG solution $(\Theta',\mu)$ and a progressively measurable functional $u\in\mathcal{U}_{fb}$ such that $\Theta'((\varphi,q)\in\X\times\V:q_t=\delta_{u(t,\varphi)})=1$ for $\mathcal L_T$-a.e. $t\in[0,T]$ and $J^{\mu}(\Theta')=J^{\mu}(\Theta)=V^{\mu}$, i.e. $(u,\mu)$ is a strict and feedback solution of the MFG as in Definition \ref{def:MFGsol}.
\eprop
\begin{proof}
\noindent We follow once more the proof of Theorem 3.7 in \cite{lacker}, highlighting the main differences with respect to our setting. The first part of the proof proceeds as in Proposition \ref{prop:existenceRelFeedlSolMFG}. Since for all $(t,\varphi)\in[0,T]\times\mathcal{X}$ the pair $(\bar{b}(t,\varphi(t),\mu_t,u),\bar{f}(t,\varphi(t),\mu_t,u))$ belongs to $K(t,\varphi,\mu)$ for all $u\in\Gamma$ and $K(t,\varphi,\mu)$ is convex, we have
\ba 
\int_{\Gamma}\round{\bar{b}\round{t,\varphi(t),\mu_t,u},\bar{f}\round{t,\varphi(t),\mu_t,u}}\lambda_{t,\varphi}(du)\,\in K\round{t,\varphi,\mu}.
\nonumber
\ea
\noindent By applying the measurable selection argument in \cite{haussmann,dufour} (with respect to the progressive $\sigma$-algebra, i.e. the $\sigma$-algebra generated by progressively measurable processes), we find a progressively measurable functional $u:[0,T]\times\mathcal{X}\rightarrow \Gamma$ such that
\ba 
\int_{\Gamma}\bar{b}\round{t,\varphi(t),\mu_t,u}\lambda_{t,\varphi}(du)=\bar{b}\round{t,\varphi(t),\mu_t,u(t,\varphi)}
\nonumber
\ea
\noindent and
\ba 
\int_{\Gamma}\bar{f}\round{t,\varphi(t),\mu_t,u}\lambda_{t,\varphi}(du)\geq \bar{f}\round{t,\varphi(t),\mu_t,u(t,\varphi)}
\label{eq:conv-f}
\ea
\noindent for all $(t,\varphi)\in[0,T]\times\mathcal{X}$. Define $\Theta'\doteq \QB'\circ(X',\Lambda')^{-1}$ where $\QB'$ is as in the proof of Proposition \ref{prop:existenceRelFeedlSolMFG} and $\Lambda'(\varphi,q)(dt,du)\doteq dt\delta_{u(t,\varphi)}(du)$. $\Theta'$ solves the same relaxed martingale problem as $\Theta$. As for the costs, we have\ba 
J^{\mu}\round{\Theta'}&=& 
\mathbb{E}^{\QB'}\brackets{\int_0^{\tau'}\int_{\Gamma}\bar{f}\round{t,X'_t,\mu_t,u}\delta_{u(t,X')}(du)dt+F\round{\tau,X'_{\tau}}}\nonumber\\
&=&\mathbb{E}^{\QB'}\brackets{\int_0^{\tau'}\bar{f}\round{t,X'_t,\mu_t,u(t,X')}dt+F\round{\tau,X'_{\tau}}}\nonumber\\
&\leq &\mathbb{E}^{\QB'}\brackets{\int_0^{\tau'}\int_{\Gamma}\bar{f}\round{t,X'_t,\mu_t,u}\lambda_{t,X'}(du)dt+F\round{\tau,X'_{\tau}}}\nonumber\\
&=&\mathbb{E}^{\QB}\brackets{\int_0^{\tau}\int_{\Gamma}\bar{f}\round{t,X_t,\mu_t,u}\lambda_{t,X}\round{du}dt+F\round{\tau,X_{\tau}}}\nonumber\\
&=&\mathbb{E}^{\QB}\brackets{\int_{[0,\tau]\times\Gamma}\bar{f}\round{t,X_t,\mu_t,u}\Lambda\round{dt,du}+F\round{\tau,X_{\tau}}}\nonumber\\
&=&J^{\mu}\round{\Theta}
\nonumber 
\ea
\noindent where the inequality above is due to Eq.\eqref{eq:conv-f}. Given the optimality of $(\Theta,\mu)$ we already have the converse inequality, i.e. $J^{\mu}(\Theta)\leq J^{\mu}(\Theta')$. Hence $J^{\mu}(\Theta)=J^{\mu}(\Theta')$.
\end{proof}

\noindent We can finally give the proof of Theorem \ref{teo:existenceRelFeedSolMFG}.
\begin{proof}[Proof of Theorem \ref{teo:existenceRelFeedSolMFG}]
\noindent Grant Assumptions \hyperlink{H1}{(\textrm{H1})}-\hyperlink{H8}{(\textrm{H8})} and \hyperlink{C1}{(\textrm{C1})}. Proposition \ref{prop:existenceRelSolMFG} guarantees existence of a relaxed MFG solution $(\Theta,\mu)$ as in Definition \ref{def:MFGsolRel2}. By Proposition \ref{prop:existenceRelFeedlSolMFG} there exists another relaxed MFG solution $(\Theta',\mu)$ together with a progressively measurable functional $\lambda:[0,T]\times\X\rightarrow\PP(\Gamma)$ such that $\Theta'((\varphi,q)\in\X\times\V:q_t=\lambda(t,\varphi))=1$ for $\mathcal L_T$-a.e. $t$ and $J^{\mu}(\Theta')=J^{\mu}(\Theta)=V^{\mu}$. Then $(\lambda,\mu)$ is a relaxed and feedback solution of the MFG as in Definition \ref{def:MFGsolRel}. 

Additionally grant Assumption \hyperlink{C2}{(\textrm{C2})}. By Proposition \ref{prop:existenceStrictFeedlSolMFG} there exists another relaxed MFG solution $(\Theta',\mu)$ and a progressively measurable functional $u\in\mathcal{U}_{fb}$ such that $\Theta'((\varphi,q)\in\X\times\V:q_t=\delta_{u(t,\varphi)})=1$ for $\mathcal L_T$-a.e. $t \in [0,T]$, and $J^{\mu}(\Theta')=J^{\mu}(\Theta)=V^{\mu}$. Then $(u,\mu)$ is a strict and feedback solution of the MFG as in Definition \ref{def:MFGsol}.
\end{proof}

\noindent\textit{Markovian MFG solutions.}\quad We conclude this part with showing that  there exist relaxed and strict feedback solutions that are Markovian up to the exit time.
\bprop[\textit{Markovian MFG solutions}]\label{lem:existencemarkovSolMFG} Grant Assumptions \hyperlink{H1}{(\textrm{H1})}-\hyperlink{H8}{(\textrm{H8})} and \hyperlink{C1}{(\textrm{C1})}. Let $(\Theta,\mu)$ be a relaxed MFG solution as in Definition \ref{def:MFGsolRel2}. Then there exists another relaxed MFG solution $(\Theta',\mu)$ and a function $\lambda:[0,T]\times\Rd\rightarrow\PP(\Gamma)$ such that
\[ \mathcal L_T \otimes \Theta^\prime (\{(t,\varphi,q): q_t =\lambda(t,\varphi(t)), t \le \tau^X (\varphi)\})=1\]
and $J^{\mu}(\Theta')=J^{\mu}(\Theta)=V^{\mu}$. Additionally, grant Assumption \hyperlink{C2}{(\textrm{C2})}. Then there exists a function $u:[0,T]\times\Rd\rightarrow\Gamma$ such that
\[ \mathcal L_T \otimes \Theta^\prime (\{(t,\varphi,q): q_t =\delta_{u(t,\varphi(t))}, t \le \tau^X (\varphi)\})=1\]
and $J^{\mu}(\Theta')=J^{\mu}(\Theta)=V^{\mu}$.
\eprop
\begin{proof}
\noindent Let us define the following processes
\ba 
Y_t\doteq (t,X_t),\quad
X^{\tau^X}_t\doteq X_{t\wedge\tau^X},\quad
Y^{\tau^X}_t\doteq  Y_{t\wedge\tau^X}
\nonumber
\ea
\noindent for $t\in[0,T]$. If $X$ satisfies Eq.\eqref{eq:MFGRel3} with flow of sub-probability measures $\mu$ and relaxed control $\Lambda$ then the SDE satisfied by $X^{\tau^X}$ is (on the same probability space)
\ba 
X^{\tau^X}_t &=& \xi+\int_{[0,t]\times\Gamma}\bar{b}\round{s,X_s^{\tau^X},\mu_s,u}\mathbf{1}_{[0,\tau^X)}(s)\Lambda_s(du)ds +\sigma\int_0^t\mathbf{1}_{[0,\tau^X)}(s)dW_s
\label{eq:stoppedProcess}
\nonumber
\ea
for $t\in[0,T]$. Notice that until $t\leq \tau^X$ the stopped process $X^{\tau^X}$ coincides pathwise with the original process $X$. We now apply the mimicking result in  Corollary 3.7 of  \cite{brunick}, to the stopped process $Y^{\tau^X}$. To this end, we follow the proof of Theorem 3.7 in \cite{lacker} and the proofs of Propositions \ref{prop:existenceRelFeedlSolMFG} and \ref{prop:existenceStrictFeedlSolMFG} in the present paper.\\
First, we claim that there exists a measurable function $\lambda:[0,T]\times\R^{d+1}\rightarrow\PP(\Gamma)$ such that
\ba 
\lambda_{t,Y_t^{\tau^X}}(\cdot)&=&\mathbb{E}^{\Theta}\brackets{\Lambda_t(\cdot)\big| Y^{\tau^X}_t},\quad\Theta\text{-a.s. and for $\mathcal L_T$-a.e.}\,t\in[0,T].
\nonumber
\ea
Such a function can be constructed by disintegration as follows. Let $\eta\in\PP([0,T]\times\R^{d+1}\times\Gamma)$ be given by\ba 
\eta (B)&\doteq &\frac{1}{T}\mathbb{E}^{\Theta}\brackets{\int_{[0,T]\times\Gamma}\mathbf{1}_C\round{t,Y_t^{\tau^X},u}\Lambda(dt,du)}.
\nonumber
\ea
\noindent We define $\lambda$ through $\eta(dt,dy,du)\doteq\tilde{\eta}(dt,dy)\lambda_{t,y}(du)$. By Corollary 3.7 in \cite{brunick} applied to $\lambda_{t,Y_t^{\tau^X}}$ there exists a weak solution $(\tilde{\Omega}',\tilde{\F}',(\tilde{\F}'_t)_{t\in[0,T]},\QB',W',X')$ of
\ba 
X'_t =\xi+\int_0^{t}\int_{\Gamma}\bar{b}\round{s,X'_s,\mu_s,u}\mathbf{1}_{[0,\tau^{X'})}(s)\lambda_{s,Y_t^{\tau^{X'}}}(du)ds +\sigma\int_0^t\mathbf{1}_{[0,\tau^{X'})}(s) dW'_s\nonumber
\ea
for $t\in[0,T]$, where $Y_t^{\tau^{X'}}\doteq(t\wedge\tau^{X'},X'_t)$ and $\QB'\circ (t\wedge\tau^{X'},X'_t)^{-1}=\QB\circ(t\wedge\tau^X,X^{\tau^X}_t)^{-1}$ for all $t\in[0,T]$, i.e. $Y^{\tau^{X'}}$ and $Y^{\tau^{X}}$ have the same time marginals. Now set $\tau'\doteq \tau^{X'}\wedge T$. Recall that $\Theta=\mathbb{Q}\circ (X,\Lambda)^{-1}$ and define $\Theta'\doteq \QB'\circ (X',\Lambda')^{-1}$ where $\Lambda'(dt,du)\doteq dt\lambda_{t,Y_t^{\tau^{X'}}}(du)$. Equality of the costs can be shown just as in the proof of Proposition \ref{prop:existenceRelFeedlSolMFG}:
\ba
J^{\mu}\round{\Theta'}&=& 
\mathbb{E}^{\mathbb{Q}'}\brackets{\int_0^{\tau'}\int_{\Gamma}\bar{f}(t,X'_t,\mu_t,u)\lambda_{t,t\wedge\tau^{X'},X'_t}(du)dt+F\round{\tau',X'_{\tau'}}}\nonumber \\
&=&\mathbb{E}^{\mathbb{Q}}\brackets{\int_0^{\tau}\int_{\Gamma}\bar{f}(t,X_t^{\tau^X},\mu_t,u)\lambda_{t,t\wedge\tau^X,X^{\tau^X}_t}(u)dt+F\round{\tau,X_{\tau}^{\tau^X}}}\nonumber \\
&=&\mathbb{E}^{\mathbb{Q}}\brackets{\int_{[0,\tau]\times\Gamma}\bar{f}(t,X_t^{\tau^X},\mu_t,u)\Lambda(dt,du)+F\round{\tau,X_{\tau}^{\tau^X}}}\nonumber\\
&=& J^{\mu}\round{\Theta}.
\nonumber
\ea
Therefore, $\lambda:[0,T]\times[0,T]\times\Rd\rightarrow\PP(\Gamma)$ satisfies $\Theta'(q\in\mathcal{V}:q_t=\lambda(t,t\wedge\tau^{\hat{X}},\hat{X}^{\tau^{\hat{X}}}_t))=1$ for $\mathcal L_T$-a.e. $t\in[0,T]$ and $J^{\mu}(\Theta')=J^{\mu}(\Theta)=V^{\mu}$.

Consider now a weak solution $(\tilde{\Omega}'',\tilde{\F}'',(\tilde{\F}''_t)_{t\in[0,T]},\QB'',W'',X'')$ of
\ba 
X''_t =\xi+\int_0^{t}\int_{\Gamma}\bar{b}\round{s,X''_s,\mu_s,u}\lambda_{s,Y_t^{\tau^{X''}}}(du)ds+\sigma W''_t,\quad t\in[0,T]
\nonumber
\ea
\noindent where $Y_t^{\tau^{X''}}=(t\wedge\tau^{X''},X''_t)$. Set $\Theta''\doteq \QB''\circ (X'',\Lambda'')^{-1}$ where  $\Lambda''(dt,du)\doteq dt\lambda_{t,Y_t^{\tau^{X''}}}(du)$. To avoid confusion between specific solutions, here $(\hat{X},\hat{\Lambda})$ denotes the canonical process on $\X\times\V$. First, $\Theta'$ solves the martingale problem associated to
\ba
\widehat{M}^{g}_t(\varphi,q)\doteq g(\varphi(t))-g(\varphi(0))
-\int_{[0,t]\times\Gamma}\bar{b}(s,\varphi(s),\mu_s,u)^\top D g(\varphi(s))\mathbf{1}_{[0,\tau^{\hat{X}})}(s)q(ds,du)&&\nonumber\\
+\frac{1}{2}\int_0^t\text{Tr}\left[{\sigma \sigma^\top} D^2 g(\varphi(s))\right]\mathbf{1}_{[0,\tau^{\hat{X}})}(s)ds,\quad t\in[0,T].&&
\nonumber 
\ea
\noindent as well as the one associated to
\ba 
M^{g}_t(\varphi,q)\doteq g(\varphi(t))-g(\varphi(0))-\int_{[0,t]\times\Gamma}\bar{b}(s,\varphi(s),\mu_s,u)^\top D g(\varphi(s))q(ds,du)&&\nonumber\\
+\frac{1}{2}\int_0^t \text{Tr}\left[{\sigma \sigma^\top} D^2 g(\varphi(s))\right] ds&&
\nonumber 
\ea
\noindent up to time $\tau^{\hat{X}}\wedge T$, i.e. the martingale property is satisfied by the processes above stopped at time $\tau^{\hat{X}}\wedge T$.
Second, $\Theta''$ solves the latter martingale problem up to time $T$. Then $\Theta'$ and $\Theta''$ solve the same martingale problem up to time $\tau^{\hat{X}}\wedge T$. Moreover, we have $\Theta''(q\in\mathcal{V}:q_t=\lambda(t,t\wedge\tau^{\hat{X}},\hat{X}_t))=1$ for $\mathcal L_T$-a.e. $t\in[0,T]$. If we set $\Theta_t\doteq\Theta\circ(\hat{X},\hat{\Lambda})_{\cdot\wedge t}^{-1}$ for all $\Theta\in\PP(\X\times\V)$ and $t\in[0,T]$, then by uniqueness of the solution of the martingale problem up to time $\tau^{\hat{X}}\wedge T$ we have 
\[ \Theta'_t (\cdot \cap \{t\le \tau^{\hat{X}} \wedge T\}) = \Theta''_t (\cdot \cap \{t \le \tau^{\hat{X}} \wedge T\}).\]
\noindent Hence $J^{\mu}(\Theta')=J^{\mu}(\Theta'')$. Now $\Theta''$ satisfies item (ii) of Definition \ref{def:MFGsolRel2}.

To conclude notice that the process $Y^{\tau^{X''}}_t=(t\wedge\tau^{X''},X''_t)$ reduces to $(t,X''_t)$ before time $\tau^{X''}\wedge T$. Hence, also $\lambda_{t,Y_t^{\tau^{X''}}}$, with a slight abuse of notation, reduces to $\lambda_{t,X''_t}$. With the additional Assumption \hyperlink{C2}{(\textrm{C2})}, the second part of this lemma follows from the proof of Proposition \ref{prop:existenceStrictFeedlSolMFG} applied to the stopped process $Y^{\tau^X}$.
\end{proof}

\section{Uniqueness of solutions of the mean-field game}\label{sec:uniqueness}
\noindent In this section we address the problem of uniqueness of MFG solutions. Precisely, under Assumptions \hyperlink{H1}{(\textrm{H1})}-\hyperlink{H8}{(\textrm{H8})} and with the additional Assumptions \hyperlink{U1}{(\textrm{U1})}-\hyperlink{U2}{(\textrm{U4})} given below, where the second one guarantees monotonicity of the running cost in the same spirit as \cite{lasry-lions2007} (see also Theorem 3.29 in \cite{carmonadelarue}), we show uniqueness of the MFG solution also in the presence of smooth dependence on past absorptions. The extra assumptions can be formulated as follows.
\begin{itemize}
\item[\hypertarget{U1}{(\textrm{U1})}] The running cost can be split in two terms:
\ba 
\bar{f}(t,x,\mu,u)=\bar{f}_0(t,x,u)+\bar{f}_1(t,x,\mu)
\nonumber
\ea
\noindent for some measurable functions $\bar{f}_0:[0,T]\times\R^d\times\Gamma\rightarrow[0,\infty)$ and $\bar{f}_1:[0,T]\times\Rd\times\mathcal{M}_{\leq 1,1}(\R^d)\rightarrow[0,\infty)$.
\item[\hypertarget{U2}{(\textrm{U2})}] Lasry-Lions monotonicity assumption: Let $\mu,\tilde{\mu}\in\mathcal{M}_{\leq 1,1}(\R^d)$, $\mu\neq\tilde{\mu}$. Then
\ba 
\int_{\Rd}\round{\bar{f}_1(t,x,\mu)-\bar{f}_1(t,x, \tilde{\mu})}(\mu-\tilde{\mu})(dx)\geq 0,\quad t\in[0,T].
\nonumber
\label{as:monotonicity}
\ea
\item[\hypertarget{U3}{(\textrm{U3})}] The drift $b$ does not depend on the measure variable.
\item[\hypertarget{U4}{(\textrm{U4})}] Let $\bar{\mu}\in\Upsilon^T_{\leq 1,1}$ be fixed. Then the following optimization problem
\ba
\inf_{\Lambda\in\tilde{\mathbb{U}}}J^{\bar\mu}\round{\Lambda} \doteq \E\Biggl[\int_{[0,\tau]\times\Gamma}\bar{f}\round{s, X_s,\bar{\mu}_s, u }\Lambda_s(du)ds+F\round{\tau,X_{\tau}}\Biggr]
\label{eq:opt}
\ea
\noindent has a unique solution $\Lambda^{\bar{\mu}}$, where $(\Omega, \F, (\F_t)_{t\in[0,T]},\PB, W, X)$ is a solution of Eq.\eqref{eq:SDEopen} under $\Lambda^{\bar{\mu}}$ with initial distribution $\nu$ and drift $b$ satisfying \hyperlink{U3}{(\textrm{U3})}.
\end{itemize}
\bteo[\textit{Uniqueness}]\label{teo:uniqueness} Under Assumptions \hyperlink{H1}{(\textrm{H1})}-\hyperlink{H8}{(\textrm{H8})} and \hyperlink{U1}{(\textrm{U1})}-\hyperlink{U4}{(\textrm{U4})}, if there exists a feedback solution of the MFG $(\lambda,\mu)$ (as in Definition \ref{def:MFGsolRel}) then it is unique. 
\eteo
\begin{proof}
\noindent By contradiction, let $(\lambda,\mu)$ and $(\tilde{\lambda},\tilde{\mu})$ be two different feedback MFG solutions (as in Definition \ref{def:MFGsolRel}). Then
\ba 
J^{\tilde{\mu}}(\lambda)-J^{\tilde{\mu}}(\tilde{\lambda})> 0\quad\text{and}\quad J^{\mu}(\tilde{\lambda})-J^{\mu}(\lambda)> 0
\nonumber
\ea
\noindent where the inequality is strict by uniqueness of the minimizer in Assumption \hyperlink{U2}{(\textrm{U4})}, and in particular
\ba 
\Delta(\mu,\tilde{\mu},\lambda,\tilde{\lambda})\doteq J^{\tilde{\mu}}(\lambda)-J^{\tilde{\mu}}(\tilde{\lambda})+ J^{\mu}(\tilde{\lambda})-J^{\mu}(\lambda)> 0.
\nonumber
\ea
\noindent However, thanks to Assumption \hyperlink{U3}{(\textrm{U3})} that grants independence of the dynamics of the state processes from the flows of measures $\mu$ and $\tilde{\mu}$
\ba 
\Delta(\mu,\tilde{\mu},\lambda,\tilde{\lambda})&=&\E^{\PB}\brackets{\int_0^T\mathbf{1}_{[0,\tau)}(t)\round{\bar{f}_1(t,X_t,\tilde{\mu}_t)-\bar{f}_1(t,X_t,\mu_t)}dt}\nonumber\\
&&+\E^{\tilde{\PB}}\brackets{\int_0^T\mathbf{1}_{[0,\tilde{\tau})}(t)\round{\bar{f}_1(t,\tilde{X}_t,\mu_t)-\bar{f}_1(t,\tilde{X}_t,\tilde{\mu}_t)}dt}
\nonumber
\ea
\noindent where $(\Omega, \F, (\F_t)_{t\in[0,T]},\PB, W, X)$ and $(\tilde{\Omega}, \tilde{\F}, (\tilde{\F}_t)_{t\in[0,T]},\tilde{\PB}, \tilde{W}, \tilde{X})$ are weak solutions of Eq.\eqref{eq:MFGRel} respectively with controls $\lambda$ and $\tilde{\lambda}$. Set $\theta\doteq\PB\circ X^{-1}$ and $\tilde{\theta}\doteq\tilde{\PB}\circ \tilde{X}^{-1}$. Then
\ba 
\Delta(\mu,\tilde{\mu},\lambda,\tilde{\lambda})&=&\int_{\X}\int_0^T\mathbf{1}_{[0,\tau(\varphi))}(t)\brackets{\bar{f}_1(t,\varphi(t),\mu_t)-\bar{f}_1(t,\varphi(t),\tilde{\mu}_t)}dt\tilde{\theta}(d\varphi)\nonumber\\
&&-\int_{\X}\int_0^T\mathbf{1}_{[0,\tau(\varphi))}(t)\brackets{\bar{f}_1(t,\varphi(t),\mu_t)-\bar{f}_1(t,\varphi(t),\tilde{\mu}_t)}dt\theta(d\varphi)\nonumber\\
&=&\int_0^T\int_{\X}\brackets{\bar{f}_1(t,\varphi(t),\mu_t)-\bar{f}_1(t,\varphi(t),\tilde{\mu}_t)}\mathbf{1}_{[0,\tau(\varphi))}(t)\tilde{\theta}(d\varphi)dt\nonumber\\
&&-\int_0^T\int_{\X}\brackets{\bar{f}_1(t,\varphi(t),\mu_t)-\bar{f}_1(t,\varphi(t),\tilde{\mu}_t)}\mathbf{1}_{[0,\tau(\varphi))}(t)\theta(d\varphi)dt\nonumber\\
&=&\int_0^T\int_{\Rd}\brackets{\bar{f}_1(t,x,\mu_t)-\bar{f}_1(t,x,\tilde{\mu}_t)}\tilde{\mu}_t(dx)dt\nonumber\\
&&-\int_0^T\int_{\Rd}\brackets{\bar{f}_1(t,x,\mu_t)-\bar{f}_1(t,x,\tilde{\mu}_t)}\mu_t(dx)dt\nonumber\\
&=&-\int_0^T\int_{\Rd}\brackets{\bar{f}_1(t,x,\mu_t)-\bar{f}_1(t,x,\tilde{\mu}_t)}(\mu_t-\tilde{\mu}_t)(dx)dt 
\nonumber
\ea
\noindent which is lower than or equal to zero by Assumption \hyperlink{U2}{(\textrm{U2})}. In the second equality we have used Fubini-Tonelli theorem, while the third one comes from the definitions of $\mu$ and $\tilde{\mu}$, i.e.
\ba 
\mu_t(B)&\doteq &\theta\round{\{X_t\in B\}\cap\{t<\tau\}}\nonumber\\
&=&\int_{\X}\mathbf{1}_B(\varphi(t))\mathbf{1}_{[0,\tau(\varphi))}(t)\theta(d\varphi)\nonumber\\
&=&\int_{\Rd}\mathbf{1}_B(x)\mu_t(dx),\quad t\in[0,T]
\nonumber
\ea
\noindent for all $B\in\mathcal{B}(\Rd)$ and similarly for $\tilde{\mu}$.
\end{proof}
\begin{example}[Non-local dependence on the measure through a weighted average] We provide and example of running cost $\bar{f}$ satisfying the monotonicity condition \hyperlink{U2}{(\textrm{U2})}, which is an assumption on the measure-dependent term $\bar{f}_1$ only.
Let $w: \R^d \to [0,\infty)$ be some measurable function with sub-linear growth so that 
\[ m_w(\mu)\doteq \int_{\R^d} w(x)\mu (dx) < \infty, \quad \textrm{for all }\mu \in \mathcal{M}_{\le 1,1}(\R^d)\]  and set
\ba 
\bar{f}_1(t,x,\mu)\doteq w(x)\int_{\Rd}w(y)\mu(dy)=w(x)m_w(\mu),\quad (t,x,\mu)\in[0,T]\times\Rd\times\mathcal{M}_{\leq 1,1}(\R^d).
\nonumber
\ea
\noindent Since
\ba 
\bar{f}_1(t,x,\mu)-\bar{f}_1(t,x,\tilde{\mu})= w(x)\int_{\Rd}w(y)(\mu-\tilde{\mu})(dy)\nonumber
\ea
we obtain
\ba
\int_{\Rd}\round{\bar{f}_1(t,x,\mu)-\bar{f}_1(t,x,\tilde{\mu})}(\mu-\tilde{\mu})(dx)&=&\int_{\Rd}w(x)\int_{\Rd}w(y)(\mu-\tilde{\mu})(dy)(\mu-\tilde{\mu})(dx),\nonumber\\
&=&\int_{\Rd}w(x)(\mu-\tilde{\mu})(dx)\int_{\Rd}w(y)(\mu-\tilde{\mu})(dy),\nonumber\\
&=&\round{\int_{\Rd}w(x)(\mu-\tilde{\mu})(dx)}^2\geq 0.
\nonumber
\ea
\end{example}

\section{Approximate Nash equilibria for the $N$-player game with finite-dimensional interaction}\label{sec:Nplayer}

\noindent In this section, we consider an important particular case of our MFG with absorption, where the mean-field interaction is finite-dimensional. This is inspired by the original model of \cite{campi2018n}. We show that any feedback solution of the MFG can be used to construct a sequence of approximate Nash equilibria for the corresponding $N$-player game.  To this end, we will need two additional assumptions (Assumptions \hyperlink{N1}{(\textrm{N1})} and \hyperlink{N2}{(\textrm{N2})} below). We focus on a finite-dimensional example first for technical reasons: this setting is very suitable to the propagation of chaos result that we use in the proofs without being too technical. Second, we think that this case is also particularly relevant for the applications as mentioned in the introduction. Overall, we believe that the finite-dimensional setting enables us to keep a good balance between abstract technicalities and modelling needs. 

The approximation result is the content of Theorem \ref{teo:approxNashRel}  and Corollary \ref{cor:approxNash}. In order to prove this, we interpret the $N$-player system as a system of $N$ interacting diffusions (as in, e.g., \cite{mckean,sznitman,gartner}). While the usual mode of convergence of an $N$-particle system is the convergence in law of the empirical measures, here we obtain a stronger form of propagation of chaos as in  \cite{lacker2018strong} but with possibly unbounded drift in the state variable. We prove that the empirical measures converge in the stronger $\tau$-topology, which is widely used in the large deviations literature (see, for instance, Chapter 6.2 in \citet{dembo}); see Subsection \ref{subsec:chaos}.\\
\medskip

\subsection{The setting with finite-dimensional interaction}
\noindent Here, we describe the MFG and the corresponding $N$-player game with smooth dependence on past absorptions, specializing them to the finite-dimensional interaction setting. In particular, we give the definition of $\epsilon$-Nash equilibrium for the $N$-player game. Then, we give the assumptions that are specific to this model. We conclude by checking that the MFG with finite-dimensional interactions satisfies the hypotheses of Theorem \ref{teo:existenceRelFeedSolMFG}, granting the existence of relaxed and strict solutions of the MFG.\\
\medskip

\noindent\textit{The mean-field dynamics.}\quad Given a feedback control $u \in \mathcal U_{fb}$ and a flow of sub-probability measures $\mu \in \Upsilon^T_{\leq 1,1}$, the representative player's state evolves according to the equation\\
\ba 
X_t = X_0 + \int_{0}^{t}\tilde{b}\round{s, X_s,L\round{\mu_s}, m_w\round{\mu_s}, u\round{s,X} }\,ds + \sigma W_t,\quad t \in \brackets{0,T}
\label{eq:MFGfindim}
\ea
\noindent where $X$ is a $d$-dimensional stochastic process starting at $X_0 \overset{d}{\sim} \nu \in \PP(\Rd)$, $W$ is a $d$-dimensional Wiener process on some filtered probability space $(\Omega, \F,(\F_t)_{t\in[0,T]},\PB)$, $\tilde{b}$ and $\sigma$ are as in the assumptions below. In addition, $m_w\round{\mu}$ and $L\round{\mu}$ are functions $m_w:\mathcal{M}_{\leq 1,1}(\Rd) \rightarrow \R^{d_0}$ and $L:\mathcal{M}_{\leq 1,1}(\Rd)\rightarrow [0,1]$ defined as
\ba 
m_w\round{\mu}\doteq \int_{\Rd}w\round{x}\mu(dx)\quad\text{and}\quad
L\round{\mu}\doteq 1-\int_{\Rd}\mu(dx)
\nonumber
\ea
where $w:\R^d\rightarrow\R^{d_0}$, $d_0\in\N$, is a fixed weight function with sub-linear growth. Again, solutions of Eq.\eqref{eq:MFGfindim} are understood in the weak sense (see Remark \ref{rem:existenceUniqWeak}). The cost associated to a strategy $u \in \mathcal{U}_{fb}$ and a flow of sub-probability measures $\mu\in\Upsilon^T_{\leq 1,1}$ is given by
\ba
J^{\mu}\round{u} \doteq \E\Biggl[\int_{0}^{\tau}\tilde{f}\round{s, X_s, L\round{\mu_s},m_w\round{\mu_s}, u\round{s,X} }ds+F\round{\tau,X_{\tau}}\Biggr]
\label{eq:MFGcostFindim}
\ea
\noindent where 
$\tau \doteq \tau^{X} \wedge T$ is the random time horizon as in the previous sections.\\
\medskip

\noindent\textit{The $N$-player dynamics.}\quad Let $N \in \mathbb{N}$ be the number of players. We assume that the players' private states evolve according to the following system of $N$ $d$-dimensional SDEs: for $i \in \braces{1,\ldots ,N}$,
\ba 
X\Ni_t = X\Ni_0 + \int_{0}^{t} \tilde{b}\round{s, X\Ni_s, L\round{\mu^N_s},m_w\round{\mu^N_s}, u\Ni\round{s,\textbf{X}^N}}\,ds + \sigma W\Ni_t
\label{eq:Ngame}
\ea
\noindent for $t\in[0,T]$, where $X^{N,i}_0\simdistr\nu$ i.i.d.,  $W^{N,1}, \ldots,W^{N,N}$ is an $N$-dimensional vector of independent $d$-dimensional Wiener processes,  $\textbf{X}^N$ denotes the vector of all players' private states, $\bu^N$ the vector of feedback strategies, $\tilde{b}$ and $\sigma$ are as in the assumptions below. We remind that $\mu^N\in\Upsilon^T_{\leq 1,1}$ is the random empirical \emph{sub-probability measures} defined as
\ba
\mu_t^N\round{\cdot} \doteq \frac{1}{N}\sum_{i = 1}^{N} \delta_{X_t^{N,i}}\round{\cdot}\mathbf{1}_{[0,\tau^{X^{N,i}})}\round{t},\quad t\in[0,T].
\label{eq:subProbN}
\ea
\noindent Solutions of the SDEs in Eq.\eqref{eq:Ngame} are understood to be in the weak sense on some filtered probability space $(\Omega^N,\F^N,(\F^N_t)_{t\in[0,T]},\PB^N)$ satisfying the usual conditions (see Remark \ref{rem:existenceUniqWeak}). \\
\indent Let $\mathcal{U}_{1}^{N}$ be the set of all progressively measurable functionals $u:[0,T]\times \mathcal{X}^N\rightarrow\Gamma$, and let $\mathcal{U}_{N}^{N}$, the set of all vectors $\bu^N$ such that $u^{N,i} \in \mathcal{U}_{1}^{N}$, $i \in \braces{1, \ldots, N}$. Each element of $\mathcal{U}_{N}^{N}$ is called \textit{feedback strategy vector}.
In this game, player $i$ evaluates a strategy vector $\bu^N\in\mathcal{U}^N_N$ according to his/her expected costs
\ba
J^{N,i}\round{\bu^N} \doteq \E\Biggl[\int_{0}^{\tau\Ni}\tilde{f}\round{s,X\Ni_s,L\round{\mu^N_s},m_w\round{\mu^N_s},u\Ni(s,\textbf{X}^N)}ds \nonumber\\
+F\round{\tau\Ni,X\Ni_{\tau\Ni}}\Biggr]
\label{eq:Ncost}
\ea
\noindent over a random time horizon, where $\textbf{X}^N$ is the $N$-player dynamics under $\bu^N$ and $\tau^{N,i} \doteq \tau^{X^{N,i}} \wedge T$. Our aim is the construction of approximate Nash equilibria for the $N$-player game from a solution of the limit problem. In the next definition, we use the standard notation $[u^{N,-i},v]$ to indicate a strategy vector equal to $\bu^N$ for all players but the $i$-th, who deviates by playing $v\in\mathcal{U}^N_1$ instead. 
\bdefi[\textit{$\epsilon$-Nash equilibrium}] Let $\epsilon \geq 0$.  A strategy vector $\bu^N\in\mathcal{U}_N^N$ is called $\epsilon$-Nash equilibrium for the $N$-player game if for every $i\in\{1,\ldots,N\}$ and for any deviation $v\in\mathcal{U}^N_1$ we have:
\ba 
J\Ni(\bu^N)\leq J\Ni\round{\brackets{u^{N,-i},v}}+\epsilon.
\nonumber
\label{eq:epsilonNash}
\ea
\edefi
\medskip

\noindent\textit{Relaxed controls.}\quad It will be very convenient to use relaxed controls also in the $N$-player case. Let $\widetilde{\mathcal{U}}^N_1$ be the set of all single-player relaxed strategies for the $N$-player game, and let $\widetilde{\mathcal{U}}^N_N$ be the set of $N$-player relaxed strategy vectors, i.e. vectors $\boldsymbol{\lambda}^N = (\lambda^{N,1}, \ldots, \lambda^{N,N})$ with $\lambda^{N,i} \in \widetilde{\mathcal{U}}^N_1$, $i\in\{1,\ldots,N\}$.
At this point, we can rewrite the dynamics and the cost functional of the $N$-player game (Eq.\eqref{eq:Ngame} and Eq.\eqref{eq:Ncost}) by using relaxed controls as
\ba 
X\Ni_t = X\Ni_0 +\! \int_{[0,t]\times\Gamma}\tilde{b}\round{s, X\Ni_s,L\round{\mu^N_s},m_w\round{\mu^N_s},u}\lambda^{N,i}\round{s,\textbf{X}^{N}}(du)ds + \! \sigma W\Ni_t 
\label{eq:NgameRel}
\ea
\noindent with associated cost
\ba
J^{N,i}\round{\boldsymbol{\lambda}^N} =\E\Biggl[\int_{[0,\tau\Ni]\times\Gamma}\tilde{f}\round{s, X\Ni_s,L\round{\mu^N_s},m_w\round{\mu^N_s},u}\lambda^{N,i}\round{s,\textbf{X}^{N}}(du)ds\nonumber\\
+F\round{\tau\Ni,X\Ni_{\tau\Ni}}\Biggr]
\label{eq:costNgameRel}
\ea
\noindent for $t\in[0,T]$, $i\in\{1,\ldots,N\}$, $\boldsymbol{\lambda}^N\in\widetilde{\mathcal{U}}^N_N$ and $\lambda^{N,i}\in\widetilde{\mathcal{U}}^N_1$ for all $i\in\{1,\ldots,N\}$.
Moreover, we extend accordingly the notion of $\epsilon$-Nash equilibrium.
\bdefi[\textit{Relaxed $\epsilon$-Nash equilibrium}]\label{def:epsilonNashRel} A strategy vector $\boldsymbol{\lambda}^N\in\widetilde{\mathcal{U}}_N^N$ is an $\epsilon$-Nash equilibrium for the $N$-player game if for every $i\in\{1,\ldots ,N\}$ and for any single-player strategy $\beta\in\widetilde{\mathcal{U}}^N_1$
\ba 
J\Ni(\boldsymbol{\lambda}^N)\leq J\Ni\round{\brackets{\boldsymbol{\lambda}^{N,-i},\beta}}+\epsilon.
\nonumber
\label{eq:epsilonNashRel}
\ea
\edefi
\noindent The drift $\tilde{b}$, the function $w$, the running cost $\tilde{f}$ and the terminal cost $F$ now satisfy the following assumptions, replacing Assumptions \hyperlink{H1}{(\textrm{H1})}-\hyperlink{H3}{(\textrm{H3})}:
\bi 
\item[\hypertarget{H1'}{\textrm{(H1')}}] The drift $\tilde{b}:[0,T]\times\R^d\times[0,1]\times\R^{d_0}\times\Gamma\rightarrow \R^d$ is jointly continuous and satisfies the following uniform Lipschitz continuity: there exists $L>0$ such that
\ba 
\modl{\tilde{b}\rbr{t,x,\ell,m,u}-\tilde{b}\rbr{t,x',\ell,m,u}}\leq L\modl{x-x'}
\nonumber
\label{eq:driftLipFindim}
\ea
\noindent for all $x,x'\in \R^d$ and all $(t,\ell,m,u)\in[0,T]\times[0,1]\times\R^{d_0}\times\Gamma$. Moreover it has sub-linear growth in $(x,m)$ uniformly in the other variables, i.e. there exists a constant $C>0$ such that
\ba 
\modl{\tilde{b}\rbr{t,x,\ell,m,u}}\leq C\rbr{1+\modl{x}+\modl{m}}
\nonumber
\label{eq:driftLinFindim}
\ea
\noindent for all $(t,x,\ell,m, u)\in[0,T]\times\Rd\times[0,1]\times\R^{d_0}\times \Gamma$.
\item[\hypertarget{H2'}{\textrm{(H2')}}] $w:\R^d\rightarrow \R^{d_0}$ is continuous and has sub-linear growth: $|w(x)| \leq C(1+|x|)$ for all $x \in \R^d$.
\item[\hypertarget{H3'}{\textrm{(H3')}}] The costs $\tilde{f}:\brackets{0,T} \times \R^{d}\times[0,1] \times  \R^{d_0} \times \Gamma \rightarrow [0,\infty)$ and $F:\brackets{0,T} \times  \R^{d} \rightarrow [0,\infty)$ are jointly continuous. Moreover, they have sub-linear growth:
\ba
\modl{\tilde{f}(t,x,\ell,m,u)}&\leq & C\rbr{1+\modl{x}+\modl{m}},\nonumber \\
\modl{F(t,x)}&\leq & C\rbr{1+\modl{x}},
\nonumber
\ea
\noindent for all $(t,x,\ell,m,u)\in[0,T]\times\Rd\times[0,1]\times\R^{d_0}\times\Gamma$.
\ei
\noindent We conclude the presentation of the finite-dimensional model by introducing the coefficients' reparametrization on $\PP_1(\X)$, by checking their joint continuity (as in Assumption \hyperlink{H3}{(H3)}), where continuity in the measure variable is in the 1-Wasserstein distance and at points $\theta\ll\wiener$.
We set $(\bar{b},\bar{f})(t,x,\mu,u)\doteq (\tilde{b},\tilde{f})(t,\varphi(t),L(\mu),m_w(\mu),u)$ for all $(t,x,\mu,u)\in[0,T]\times\Rd\times\M_{\le 1,1}(\Rd)\times\Gamma$ and define the reparametrization $(b,f)$ as in Section \ref{sec:preliminaries}. Then
\begin{eqnarray}
(b,f)(t,\varphi,\theta,u)&= &(\tilde{b},\tilde{f})(t,\varphi(t),L(t;\theta),m_w(t;\theta),u)
\nonumber
\end{eqnarray}
\noindent where
\ba 
m_w(t;\theta)&\doteq &\int_{\X}w\round{\varphi(t)}\mathbf{1}_{[0,\tau(\varphi))}(t)\theta(d\varphi),\nonumber\\
L(t;\theta)&\doteq &1-\int_{\X}\mathbf{1}_{[0,\tau(\varphi))}(t)\theta(d\varphi)
\nonumber
\ea
\noindent are called the average and loss process and they equal $m_w(\mu_t)$ and $L(\mu_t)$ in case $\mu_t =g(t,\theta)$ where $g$ is defined as in Eq.\eqref{eq:defGtheta}.\\
\indent Joint continuity of $b$ and $f$ follows from joint continuity of $\tilde{b}$ and $\tilde{f}$ and from the following lemma.
\blem[\textit{Continuity of the average and loss processes}] \label{lem:cont-l-m} Grant Assumptions \hyperlink{H1'}{(\textrm{H1'})}-\hyperlink{H3'}{(\textrm{H3'})} and \hyperlink{H4}{(\textrm{H4})}-\hyperlink{H8}{(\textrm{H8})}. Let $(\theta_n)_{n\in\N}\subset\PP_1(\X)$ converge to $\theta\in\PP_1(\X)$, $\theta\ll\wiener$, in the 1-Wasserstein distance, then
\bi 
\item[(i)]\quad $L(t;\theta^n)\rightarrow L(t;\theta)$ as $n\rightarrow\infty$.
\item[(ii)]\quad $m_w(t;\theta^n)\rightarrow m_w(t;\theta)$ as $n\rightarrow\infty$.
\ei
\elem
\begin{proof}\noindent\textit{(i).}\quad Denote by $\mathbb D_{\tau}(t)$ the set of discontinuity points of the map $\varphi\mapsto\mathbf{1}_{[0,\tau(\varphi))}(t)$ for $t\in[0,T]$. In particular $\theta^n\weakconv\theta$. Then:
\ba 
L(t;\theta^n)-L(t;\theta)=-\int_{\X}\mathbf{1}_{[0,\tau(\varphi))}(t)\round{\theta^n-\theta}(d\varphi)\conv 0
\nonumber
\ea
\noindent for all $t\in[0,T]$. This follows from the definition of weak convergence of measures, the fact that $\theta(\mathbb D_{\tau}(t))=0$ for all $t\in[0,T]$ (due to $\theta\ll\wiener$) and by Lemma \ref{lem:regularity}.(d).\\
\indent\textit{(ii).}\quad Now we have:
\ba 
\modl{m_w(t;\theta^n)-m_w(t;\theta)}&\leq &\modl{\int_{\X}w(\varphi(t))\mathbf{1}_{[0,\tau(\varphi))}(t)\round{\theta^n-\theta}(d\varphi)}\conv 0
\nonumber
\ea
\noindent for all $t\in[0,T]$ as a consequence of the convergence in the 1-Wasserstein distance, the fact that $\theta(\mathbb D_{\tau}(t))=0$ for all $t\in[0,T]$ and by Lemma \ref{lem:regularity}.(d) together with Lemma \ref{lem:regularityConvergence}. 
\end{proof}
\noindent We conclude by proving that we can use Theorem \ref{teo:existenceRelFeedSolMFG} and get existence of a feedback relaxed and strict solutions of the MFG with smooth dependence on past absorptions and finite-dimensional dependence on the measure.
\bcor[\textit{Existence of relaxed and strict feedback MFG solutions}]\label{teo:existenceSolFinite} Under Assumptions \hyperlink{H1'}{(\textrm{H1'})}-\hyperlink{H3'}{(\textrm{H3'})}, \hyperlink{H4}{(\textrm{H4})}-\hyperlink{H8}{(\textrm{H8})} and \hyperlink{C1}{(\textrm{C1})} , there exists a relaxed feedback solution $(\lambda,\mu)$ of the MFG with finite dimensional interaction. Moreover, under the additional Assumption \hyperlink{C2}{(\textrm{C2})} , there exists a strict feedback MFG solution $(u,\mu)$.
\ecor
\begin{proof}
\noindent Assumptions \hyperlink{H1'}{(\textrm{H1'})}-\hyperlink{H3'}{(\textrm{H3'})} imply Assumptions  \hyperlink{H1}{(\textrm{H1})}-\hyperlink{H3}{(\textrm{H3})} of Theorem \ref{teo:existenceRelFeedSolMFG}. Indeed,  \hyperlink{H1}{(\textrm{H1})}-\hyperlink{H2}{(\textrm{H2})} follow from the definition of the coefficients $\tilde{b}$ and $\tilde{f}$. Assumption \hyperlink{H3}{(\textrm{H3})}, i.e. joint continuity of the reparametrized coefficients, is a consequence of joint continuity of $\tilde{b}$ and $\tilde{f}$ and Lemma \ref{lem:cont-l-m}.
\end{proof}

\subsection{The $N$-player approximation theorem}\label{subsec:Nplayer}
\noindent In order to state the $N$-player approximation results, we need the following two additional assumptions \hyperlink{N1}{(\textrm{N1})}-\hyperlink{N2}{(\textrm{N2})}, whose formulation requires some more terminology.\\
\indent We set
\ba 
d_t^{TV}(\theta, \tilde{\theta}) \doteq \text{sup}_{B \in \mathcal{F}_t}|\theta(B)-\tilde{\theta}(B)|,
\nonumber
\ea
\noindent for all $\theta,\,\tilde{\theta} \in \mathcal{P}(\mathcal{X})$ and we note that for $t \in [0,T)$, $d_t$ is only a pseudo-metric, whereas for $t = T$ it is a proper metric; $d_{T}^{TV}$ is called the total variation distance. However, with a slight abuse of terminology, we will often refer to $d^{TV}_t$ as the total variation distance for each $t \in \brackets{0,T}$.
\begin{itemize}
\item[\hypertarget{N1}{\textrm{(N1)}}] The function $w:\Rd\rightarrow\R^{d_0}$ is bounded.
\item[\hypertarget{N2}{\textrm{(N2)}}] The drift $\tilde{b}$ satisfies the following Lipschitz continuity:
\ba 
\modl{\tilde{b}\rbr{t,x,\ell,m,u}-\tilde{b}\rbr{t,x',\ell',m',u}}\leq L\left(\modl{x-x'} +\modl{\ell-\ell'}+\modl{m-m'}\right)
\nonumber
\label{eq:driftLipFindim2}
\ea
\noindent for all $(x,\ell,m),(x',\ell',m')\in \R^d\times[0,1]\times\R^{d_0}$ and all $(t,u)\in[0,T]\times\Gamma$, with Lipschitz constant $L>0$.
The running cost $\tilde{f}$ can be decomposed as
\ba
\tilde{f}(t, x,\ell,m, u) = \tilde{f}_0(t,x, u) + \tilde{f}_1(t, x,\ell,m),
\nonumber
\ea
\noindent where
\ba 
|\tilde{f}_0(t,x, u)| \leq K \quad \text{and} \quad | \tilde{f}_1(t, x,\ell,m)| \leq C (1 + |x|),
\nonumber
\ea
\noindent for all $(t, x,\ell,m,u) \in [0,T] \times\Rd\times[0,1]\times\R^{d_0} \times \Gamma $ and some constants $C,K>0$.
\end{itemize}
\noindent From Assumptions \hyperlink{N1}{(\textrm{N1})}-\hyperlink{N2}{(\textrm{N2})}, the reparametrizations $b$ and $f$ inherit a series of properties that are fundamental in the proof of the approximation result. First, being $w:\R^d\rightarrow\R^{d_0}$ bounded, the drift $b$ is Lipschitz continuous with respect to the total variation distance, which is a key assumption in Lemma \ref{lem:mcKeanEU}. Indeed
\ba 
\modl{b(t, \varphi, \theta, u)-b(t, \varphi, \theta', u)}&\leq &
L\round{\modl{L(t;\theta)-L(t;\theta')}+\modl{m_w(t;\theta)-m_w(t;\theta')}}\nonumber\\
&\leq &L (1+\infnorm{w})d^{TV}_T(\theta,\theta')\doteq L^{TV}_b d^{TV}_T(\theta,\theta')
\nonumber
\ea
because
\ba 
\modl{L(t;\theta)-L(t;\theta')}&=& \modl{\int_{\X}\mathbf{1}_{[0,\tau(\varphi))}(t)(\theta'-\theta)(d\varphi)}\leq d^{TV}_T(\theta,\theta')\quad\text{and}\nonumber\\
\modl{m_w(t;\theta)-m_w(t;\theta')}&=& \modl{\int_{\X}w(\varphi(t))\mathbf{1}_{[0,\tau(\varphi))}(t)(\theta-\theta')(d\varphi)}\leq \infnorm{w}d^{TV}_T(\theta,\theta').
\nonumber
\ea
\noindent Second, the sub-linear growth property
\ba
|b(t, \varphi, \theta, u)| \leq C(1+\infnorm{w}+\|\varphi\|_{\infty, t}),\quad(t, \varphi) \in [0,T] \times \mathcal{X}
\nonumber
\ea
is uniform in $\theta\in \mathcal{P}(\mathcal{X})$ and in $u \in \Gamma$, implying that $b$ is bounded in the measure and control variables (and analogously $f$). This means that $b$ and $f$ are well defined on all $\PP(\X)$ not only on $\PP_1(\X)$, which is fundamental to apply the fixed point theorem in Lemma \ref{lem:mcKeanEU}. Finally, the running cost $f$ can be decomposed as
\ba
f(t, \varphi, \theta, u) = f_0(t, \varphi, u) + f_1(t, \varphi, \theta)
\nonumber
\ea
\noindent where its components are
\ba
f_0(t, \varphi, u)\doteq \tilde{f}_0(t,\varphi(t),u)\quad\text{and}\quad
f_1(t, \varphi, \theta)\doteq \tilde{f}_1(t,\varphi(t),L(t;\theta),m_w(t;\theta))
\nonumber
\ea
which inherit from $\tilde f_0$ and $\tilde f_1$ the properties
\ba 
|f_0(t, \varphi, u)| \leq K \quad \text{and} \quad |f_1(t, \varphi, \theta)| \leq C (1 + \|\varphi\|_{\infty,t})
\nonumber
\ea
for all $(t, \varphi, \theta, u) \in [0,T] \times \mathcal{X} \times \mathcal{P}(\mathcal{X}) \times \Gamma $. This is a key assumption to perform the passage to the many-player limit in Theorem \ref{teo:approxNashRel}. Indeed, boundedness in the control of $f_0$ enables us to exploit convergence in the $\tau$-topology while sub-linearity in the state variable $\varphi$ uniformly in the measure variable $\theta$ makes $f_1$ a good test function for the convergence in the 1-Wasserstein distance.

\bteo[\textit{Approximate Nash equilibria - relaxed}] Let $(\lambda,\mu)$ be a relaxed feedback MFG solution. For all $N \ge 2$, define $\blambda^N = (\lambda^{N,1}, \ldots, \lambda^{N,N}) \in\UTNN$  where $\lambda\Ni(t,\varphi^N)\doteq \lambda(t,\varphi\Ni)$ for all $i \in \{1,\ldots,N\}$, $t\in[0,T]$ and $\varphi^N\in\X^N$.\\
Then under Assumptions \hyperlink{H1'}{(H1')}-\hyperlink{H3'}{(H3')}, \hyperlink{H4}{(H4)}-\hyperlink{H8}{(H8)} and \hyperlink{N1}{(N1)}-\hyperlink{N2}{(N2)}, for every $\epsilon>0$ there exists $N^{\epsilon}\in\mathbb{N}$ such that $\blambda^N$ is an $\epsilon$-Nash equilibrium for the $N$-player game whenever $N\geq N^{\epsilon}$, i.e. for every $i \in \{1,\ldots,N\}$ and for any deviation $\beta\in\UTNone$
\ba 
J\Ni\round{\blambda^{N}}\leq J\Ni\round{\brackets{\lambda^{N,-i},\beta}}+\epsilon 
\nonumber
\label{eq:approxNashRel}
\ea
\noindent for all $N\geq N^{\epsilon}$.
\label{teo:approxNashRel}
\eteo

\bcor[\textit{Approximate Nash equilibria - strict}] Let $(u,\mu)$ be a strict feedback MFG solution. For all $N \ge 2$, define $\bu^N = (u^{N,1},\ldots, u^{N,N})\in\UNN$ where $u\Ni(t,\varphi^N)\doteq u(t,\varphi\Ni)$ for all $i\in \{1,\ldots,N\}$, $t\in[0,T]$ and $\varphi^N\in\X^N$.\\
Then under Assumptions \hyperlink{H1'}{(\textrm{H1'})}-\hyperlink{H3'}{(\textrm{H3'})}, \hyperlink{H4}{(\textrm{H4})}-\hyperlink{H8}{(\textrm{H8})} and \hyperlink{N1}{(\textrm{N1})}-\hyperlink{N2}{(\textrm{N2})}, for every $\epsilon>0$ there exists a $N^{\epsilon}\in\N$ such that $\bu^N$ is an $\epsilon$-Nash equilibrium for the $N$-player game whenever $N\geq N^{\epsilon}$, i.e. for every $i\in \{1,\ldots,N\}$ and for any deviation $v\in\UNone$
\ba 
J\Ni\round{\bu^{N}}\leq J\Ni\round{\brackets{u^{N,-i},v}}+\epsilon 
\nonumber
\label{eq:approxNash}
\ea
\noindent for all $N\geq N^{\epsilon}$.
\label{cor:approxNash}
\ecor
\noindent
Before proceeding, we define the empirical measure $\zeta^N$ of the $N$-player system (Eq.\eqref{eq:NgameRel}) as
\ba 
\zeta^N\round{\cdot}\doteq\frac{1}{N}\sum_{i=1}^N\delta_{X^{N,i}}\round{\cdot}
\label{eq:empiricalMeasure}
\ea
which is a $\PP(\X)$-valued random variable. Moreover, we fix a relaxed feedback MFG solution $(\lambda,\mu)$ and define (cfr. Theorem \ref{teo:approxNashRel} and Corollary \ref{cor:approxNash}) $\blambda^N\in\UTNN$ as $\blambda^N\doteq (\lambda\Ni)_{i=1,\ldots,N}$ where $\lambda\Ni(t,\varphi^N)\doteq \lambda(t,\varphi\Ni)$ for all $i=1,\ldots,N$, $t\in[0,T]$ and $\varphi^N\in\X^N$. In the next two subsections we consider the following $N$-particle system:
\ba 
X\None_t &=& X\None_0 + \int_{[0,t]\times\Gamma}b\round{s, X\None,\zeta^N,u}\beta\round{s,\textbf{X}^N}(du)ds+\sigma W\None_t, \label{eq:deviation}\\
X\Ni_t &=& X\Ni_0 + \int_{[0,t]\times\Gamma}b\round{s, X\Ni,\zeta^{N},u}\lambda\round{s,X^{N,i}}(du)ds+\sigma W\Ni_t
\label{eq:Npart}
\ea
\noindent for $i = 2,\ldots,N$, $t\in[0,T]$ and where $\beta\in\UTNone$ is a generic single-player control. Precisely, in Subsection \ref{subsec:chaos} we set $\beta(t,\varphi^N)\doteq\lambda(t,\varphi^{N,1})$ for $t\in[0,T]$ and $\varphi^N\in\X^N$ (we say that $\beta=\lambda$ for short); whereas, in Subsection \ref{subsec:approxNash} we let $\beta$ be generic (unless differently specified), which means that we allow the first player to deviate from the MFG solution $\lambda$.

\subsection{Propagation of chaos}\label{subsec:chaos}
\noindent In this subsection we consider the system of $N$ interacting symmetric diffusions given by Eq.s \eqref{eq:deviation} and \eqref{eq:Npart} with $\beta = \lambda$. We associate to this system a suitable McKean-Vlasov equation (Eq.\eqref{eq:mcKean} below) and show a propagation of chaos result, that we will need in the proofs of Theorem \ref{teo:approxNashRel} and Corollary \ref{cor:approxNash}.
\bdefi[\textit{McKean-Vlasov solution}] A law $\theta^*\in\PP(\X)$ is a McKean-Vlasov solution of equation
\ba 
X_t = X_0 + \int_{[0,t]\times\Gamma}b\round{s, X,\theta^*,u}\lambda\round{s,X}(du)ds + \sigma W_t,\quad t\in[0,T], \quad X_0 \simdistr \nu
\label{eq:mcKean}
\ea
\noindent if there exists a weak solution $(\Omega,\F,(\F_t)_{t\in[0,t]},\PB,X,W)$ with $\PB\circ X^{-1}=\theta^*$ and $\PB\circ X^{-1}_0=\nu$.
\label{def:mcKean}
\edefi
\noindent The following lemma ensures the well-posedness of Eq.\eqref{eq:mcKean}.
\blem[\textit{Existence and uniqueness of McKean-Vlasov solutions}]\label{lem:mcKeanEU} Grant Assumptions  \hyperlink{H1'}{(\textrm{H1'})}-\hyperlink{H3'}{(\textrm{H3'})}, \hyperlink{H4}{(\textrm{H4})}-\hyperlink{H8}{(\textrm{H8})} and \hyperlink{N1}{(\textrm{N1})}-\hyperlink{N2}{(\textrm{N2})}. Then, there exists a unique McKean-Vlasov solution for Eq.\eqref{eq:mcKean}.
\elem
\begin{proof}
\noindent We follow \cite{lacker2018strong}, proof of Theorem 2.4. 
Precisely, we apply Banach fixed point theorem on the complete metric space $(\mathcal{P}(\mathcal{X}), d_T)$ together with Picard iterations. To this end, we start by defining, for any $\alpha > 0$, the following distance:
\ba 
d^{\alpha}(\theta, \theta')^{2} \doteq \int_{0}^{T} e^{-\alpha t}d_t(\theta,\theta')^{2}\,dt,\quad \theta,\theta' \in \mathcal{P}(\mathcal{X}).
\nonumber
\ea
\noindent We note that $d^{\alpha}(\cdot, \cdot)$ is a complete metric on $\mathcal{P}(\mathcal{X})$. 
We now define $\Psi:\mathcal{P}\left(\mathcal{X}\right)\rightarrow\PP(\X)\subset \mathcal{P}\left(\mathcal{X}\right)$ as the map $\theta\mapsto \Psi(\theta) \doteq \PB^{\theta}\circ (X^{\theta})^{-1}$ where $(\Omega^{\theta},\F^{\theta},\PB^{\theta},X^{\theta},W^{\theta})$ is a weak solution of Eq.\eqref{eq:mcKean} with $\theta$ in the drift, which is well defined (see Remark \ref{rem:existenceUniqWeak}).\\
\indent We show that $\Psi$ is a contraction on $ \mathcal{P}\left(\mathcal{X}\right)$ with respect to the distance $d^{\alpha}$ for a sufficiently large $\alpha>0$. Let $\mathcal{H}(\theta |\theta')$ denote the relative entropy of $\theta$ with respect to $\theta'$ for $\theta,\theta'\in\PP(\X)$, and let $\mathcal{H}_t(\theta |\theta')=\mathcal{H}(\theta_t |\theta'_t)$, $\theta_t\doteq \PB^{\theta}\circ (X^{\theta}_{\cdot\wedge t})^{-1}$. By Pinsker's inequality, there exists a constant $C_H >0$ such that
\ba 
d_t(\Psi(\theta),\Psi(\theta'))^2&\leq & C_H\mathcal{H}_t(\Psi(\theta),\Psi(\theta'))\nonumber\\
&\le &\frac{1}{2}C_H |\sigma^{-1}|^2\tilde{L}^2\int_0^td_s(\theta,\theta')^2ds 
\nonumber
\ea
where we set $\tilde{L}\doteq L^{TV}_b$. Therefore, we have
\ba
d^{\alpha}(\Psi(\theta),\Psi(\theta'))^2&=&  \int_0^Te^{-\alpha t}d_t(\Psi(\theta),\Psi(\theta'))^2dt \nonumber\\
&\leq &\frac{1}{2}C_H |\sigma^{-1}|^2\tilde{L}^2\int_0^Te^{-\alpha t}\int_0^td_s(\theta,\theta')^2ds\,dt \nonumber\\
&=&\frac{1}{2}C_H|\sigma^{-1}|^2\tilde{L}^2\int_0^Td_t(\theta,\theta')^2 \int_t^Te^{-\alpha s}ds\,dt \nonumber\\ 
&\leq & \frac{1}{2}\frac{C_H}{\alpha}|\sigma^{-1}|^2\tilde{L}^2\int_0^Te^{-\alpha t}d_t(\theta,\theta')^2dt = \frac{1}{2}\frac{C_H}{\alpha}|\sigma^{-1}|^2\tilde{L}^2d^{\alpha}(\theta,\theta')^2 
\nonumber
\ea
which shows that $\Psi$ is a contraction whenever $\frac{1}{2}\frac{C_H}{\alpha}|\sigma^{-1}|^2\tilde{L}^2<1$. Thanks to the arbitrariness of $\alpha>0$, we conclude that $\Psi$ has a unique fixed-point in $\PP(\X)$.
\end{proof}

\noindent We consider the sequence of empirical measures $(\zeta^N)_{N\in\N}$ in Eq.\eqref{eq:empiricalMeasure} associated to the $N$-particle systems in Eq.s \eqref{eq:deviation} and \eqref{eq:Npart} (with $\beta=\lambda$). We follow \cite{lacker2018strong} and we prove the convergence, both in law and in probability in the $\tau$-topology, of $(\zeta^N)_{N\in\N}$ to the McKean-Vlasov solution $\theta^*\in\PP(\X)$ of Eq.\eqref{eq:mcKean}. We remind that the $\tau$-topology on $\mathcal P(\X)$, denoted with $\tau(\PP(\X))$, is the topology generated by the sets
\ba
B_{f,x,\delta} \doteq\cbr{\pi\in\PP(\X):\modl{\int_{\X}f(y)\pi(dy)-x}<\delta}
\nonumber
\ea
\noindent where $f:\X\rightarrow\R$ is any measurable bounded function, $x\in\R$ and $\delta$ is any strictly positive constant. In particular, the $\tau$-topology is the coarsest topology that makes the maps $\pi\mapsto \int_{\X}f(y)\pi(dy)$ continuous for all measurable bounded functions $f:\X\rightarrow\R$ (see, for instance, Chapter 6.2 in \citet{dembo}). \\
\indent Moreover, we denote by $w(\PP(\X))$ the weak topology on $\PP(\X)$ and with $\B(\PP(\X))$ the Borel $\sigma$-algebra on $\X$ generated by the open sets of the weak topology. The following lemma adapts Theorem 2.6.1-2 in \cite{lacker2018strong} to our framework, in particular to the case of diffusions with possibly unbounded drift.
\blem[\textit{Propagation of chaos}]\label{lem:propagationChaos} Grant Assumptions \hyperlink{H1'}{(\textrm{H1'})}-\hyperlink{H3'}{(\textrm{H3'})}, \hyperlink{H4}{(\textrm{H4})}-\hyperlink{H8}{(\textrm{H8})} and \hyperlink{N1}{(\textrm{N1})}-\hyperlink{N2}{(\textrm{N2})}. Let $\theta^*\in\PP(\X)$ be the unique McKean-Vlasov solution of Eq.\eqref{eq:mcKean}. Then the sequence $(\zeta^N)_{N\in\N}$ converges in law to $\theta^*$, i.e. $\zeta^N\overset{\mathcal{L}}{\longrightarrow} \theta^*$, as $N\rightarrow \infty$. Moreover
\ba 
\lim_{N\rightarrow\infty}\PB^N\round{\zeta^N \not\in B}=0
\nonumber
\ea
\noindent for all open neighbourhoods $B$ of $\theta^*$ in the $\tau$-topology that are in $\mathcal{B}(\PP(\mathcal{X}))$.
\elem
\begin{proof}
\noindent Let $(\Omega,\F,\PB)$ be a probability space that supports an i.i.d. sequence of $\X$-valued random variables with law $\theta^*$. For each $N\in\N$, set $(\F^N_t)_{t\in[0,T]}$ to be the filtration generated by $X^1,\ldots,X^N$. Define
\ba 
W^i_t\doteq \sigma^{-1}\round{X^i_t-\xi-\int_{[0,t]\times\Gamma}b(s,X^i,\theta^*,u)\lambda(s,X^i)(du)ds},\,\, t\in[0,T],\,\, i\in\{1,\ldots,N\}.
\nonumber
\ea
\noindent In particular, $W^1,\ldots,W^N$ are independent Wiener processes on $(\Omega,\F,(\F^N_t)_{t\in[0,T]},\PB)$. Fix $N\in\N$, and consider the tuple $(\Omega,\F,(\F^N_t)_{t\in[0,T]},\PB,(X^{N,1},\ldots,X^{N,N}),(W^1,\ldots,W^N))$, with $X\Ni\doteq X^i$, for all $i \in \{1,\ldots,N\}$. This is a weak solution of
\ba 
X\Ni=\xi+\int_{[0,t]\times\Gamma}b(s,X\Ni,\theta^*,u)\lambda(s,X\Ni)(du)ds+\sigma W^i_t,\quad t\in[0,T],\quad i\in\{1,\ldots,N\}.
\nonumber
\ea
\noindent Now, define the probability $\PB^N$ via its density with respect to $\PB$, $\frac{d\PB^N}{d\PB}\doteq Z^N_T$, where, for all $t\in [0,T]$
\ba 
Z^N_t\doteq\mathcal{E}_t\round{\int_0^{\cdot}\sum_{i=1}^N\int_{\Gamma}\sigma^{-1}\round{b(s,X\Ni,\zeta^N,u)-b(s,X\Ni,\theta^*,u)}\lambda(s,X\Ni)(du)dW^i_s}.
\nonumber
\ea
A standard application of Girsanov's theorem gives
\ba 
X\Ni_t=\xi+\int_{[0,t]\times\Gamma}b(s,X\Ni,\zeta^N,u)\lambda(s,X\Ni)(du)ds+\sigma W\Ni_t,\quad t\in[0,T],\,\,\, i\in\{1,\ldots,N\}
\nonumber
\ea
\noindent for some $\PB^N$-Wiener process $W^N$. Notice that $(\Omega,\mathcal{F},(\mathcal{F}^N_t)_{t\in[0,T]},\mathbb{P}^N,X^N,W^N)$ is a weak solution of the $N$-particle system in Eq.s \eqref{eq:deviation} and \eqref{eq:Npart}, with $\beta(t,\varphi^N)\doteq\lambda(t,\varphi^{N,1})$ for $t\in[0,T]$ and $\varphi^N\in\X^N$. 

At this point, the rest of the proof can be performed as in \cite{lacker2018strong}, Theorem 2.6.1-2, along the following steps:
\bi 
\item[(i)] Show that $F_{t_1,t_2}:\PP(\X)\rightarrow\R$ defined as
\ba \label{eq:F12}
F_{t_1,t_2}(\theta)\doteq \int_{\X}\int_{t_1}^{t_2}\modl{\int_{\Gamma}\sigma^{-1}\round{b(s,\varphi,\theta,u)-b(s,\varphi,\theta^*,u)}\lambda(s,\varphi)(du)}^2ds\theta(d\varphi)
\ea
is $\tau$-continuous for all $t_1,t_2\in[0,T]$, $t_1<t_2$ and $\mathcal{B}(\PP(\X))$-measurable, which is done aside at the end of this proof. Moreover $F_{t_1,t_2}(\theta)\leq \tilde{L}(t_2-t_1)\mathcal{H}(\theta\vert\theta^*)$ for all $t_1,t_2\in[0,T]$, $t_1<t_2$ and for all $\theta\in\PP(\X)$, which is a straightforward consequence of the Lipschitz continuity in the total variation distance.
\item[(ii)] Since $X^{N,1},X^{N,2},\ldots X^{N,N}$ are i.i.d. under $\mathbb{P}$, Sanov's Theorem (e.g. Theorem 6.2.10 in \citet{dembo}) can be applied to $\PB\circ(\zeta^N)^{-1}$.
\item[(iii)] Derive a large deviation principle for $\PB^N\circ(\zeta^N)^{-1}$, precisely
\ba 
\underset{N\rightarrow\infty}{\lim\sup}\frac{1}{N}\log\PB^N\round{\zeta^N\not\in B}\leq -\e^{-\tilde{L}T}\inf_{\theta\not\in B}\mathcal{H}\round{\theta\vert\theta^*}
\nonumber
\ea
\noindent for all open neighbourhoods $B$ of $\theta$ in the $\tau$-topology that are in $\mathcal{B}(\PP(\X))$, for some constant $\tilde{L}>0$.\\
\noindent To this aim, we stress that we can proceed just as in \cite{lacker2018strong}\footnote{Precisely we can show by induction that Eq.(4.1) in \cite{lacker2018strong} holds also in this case, then conclude observing that $\PB^N$ and $\PB$ agree on $\F_0$.}. Indeed, regardless of the sub-linear growth of the drift, we can adapt Lacker's estimates thanks to
\ba 
\modl{b\round{t,\varphi,\theta,u}-b\round{t,\varphi,\theta',u}}\leq 2\tilde{L}.
\nonumber\ea
\noindent Moreover we can apply Varadhan's integral lemma \citep[Theorem 4.3.1]{dembo} thanks to the continuity of $F_{t_1,t_2}$.
\item[(iv)] Conclude by showing that $\inf_{\theta\not\in B}\mathcal{H}(\theta\vert\theta^*)> 0$ so that 
\ba 
\lim_{N\rightarrow\infty}\PB^N\round{\zeta^N\not\in B}=0
\nonumber
\ea
\noindent which can be performed as in \cite{lacker2018strong}.
\ei
\textit{Proof of the continuity of $F_{t_1,t_2}$ in the $\tau$-topology.} We actually prove the stronger claim that the functional $F_{t_1,t_2}$ in Eq.\eqref{eq:F12} is continuous in the weak topology ($w$-topology for short). First, we can write $F_{t_1,t_2}(\theta) = \int_{\X} f_{t_1,t_2}(\varphi,\theta) \theta(d\varphi)$ for $\theta \in \mathcal P(\X)$, where
\[ f_{t_1,t_2}(\varphi,\theta) \doteq  \int_{t_1} ^{t_2} \modl{\int_{\Gamma}\sigma^{-1}\round{b(s,\varphi,\theta,u)-b(s,\varphi,\theta^*,u)}\lambda(s,\varphi)(du)}^2 ds\]
which is a real-valued bounded measurable function defined on $\X \times \mathcal P(\X)$. Let $(\theta^n)_{n\in\N},\theta\in\mathcal{P}(\X)$ be such that  $\theta^n\weakconv \theta$. We want to show that $F_{t_1,t_2}(\theta^n)\rightarrow F_{t_1,t_2}(\theta)$ as $n \to \infty$.\\
\indent Set $f_n(\varphi)\doteq f_{t_1,t_2}(\varphi,\theta^n)$ and $f(\varphi)\doteq f_{t_1,t_2}(\varphi,\theta)$. They are all in $C_b(\X)$ with uniform bound in $n\in\N$. Moreover, $f_n\rightarrow f$ in the sup-norm. Indeed
\begin{eqnarray}
\sup_{\varphi\in\X}| f_n(\varphi)- f(\varphi)|
\leq 4 L^{TV}_b L\int_{t_1}^{t_2}\left| L(s;\theta^n)-L(s;\theta) \right|+\left| m_w(s;\theta^n)-m_w(s;\theta) \right| ds 
\nonumber
\end{eqnarray}
\noindent which vanishes in the limit for $n\rightarrow\infty$ due to Lemma \ref{lem:cont-l-m}. As a consequence, we obtain
\ba 
F_{t_1,t_2}(\theta^n)=\int_{\X}f_n(\varphi)\theta^n(d\varphi)\conv \int_{\X} f(\varphi)\theta(d\varphi)=F_{t_1,t_2}(\theta).
\nonumber
\ea
\end{proof}
\subsection{Proof of the The $N$-player approximation theorem}\label{subsec:approxNash}
\noindent This section is devoted to the construction of approximate Nash equilibria for the $N$-player game from a solution of the limit problem, in the particular case of finite-dimensional interaction as described before. The results of previous Subsection \ref{subsec:chaos} allow us to pass to the many-player limit even if feedback MFG strategies are discontinuous in the state variable. We have observed in the introduction that the construction of approximated Nash equilibria for the $N$-player games in \cite{campi2018n} was crucially based on the continuity of the limit optimal control for almost every paths of the state variable with respect to the Wiener measure. In our setting, such a regularity property is no longer feasible due to the possible unboundedness of the coefficients, which makes it difficult to apply PDE-based estimates as in  \cite{campi2018n} to get the needed continuity. Therefore, in order to overcome this obstacle, we will use the strong form of propagation of chaos in Lemma \ref{lem:propagationChaos}, which allows to pass to the limit even through possibly discontinuous MFG optimal controls.\\
\indent In this part, we consider the dynamics in Eq.\eqref{eq:deviation} and Eq.\eqref{eq:Npart} without necessarily taking $\beta=\lambda$, unless differently specified. We start with some preliminary estimates ensuring that the costs remain bounded in the mean-field limit despite the sub-linear growth.
\blem[\textit{A-priori estimates}]\label{lem:NgameEstimates} Grant Assumptions \hyperlink{H1'}{(\textrm{H1'})}-\hyperlink{H3'}{(\textrm{H3'})}, \hyperlink{H4}{(\textrm{H4})}-\hyperlink{H8}{(\textrm{H8})} and \hyperlink{N1}{(\textrm{N1})}-\hyperlink{N2}{(\textrm{N2})}. Consider the dynamics in Eq.s \eqref{eq:deviation} and \eqref{eq:Npart}. Then for any $\alpha\geq 1$
\ba 
\sup_{N\in\N}\E^{\PB^{N}}\brackets{\infnorm{X\Ni}^{\alpha}}&\leq & K(\alpha)
\nonumber
\ea
\noindent for $i\in\{1,\ldots,N\}$ and where $K(\alpha)<\infty$ is a positive constant independent of $N$.
\elem
\begin{proof}
\noindent This is a consequence of Gr\"onwall's lemma together with uniform boundedness of the drift in the measure and control variables.
\end{proof}
\noindent Now, we prove the tightness of the sequence of laws $(\PB^N\circ (\zeta^N)^{-1})_{N\in\N}$ when $\beta=\lambda$ in Eq.\eqref{eq:deviation}, i.e. when the dynamics are symmetric. Then, thanks to Lemma \ref{lem:propagationChaos}, we characterize the limit points of $(\PB^N\circ (\zeta^N)^{-1})_{N\in\N}$ as McKean-Vlasov solutions of Eq.\eqref{eq:mcKean}; see Lemma \ref{lem:nplayerLimitPoints}.
\blem[\textit{Tightness}]\label{lem:nplayerTight} Grant Assumptions \hyperlink{H1'}{(\textrm{H1'})}-\hyperlink{H3'}{(\textrm{H3'})}, \hyperlink{H4}{(\textrm{H4})}-\hyperlink{H8}{(\textrm{H8})} and \hyperlink{N1}{(\textrm{N1})}-\hyperlink{N2}{(\textrm{N2})}. Let $\zeta^N$ be the empirical measure of the system given by Eq.s \eqref{eq:deviation} and \eqref{eq:Npart} with $\beta=\lambda$. Then the sequence $(\PB^N\circ (\zeta^N)^{-1})_{N\in\N}$ is tight in $\PP(\PP(\mathcal{X}))$.
\elem
\begin{proof}
\noindent The tightness of such a sequence follows from \cite{sznitman}, Proposition 2.2, combined with Kolmogorov-Chentsov criterion (see, for instance, Corollary 14.9 in \citet{kallenberg}).
\end{proof}

\blem[\textit{Characterization of limit points}]\label{lem:nplayerLimitPoints} Grant Assumptions (H1')-(H3'), (H4)-(H8) and (N1)-(N2). Let $\zeta^N$ be the empirical measure of the system given by Eq.s \eqref{eq:deviation} and \eqref{eq:Npart} with $\beta=\lambda$. Let $(\PB^{N_k}\circ (\zeta^{N_k})^{-1})_{k\in\N}$ be a convergent subsequence of $(\PB^N\circ (\zeta^N)^{-1})_{N\in\N}$. Let $\zeta$ be a random variable defined on some probability space $(\Omega,\F,\PB)$ with values in $\PP(\X)$ such that $\zeta^{N_k}\overset{\mathcal{L}}{\longrightarrow}\zeta$. Then
\bi 
\item[(i)] $\zeta$ coincides $\PB$-a.s. with the unique McKean-Vlasov solution $\theta^*$ of Eq.\eqref{eq:mcKean}.
\item[(ii)] The sequence $(\zeta^N)_{N\in\mathbb{N}}$  converges in probability (hence also in law) to $\theta^*$ when $\mathcal P(\X)$ is equipped with the $\tau$-topology.
\ei
\elem
\begin{proof}
\noindent By Lemma \ref{lem:nplayerTight} there exists a subsequence $(\PB^{N_k}\circ (\zeta^{N_k})^{-1})_{k\in\N}\subset \PP(\PP(\X))$ converging to $\PB\circ\zeta^{-1}\in\PP(\PP(\X))$. 
Lemma \ref{lem:propagationChaos} guarantees the convergence in law of the whole sequence $(\zeta^N)_{N\in\mathbb{N}}$ to the deterministic limit $\theta^*$, which is the unique McKean-Vlasov solution of Eq.\eqref{eq:mcKean}. By uniqueness in law of the weak limit we have $\PB\circ\zeta^{-1}=\delta_{\theta^*}$, yielding $\zeta=\theta^*$ $\PB$-a.s.. Lemma \ref{lem:propagationChaos} also gives convergence in probability in the $\tau$-topology of $(\zeta^N)_{N\in\mathbb{N}}$ to $\theta^*$.
\end{proof}
\bcor[\textit{Characterization of the convergence}]\label{cor:nplayerLimitPoints} Under the assumptions of Lemma \ref{lem:nplayerLimitPoints}, the following properties hold:
\bi 
\item[(i)] For all Borel-measurable bounded function $f:\mathcal{X}\rightarrow\mathbb{R}$ such that $\theta\mapsto\int_{\mathcal{X}}f(\varphi)\theta(d\varphi)$ is $\tau(\PP(\mathcal{X}))$-continuous
\ba 
\E^{\PB^N}\brackets{\int_{\mathcal{X}}f(\varphi)\zeta^N(d\varphi)}\underset{N\rightarrow\infty}{\longrightarrow}\mathbb{E}^{\PB}\brackets{\int_{\mathcal{X}}f(\varphi)\zeta(d\varphi)}\equiv \mathbb{E}^{\PB}\brackets{\int_{\mathcal{X}}f(\varphi)\theta^*(d\varphi)}.
\nonumber
\ea
\item[(ii)] $\PB^N\circ(X^{N,1},\zeta^N)^{-1}\weakconv \theta^*\otimes\delta_{\theta^*}$. Moreover, $\PB^N\circ(X^{N,1})^{-1}\weakconv \theta^*$ and $\PB^N\circ(\zeta^N)^{-1}\weakconv \delta_{\theta^*}$.
\item[(iii)] For all $f\in C(\X)$ with sub-linear growth, i.e. $|f(\varphi)|\leq C_f(1+\infnorm{\varphi})$ for some $C_f>0$ and all $\varphi \in \X$, we have
\ba 
\E^{\PB^N}\brackets{\int_{\mathcal{X}}f(\varphi)\zeta^N(d\varphi)}\underset{N\rightarrow\infty}{\longrightarrow}\mathbb{E}^{\PB}\brackets{\int_{\mathcal{X}}f(\varphi)\zeta(d\varphi)}\equiv \mathbb{E}^{\PB}\brackets{\int_{\mathcal{X}}f(\varphi)\theta^*(d\varphi)}.
\nonumber
\ea
\ei
\ecor
\begin{proof}
\noindent (i)\, This is a consequence of Lemma \ref{lem:propagationChaos}, Lemma \ref{lem:nplayerLimitPoints} and of the almost sure equality $\zeta=\theta^*$.\\
\indent (ii)\, We already know that $\PB^N\circ(\zeta^N)^{-1}\weakconv \delta_{\theta^*}$ from Lemma \ref{lem:nplayerLimitPoints}. Therefore, the convergence of $\PB^N\circ(X^{N,1})^{-1}$ to $\theta^*$  follows from \cite{sznitman}, Proposition 2.2, and the symmetry of the system.\\
\indent (iii)\, Let $f\in C(\X)$ with sub-linear growth. It is enough to show that
\ba 
\E^{\PB^N}\brackets{\int_{\mathcal{X}}\infnorm{\varphi}\zeta^N(d\varphi)}\underset{N\rightarrow\infty}{\longrightarrow} \int_{\mathcal{X}}\infnorm{\varphi}\theta^*(d\varphi).
\nonumber
\ea
\noindent To this aim, for fixed $R>0$, we consider the decomposition
\ba 
\E^{\PB^N}\brackets{\int_{\mathcal{X}}\infnorm{\varphi}(\zeta^N-\theta^*)(d\varphi)}&\leq &\E^{\PB^N}\brackets{\int_{\mathcal{X}}(\infnorm{\varphi}\wedge R) (\zeta^N-\theta^*)(d\varphi)}\nonumber\\
&&+\E^{\PB^N}\brackets{\int_{\mathcal{X}}\infnorm{\varphi}\mathbf{1}_{\{\infnorm{\varphi}\geq R\}}(\zeta^N+\theta^*)(d\varphi)}.
\nonumber
\ea
\noindent By property (i), for any fixed $R>0$, we have
\ba 
\lim_{N\rightarrow\infty}\E^{\PB^N}\brackets{\int_{\mathcal{X}}(\infnorm{\varphi}\wedge R) (\zeta^N-\theta^*)(d\varphi)}=0
\nonumber
\ea
\noindent so that 
\ba 
\underset{N\rightarrow\infty}{\lim\sup}\E^{\PB^N}\brackets{\int_{\mathcal{X}}\infnorm{\varphi}(\zeta^N-\theta^*)(d\varphi)}
\leq \underset{N\rightarrow\infty}{\lim\sup}\E^{\PB^N}\brackets{\int_{\mathcal{X}}\infnorm{\varphi}\mathbf{1}_{\{\infnorm{\varphi}\geq R\}}(\zeta^N+\theta^*)(d\varphi)}.
\nonumber
\ea
\noindent Now, we let $R\rightarrow\infty$ and we show that the RHS vanishes in the limit. 
To do so, recall that, due to Lemma \ref{lem:NgameEstimates}, there exist constants $K(\alpha),K>0$ such that
\ba  \sup_{N\in\N}E^{\PB^N}\brackets{\infnorm{X^{N,i}}^{\alpha}}\leq K(\alpha)\quad \text{and}\quad\sup_{N\in\N}E^{\PB^N}\brackets{\infnorm{X^{N,i}}}\leq K
\nonumber
\ea
independently of $i\in\{1,\ldots,N\}$. Then, set $\alpha,\beta> 1$ such that $\frac{1}{\alpha}+\frac{1}{\beta}=1$ and let $\epsilon>0$. By definition of $\zeta^N$ and by Young's and Markov's inequalities, we have
\ba 
\underset{N\rightarrow\infty}{\lim\sup}\E^{\PB^N}\brackets{\int_{\mathcal{X}}\infnorm{\varphi}\mathbf{1}_{\{\infnorm{\varphi}\geq R\}}\zeta^N(d\varphi)}&=&\underset{N\rightarrow\infty}{\lim\sup}\frac{1}{N}\sum_{i=1}^N\E^{\PB^N}\brackets{\infnorm{X^{N,i}}\mathbf{1}_{\{\infnorm{X^{N,i}}\geq R\}}}\nonumber\\
&\leq &\round{\epsilon^{\alpha}\frac{K(\alpha)}{\alpha}+\frac{K}{\epsilon^{\beta}\beta R}}\\
\nonumber
\ea
which converges to zero by letting $R \to \infty$ and then $\epsilon \to 0$. A similar reasoning applies to the same expectation with $\theta^*$ instead of $\zeta^N$.
\end{proof}
\brem 
\noindent Let $\mathbb D\doteq\{\varphi\in\mathcal{X}:\tau(\varphi)\text{ is discontinuous at }\varphi\}$. Since $\zeta\overset{a.s.}{=}\theta^*\in\mathcal{Q}$, Lemma \ref{lem:regularity} implies $\theta^*(\mathbb D)=0$ and the statement of Corollary \ref{cor:nplayerLimitPoints} holds for $f=\mathbf{1}_{\mathbb D}$ as well.
\erem
\noindent Finally, we conclude this section with the proof of Theorem \ref{teo:approxNashRel}, which leads immediately to Corollary \ref{cor:approxNash}. \\
\begin{proof}[Proof of Theorem \ref{teo:approxNashRel}] The proof is structured in three steps.
\bi 
\item[(j)]\quad $\lim_{N\rightarrow\infty}J^{N,1}(\blambda^N)=J^{\mu}(\lambda)$.
\item[(jj)]\quad Let $\beta^{N,1}\in\mathcal{U}^N_1$ be such that
\ba
J^{N,1}([\blambda^{N,-1},\beta^{N,1}])\leq\inf_{\beta\in\mathcal{U}^N_1}J^{N,1}([\blambda^{N,-1},\beta])+\frac{\epsilon}{2}.
\nonumber
\ea
\noindent Then
\ba 
\underset{N\rightarrow\infty}{\lim\inf}J^{N,1}\round{\brackets{\blambda^{N,-1},\beta^{N,1}}}\geq  J^{\mu}(\lambda).
\nonumber
\ea
\item[(jjj)]\quad $J^{N,1}(\blambda^N)\leq\inf_{\beta\in\mathcal{U}^N_1}J^{N,1}([\blambda^{N,-1},\beta])+\epsilon$.
\ei
\noindent We consider the dynamics in Eq.\eqref{eq:NgameRel}. In (j) we set $\lambda^{N,1}(t,\varphi^N)=\lambda(t,\varphi^{N,i})$ for all $(t,\varphi^N)\in[0,T]\times\X^N$ and prove convergence of the first-player cost functional to the cost functional of the MFG. In (jj) instead we allow the first player to deviate and choose $\lambda^{N,1}(t,\varphi^N)=\beta^{N,1}(t,\varphi^{N})$ for all $(t,\varphi^N)\in[0,T]\times\X^N$ where $\beta^{N,1}\in\UTNone$ is a generic single-player relaxed control. We conclude the proof in (jjj) by combining the results in (j) and (jj).\\
\indent\textit{Proof of (j).}\quad To prove that $J^{N,1}(\blambda^N)\rightarrow J^{\mu}(\lambda)$, as $N\rightarrow \infty$, we split each cost functional in the sum of two terms:
\ba 
J^{N,1}(\blambda^N)&=&\E^{\PB^N}\brackets{\int_{[0,T]\times\Gamma}\int_{\X}\mathbf{1}_{[0,\tau(\varphi))}(t)f_0(t,\varphi,u)\lambda(t,\varphi)(du)\zeta^N(d\varphi)dt}\nonumber\\
&&+\,\E^{\PB^N}\brackets{\int_{0}^{T}\mathbf{1}_{[0,\tau^{N,1})}(t)f_1(t,X^{N,1},\zeta^N)dt+F(\tau^{N,1},X^{N,1}_{\tau^{N,1}})}
\nonumber
\ea
\noindent and
\ba
J^{\mu}(\lambda)&=&\E^{\PB}\brackets{\int_{[0,T]\times\Gamma}\int_{\X}\mathbf{1}_{[0,\tau(\varphi))}(t)f_0(t,\varphi,u)\lambda(t,\varphi)(du)\zeta(d\varphi)dt}\nonumber\\
&& +\,\E^{\PB}\brackets{\int_{0}^{T}\mathbf{1}_{[0,\tau)}(t)f_1(t,X,\zeta)dt+F(\tau,X_{\tau})}.
\nonumber
\ea
\noindent Since $f_0$ is bounded, the convergence of the first summand in the decomposition of $J^{N,1}(\blambda^N)$ to the corresponding term in $J^{\mu}(\lambda)$ is a consequence of Corollary \ref{cor:nplayerLimitPoints}(i) and of Lemma \ref{lem:nplayerLimitPoints}. On the other hand, since both $f_1$ and $F$ have sub-linear growth, the convergence of the second summand in $J^{N,1}(\blambda^N)$ follows from Corollary \ref{cor:nplayerLimitPoints}(iii), Lemma \ref{lem:nplayerLimitPoints} and the fact that $\theta^*\in\Q$ together with Lemma \ref{lem:regularityConvergence}.\\
\indent\textit{Proof of (jj).}\quad We follow the proof of Theorem 3.10 in \cite{lacker2020convergence} with suitable modifications due to the possibly unbounded drift and the dependence on the first exit time from the set $\OO$.\\
\noindent Let $(\Omega^N,\F^N,(\F^N_t)_{t\in[0,T]},\QB^N,Y^N,W^N)_{N\in\mathbb{N}}$ be a weak solutions of the $N$-player system. Let $(\zeta^N)_{N\in\N}$ be the associated empirical measures. Under $\QB^N$ the first player's dynamics is
\ba 
Y^{N,1}_t &=& Y^{N,1}_0 + \int_{[0,t]\times\Gamma}b(s,Y^{N,1},\zeta_Y^N,u)\beta^{N,1}(s,\textbf{Y}^{N})(du)ds + \sigma W^{N,1}_t,\quad t\in[0,T].
\nonumber
\ea
\noindent Now, let $\PB^N$ be the probability measure under which the first player's dynamics becomes
\ba 
Y^{N,1}_t &=& Y^{N,1}_0 + \int_{[0,t]\times\Gamma}b(s, Y^{N,1},\zeta_Y^N,u)\lambda(s,Y^{N,1})(du)ds + \sigma \tilde{W}^{N,1}_t,\quad t\in[0,T]
\nonumber
\ea
where $ \tilde{W}^{N,1}$ is a $\PB^N$-Wiener process. In other terms, $\PB^N$ satisfies $\frac{d\QB^N}{d\PB^N}= Z^N_T$ where
\ba 
Z^N_t=\mathcal{E}_{t}\round{\int_0^{\cdot}\int_{\Gamma}b(s,Y^{N,1},\zeta_Y^N,u)(\beta^{N,1}(s,\textbf{Y}^{N})-\lambda(s,Y^{N,1}))(du)d\tilde{W}_s}, \quad t\in [0,T].
\nonumber
\ea
\noindent By inspection of the proofs of Lemma \ref{lem:benevs} and Corollary \ref{cor:expmartMoments}, all bounds are uniform in $N\in\N$, hence Corollary \ref{cor:expmartMoments} gives the uniform integrability of the sequence of  exponential martingales $(Z^N)_{N\in\N}$. More in detail, we apply Corollary \ref{cor:expmartMoments} to the drift
\ba 
b(t,\varphi^{N})\doteq \int_{\Gamma}b(t,\varphi^{N,1},\zeta_{\varphi^N},u)(\beta^{N,1}(t,\varphi^{N})-\lambda(t,\varphi^{N,1}))(du)
\nonumber
\ea
\noindent for $(t,\varphi^{N})\in[0,T]\times\X^N$. Notice that this drift is sublinear in $\varphi^N$. Therefore convergence of the empirical measures to $\theta^*$ in probability in the $\tau$-topology under $\PB^N$ implies convergence of the empirical measures to the same limit in probability in the $\tau$-topology under $\QB^N$. Hence $\zeta_Y^N\overset{\mathcal{L}}{\longrightarrow}\theta^*$ under $\QB^N$ and
\ba 
\lim_{N\rightarrow\infty}\QB^N\round{\zeta_Y^N\not\in B}=0
\nonumber
\ea
\noindent for all neighbourhoods $B$ of $\theta$ in the $\tau$-topology which belong to $\mathcal{B}(\PP(\mathcal{X}))$. The tightness of $(Y^{N,1})_{N\in\N}$ under $\QB^N$ still follows from their tightness under $\PB^N$. Consider $(\beta^{N,1}(t,\textbf{Y}^N))_{t\in[0,T]}$ as a single-player relaxed stochastic open-loop control and denote it simply by $(\beta^{N,1}_t)_{t\in[0,T]}$. 
Interpret $(Y^{N,1},\beta^{N,1},\zeta^N_Y)_{N\in\mathbb{N}}$ as a sequence of random variables with values in $\mathcal{X}\times\mathcal{V}\times\PP(\X)$. 
Compactness of $\mathcal{V}$ and tightness of $(Y^{N,1},\zeta^N_Y)_{N\in\mathbb{N}}$ imply the tightness of $(Y^{N,1},\beta^{N,1},\zeta^N_Y)_{N\in\mathbb{N}}$ under $\QB^N$. 

Let $(Y,\beta,\theta^*)$ be a limit point of the sequence $(Y^{N,1},\beta^{N,1},\zeta^N_Y)_{N\in\mathbb{N}}$, defined on some probability space with probability measure $\QB$. Then by a standard martingale argument it can be shown to satisfy
\ba 
Y_t=\xi+\int_{[0,t]\times\Gamma}b(s,Y,\theta^*,u)\beta_t(du)ds+\sigma W_t,\quad t\in[0,T]
\label{eq:representDeviation}
\ea
\noindent where $W$ is a $\QB$-Wiener process. As in (j) we split $J^{N,1}([\lambda^{N,-1},\beta^{N,1}])$ in two terms as
\ba 
J^{N,1}([\lambda^{N,-1},\beta^{N,1}])
&=&\mathbb{E}^{\QB^N}\brackets{\int_{[0,T]\times\Gamma}\mathbf{1}_{[0,\tau^{N,1})}(t)f_0(t,Y^{N,1},u)\beta^{N,1}_t(du)dt}\nonumber\\
&&+ \, \mathbb{E}^{\QB^N}\brackets{\int_{0}^T\mathbf{1}_{[0,\tau^{N,1})}(t)f_1(t,Y^{N,1},\zeta^N_Y)dt+F(\tau^{N,1},Y^{N,1}_{\tau^{N,1}})}.
\nonumber
\ea
\noindent We move along a weakly converging subsequence of $(Y^{N,1},\beta^{N,1},W^{N,1})_{N\in\mathbb{N}}$ under $\QB^N$ to the limit point $(Y,\beta,W)$ in Eq.\eqref{eq:representDeviation}. Convergence of the first and second summands above now works as in the proof of (j). Considering again the whole sequence, we obtain
\ba 
\underset{N\rightarrow\infty}{\lim\inf}J^{N,1}([\lambda^{N,-1},\beta^{N,1}])
&\geq &\inf_{\beta}\E^{\QB^N}\brackets{\int_{[0,T]\times\Gamma}\mathbf{1}_{[0,\tau)}(t)f(t,Y,\theta^*,u)\beta_t(du)dt+F(\tau,Y_{\tau})}\nonumber\\
&= &V^{\mu}
\nonumber
\ea
\noindent where the infimum on the RHS above is taken over all relaxed stochastic open-loop controls and the last equality follows from embedding the set of strict controls into the set of relaxed controls combined with the chattering lemma \citep{karoui,fleming,bahlali}.\\
\indent\textit{Proof of (jjj).}\quad This is a consequence of steps (j) and (jj). Indeed 
\ba 
J^{N,1}(\blambda^N)-\inf_{\beta\in\mathcal{U}^N_1}J^{N,1}([\blambda^{N,-1},\beta])\leq
J^{N,1}(\blambda^N)-J^{\mu}(\lambda)+J^{\mu}(\lambda)-J^{N,1}([\blambda^{N,-1},\beta^{N,1}])+\frac{\epsilon}{2}.
\nonumber
\ea
\noindent Now by steps (j) and (jj) there exists $N^{\epsilon}\in\mathbb{N}$ such that for all $N\geq N^{\epsilon}$
\ba 
J^{N,1}(\blambda^N)-J^{\mu}(\lambda)\leq\frac{\epsilon}{4}\quad\text{and}\quad J^{\mu}(\lambda)-J^{N,1}([\blambda^{N,-1},\beta^{N,1}])\leq\frac{\epsilon}{4}.
\nonumber
\ea
Therefore, we can conclude that $J^{N,1}(\blambda^N)\leq\inf_{\beta\in\mathcal{U}^N_1}J^{N,1}([\blambda^{N,-1},\beta])+\epsilon$ for all $N\geq N^{\epsilon}$, which establishes the statement of Theorem \ref{teo:approxNashRel}.
\end{proof}

\appendix
\begin{center}
\section{Appendix}\label{app:Appendix}
\end{center}
\noindent This appendix provides some of the technical results used in the paper. More in detail, we state existence and uniqueness of weak solutions of SDEs with sub-linear drift. We characterize the space of laws of processes with sub-linear drift and initial condition $\nu$ ($\mathcal{Q}$ defined below). We prove some regularity results on the exit time $\tau^{X}$ with respect to measures in $\Q$. Finally, we discuss the convergence of measures in the 1-Wasserstein distance along test functions with sub-linear growth and possibly discontinuous over a set of limit measure zero.
\\

\subsection{Existence and uniqueness of solution of SDEs with sub-linear drift}
\noindent In this subsection we prove a slight variation of the well-known Bene\v{s}' condition (\citet{benevs}), leading to an existence and uniqueness result for weak solutions of SDEs with a sub-linear drift. More precisely, we allow the drift to depend on a rescaled Wiener process with a independent random initial condition. We recall that $\mathcal{E}_t(\cdot)$ denotes the Dol\'eans-Dade stochastic exponential. Moreover, given a function $f:E\rightarrow \R$ where $E$ is a Polish space, we denote by  $\mathbb D_f$ the set of its discontinuity points.\\
\indent As a preliminary, we introduce the set $\Q$ of laws of stochastic processes with sub-linear drift in the sense of Bene\v{s} to which these results apply.\\

\noindent \textit{Laws of processes with sub-linear drift.}\quad Let $\beta:[0,T]\times \X\rightarrow \mathbb{R}^d$ be a progressively measurable functional such that
\ba 
\modl{\beta\round{t,\varphi}}\leq C\round{1+\infnorm{\varphi}},\quad (t,\varphi)\in[0,T]\times\X
\nonumber
\ea
\noindent for some constant $C>0$. Let $(\Omega,\mathcal{F},(\mathcal{F}_t)_{t\in[0,T]},\mathbb{P},X)$ be a weak solution of the following SDE
\ba
X_t&=&\xi+\int_0^t\beta(s,X)ds+\sigma W_t,\quad \xi\overset{d}{\sim}\nu,\quad t\in[0,T]
\nonumber 
\ea 
\noindent where $W$ is a Wiener process independent of $\xi$. Existence and uniqueness of a weak solution follows from an application of Girsanov's theorem and Bene\v{s}' condition (see Lemma \ref{lem:benevs} and Lemma \ref{lem:existenceUniquenessWeak}). Moreover such laws turn out to be absolutely continuous with respect to the Wiener measure $\wiener$ (Lemma \ref{lem:regularitySubLinear}). Then, we denote by $\mathcal{Q}$ the set of laws $\theta\in\PP(\X)$ of all continuous processes $X$ solving the SDE above.\medskip

\blem[\textit{Bene\v{s}' condition}]\label{lem:benevs} Let $b:[0,T]\times \X\rightarrow \mathbb{R}^d$ be a progressively measurable functional such that
\ba 
\modl{b\round{t,\varphi}}\leq C\round{1+\infnorm{\varphi}},\quad (t,\varphi)\in[0,T]\times\X.
\nonumber
\ea
\noindent Let $\sigma\in\R^{d\times d}$ be a full rank matrix. Let $(\Omega,\mathcal{F},(\mathcal{F}_t)_{t\in[0,T]},\mathbb{P})$ be a filtered probability space satisfying usual conditions, supporting a random variable $\xi\overset{d}{\sim}\nu$ and a Wiener process $W$ independent of $\xi$. Set
\ba
X_t \doteq \xi+\sigma W_t,\quad t\in[0,T].
\nonumber 
\ea 
Then
\ba 
Z_t\doteq \mathcal{E}_t\round{\int_0^{\cdot}\sigma^{-1} b(s,X)dW_s},\quad t\in[0,T]
\nonumber
\ea
\noindent is a martingale.
\elem
\begin{proof} We follow the proof of Corollary 3.5.16 in \cite{karaztas}. Precisely let 
$t_0=0<t_1<\ldots<t_{n-1}<t_n=T$ be a partition of the interval $\brackets{0,T}$. Then thanks to the sub-linearity of the  drift
\ba 
\int_{t_{n-1}}^{t_n}\modl{b(s,X)}^2ds\leq(t_n-t_{n-1})C^2\round{1+\infnorm{X}}^2.
\nonumber
\ea
\noindent 
Let $Y^n\doteq(Y^n_t)_{t\in[0,T]}$ be defined by
\ba 
Y^n_t\doteq \e^{\frac{1}{4}(t_n-t_{n-1})C^2(1+|X_t|)^2}.
\nonumber
\ea Notice that $Y^n$ is a sub-martingale and that by Doob's maximal inequality \citep[Theorem 1.3.8.iv]{karaztas} we have $\E[\infnorm{Y^n}^2]\leq 4\E[(Y^n _T)^2]$. Moreover
\ba 
\E\brackets{(Y^n _T)^2}& \leq &\E\brackets{\e^{\frac{1}{2}(t_n-t_{n-1})C^2(1+2|\xi|^2+2|\sigma|^2|W_T|^2)}}\nonumber \\
&=&\E\brackets{\e^{(t_n-t_{n-1})C^2|\sigma|^2|W_T|^2}}\E\brackets{\e^{\frac{1}{2}(t_n-t_{n-1})C^2(1+2|\xi|^2)}}
\nonumber
\ea
\noindent where in the equality we have used the independence between $\xi$ and $W$. To conclude, it is sufficient to choose $(t_k-t_{k-1})$, $k = 1,\ldots,n$, sufficiently small, for instance $(t_k-t_{k-1})<\min\{\frac{1}{2C^2|\sigma|^2},\frac{\lambda}{C^2}\}$, and to apply Corollary 3.5.14 in \cite{karaztas}. 
\end{proof}
\bcor[\textit{Moments of the stochastic exponential}]\label{cor:expmartMoments} Under the assumptions of Lemma \ref{lem:benevs}, the process $Z=(Z_t)_{t\in[0,T]}$ has finite moments of any order $p\in[1,\infty)$, i.e. $\mathbb{E}\brackets{Z_T^p}<\infty$ for all $p\in[1,\infty)$.
\ecor
\begin{proof}
\noindent The proof follows directly from Lemma \ref{lem:benevs} combined with Corollary 2 in \cite{grigelionis}.
\end{proof}
\blem[\textit{Existence and uniqueness of weak solutions}]\label{lem:existenceUniquenessWeak} Let $b:[0,T]\times \X\rightarrow \mathbb{R}^d$ be a progressively measurable functional such that
\ba 
\modl{b\round{t,\varphi}}\leq C\round{1+\infnorm{\varphi}},\quad (t,\varphi)\in[0,T]\times\X.
\nonumber
\ea
\noindent Let $\sigma\in\R^{d\times d}$ a full rank matrix. Then there exists a weak solution $(\Omega,\F,(\F_t)_{t\in[0,T]},\PB,X,W)$ of
\ba 
X_t=\xi+\int_0^tb(s,X)ds+\sigma dW_t,\quad \xi\simdistr\nu,\quad t\in[0,T].
\nonumber
\ea
\noindent Additionally, this solution is unique in law. 
\elem
\begin{proof}
\noindent The proof follows directly from Lemma \ref{lem:benevs} and Girsanov's theorem \citep[see][Propositions 5.3.6 and 5.3.10]{karaztas}.
\end{proof}

\subsection{Characterization of the set $\Q$}

\blem[\textit{Laws of processes with sub-linear drift}]\label{lem:regularitySubLinear} Let $\theta\in \Q$. Then $\theta\sim\wiener$, i.e. $\theta$ is equivalent to the Wiener measure $\wiener$.
\elem
\begin{proof}
\noindent The proof follows directly from Lemma \ref{lem:benevs}, Girsanov's theorem and Bayes' rule to ensure that $Z^{-1}$ given by Lemma \ref{lem:benevs} is still a martingale.
\end{proof}

\noindent Before proceeding further, we recall that $\tau^{X}$ is the first exit time from $\mathcal{O}$ in the path space, i.e. 
\ba
\tau^{X}(\varphi) = \inf\braces{t \geq 0\,:\,\varphi(t) \not\in \mathcal{O}},\quad \varphi \in \X,
\nonumber
\ea
where $\OO\subset \Rd$ satisfies Assumption (H4).
\blem[\textit{Regularity results}]\label{lem:regularity} Let $\theta\in\Q$. Let $\OO\subset \Rd$ satisfy Assumption (H4) and let $X$ be the identity process on $\X$. Then
\bi 
\item[(a)] $ \tau^{X}<\infty$, $\theta$-almost surely.
\item[(b)] The mapping $\varphi\mapsto \tau^{X}(\varphi)$, from $\X$ to $[0,\infty]$, is $\theta$-a.s. continuous.
\item[(c)] $\theta(\tau^{X}=t)=0$ for all $t\in[0,T]$.
\item[(d)] The mapping $\varphi\mapsto \mathbf{1}_{[0,\tau^{X}(\varphi))}(t)$, from $\X$ to $\R$, is $\theta$-a.s. continuous for all $t\in[0,T]$.
\item[(e)] Properties (a)-(d) hold for $\OO=(0,\infty)^{\times d}$ as well.
\ei
\elem
\begin{proof}
\noindent The proof is similar to the one of Lemma D.3 in \cite{campi2018n}. 
Notice that by Lemma \ref{lem:regularitySubLinear} each $\theta \in\Q$ is equivalent to $\wiener$. So, it is sufficient to check properties (a)-(d) for $\wiener$.\\
\indent (a)\, This is a consequence of the law of iterated logarithms (as time tends to infinity) and the fact that $\OO$ is strictly included in $\mathbb{R}^d$.\\
\indent (b)\, This, again,  is a consequence of the law of iterated logarithms (as time tends to zero), the smoothness of $\OO$'s boundary, the non-degeneracy of $\sigma$ and the fact that $\mathcal{O}$ is strictly included in $\mathbb{R}^d$ (\citet{kushner}, pp. 260-261).\\
\indent (c)\, This is a consequence of the following relations
$$
\wiener(\tau^X=t)\leq \wiener(X_t\in\partial\mathcal{O})=0\quad \textrm{for all } t\in[0,T]
$$
where in the last equality we use the fact that the Lebesgue measure of the boundary of a convex subset of $\mathbb{R}^d$ is identically zero (\citet{lang}), and that $\wiener\circ X_t^{-1}$ is absolutely continuous with respect to the Lebesgue measure for all $t\in[0,T]$.\\
\indent (d)\, This is a consequence of properties (b) and (c) above.\\
\indent (e)\, When $\mathcal{O}=(0,\infty)^{\times d}$ it turns out that
\ba 
\tau^X (\varphi)=\min_{i=1,\ldots,d}\tau^{i}(\varphi),\quad\varphi\in\mathcal{X}
\nonumber
\ea
\noindent where $\tau^i(\varphi)\doteq\inf\{t\in[0,T]:\varphi_i(t)\leq 0\}$, for $i \in \{1,\ldots,d\}$ and $\varphi\in\mathcal{X}$. Then the conclusion follows from the continuity result in dimension $d=1$ (\citet{kushner}, pp. 260-261) applied to each $\tau^{i}$.
\end{proof}
\subsection{Additional convergence results}
\blem[\textit{Convergence in the 1-Wasserstein distance}]\label{lem:regularityConvergence} Let $E$ be a Polish space with a complete metric $d_E$. Let $\theta,(\theta^n)_{n\in\mathbb{N}}\subset\PP_1(E)$ such that $W_1(\theta^n,\theta)\rightarrow 0$ as $n\rightarrow \infty$. Let $f:E\rightarrow \R$ be a measurable function such that $|f(x)|\leq C(1+d_E(x,x_0))$ for all $x\in E$, for some $x_0\in E$ and for some constant $C>0$. Let $\mathbb D_f$ be the set of its discontinuity points and assume $\theta(\mathbb D_f)=0$. Then
\ba 
\int_E f(x)\theta^n(dx)\conv \int_E f(x)\theta(dx).
\nonumber
\ea
\elem
\begin{proof}
\noindent The proof works as in \cite{villani}, proof of Theorem 7.12.iv, the only difference being that here $f$ can have discontinuities with $\theta( \mathbb D_f)=0$. In particular, we perform the same decomposition as in \cite{villani}, i.e. $f(x)=f_R^{1}(x)+f_R^{2}(x)$ with $f_R^{1}(x)\doteq f(x)\wedge (C(1+R))$ and $f_R^{2}(x)\doteq f(x)-f_R^{1}(x)$ for all $x\in E$ and for some $R>0$. We have that $|f^{1}_R|$ is bounded by $C(1+R)$ and  $\theta(\mathbb D_{f^{1}_R})=0$ since $\mathbb D_{f^{1}_R}\subset \mathbb D_f$. Then all limits can be performed just as in \cite{villani}, proof of Theorem 7.12.iv.
\end{proof}

{
\bibliographystyle{Chicago}
\bibliography{Reno_alt}

\begin{thebibliography}{}

\bibitem[\protect\citeauthoryear{Aliprantis and Border}{Aliprantis and
  Border}{1994}]{aliprantis}
Aliprantis, C. and K.~Border (1994).
\newblock {\em Infinite {D}imensional {A}nalysis}.
\newblock Springer-Verlag, Berlin.

\bibitem[\protect\citeauthoryear{Ambrosio, Gigli, and Savar{\'e}}{Ambrosio
  et~al.}{2008}]{ambrosiogigli}
Ambrosio, L., N.~Gigli, and G.~Savar{\'e} (2008).
\newblock {\em Gradient flows: in metric spaces and in the space of probability
  measures}.
\newblock Springer Science \& Business Media, Basel.

\bibitem[\protect\citeauthoryear{Bahlali, Mezerdi, and Djehiche}{Bahlali
  et~al.}{2006}]{bahlali}
Bahlali, S., B.~Mezerdi, and B.~Djehiche (2006).
\newblock Approximation and optimality necessary conditions in relaxed
  stochastic control problems.
\newblock {\em International Journal of Stochastic Analysis\/}~{\em 2006}.

\bibitem[\protect\citeauthoryear{Bene{\v{s}}}{Bene{\v{s}}}{1971}]{benevs}
Bene{\v{s}}, V. (1971).
\newblock Existence of optimal stochastic control laws.
\newblock {\em SIAM Journal on Control\/}~{\em 9\/}(3), 446--472.

\bibitem[\protect\citeauthoryear{Bertucci}{Bertucci}{2018}]{bertucci2018optimal}
Bertucci, C. (2018).
\newblock Optimal stopping in mean field games, an obstacle problem approach.
\newblock {\em Journal de Math{\'e}matiques Pures et Appliqu{\'e}es\/}~{\em
  120}, 165--194.

\bibitem[\protect\citeauthoryear{Billingsley}{Billingsley}{1999}]{billingsley}
Billingsley, P. (1999).
\newblock {\em Convergence of probability measures}.
\newblock John Wiley \& Sons, New York.

\bibitem[\protect\citeauthoryear{Bouveret, Dumitrescu, and Tankov}{Bouveret
  et~al.}{2020}]{bouveret2020mean}
Bouveret, G., R.~Dumitrescu, and P.~Tankov (2020).
\newblock Mean-field games of optimal stopping: a relaxed solution approach.
\newblock {\em SIAM Journal on Control and Optimization\/}~{\em 58\/}(4),
  1795--1821.

\bibitem[\protect\citeauthoryear{Brunick and Shreve}{Brunick and
  Shreve}{2013}]{brunick}
Brunick, G. and S.~Shreve (2013).
\newblock Mimicking an {I}t{\^o} process by a solution of a stochastic
  differential equation.
\newblock {\em The Annals of Applied Probability\/}~{\em 23\/}(4), 1584--1628.

\bibitem[\protect\citeauthoryear{Campi and Fischer}{Campi and
  Fischer}{2018}]{campi2018n}
Campi, L. and M.~Fischer (2018).
\newblock {N}-player games and mean-field games with absorption.
\newblock {\em The Annals of Applied Probability\/}~{\em 28\/}(4), 2188--2242.

\bibitem[\protect\citeauthoryear{Cardaliaguet}{Cardaliaguet}{2012}]{cardaliaguet}
Cardaliaguet, P. (2012).
\newblock Notes on mean field games (from {P}.{L}. {L}ions' lecture notes at
  {C}oll\`{e}ge de {F}rance).

\bibitem[\protect\citeauthoryear{Carmona and Delarue}{Carmona and
  Delarue}{2018}]{carmonadelarue}
Carmona, R. and F.~Delarue (2018).
\newblock {\em Probabilistic Theory of Mean Field Games with Applications
  I-II}.
\newblock Springer.

\bibitem[\protect\citeauthoryear{Carmona, Delarue, and Lacker}{Carmona
  et~al.}{2017}]{carmona2017mean}
Carmona, R., F.~Delarue, and D.~Lacker (2017).
\newblock Mean field games of timing and models for bank runs.
\newblock {\em Applied Mathematics \& Optimization\/}~{\em 76\/}(1), 217--260.

\bibitem[\protect\citeauthoryear{Carmona, Fouque, and Sun}{Carmona
  et~al.}{2015}]{carmona2015mean}
Carmona, R., J.-P. Fouque, and L.-H. Sun (2015).
\newblock Mean field games and systemic risk.
\newblock {\em Communications in Mathematical Sciences\/}~{\em 13\/}(4),
  911--933.

\bibitem[\protect\citeauthoryear{Carmona and Lacker}{Carmona and
  Lacker}{2015}]{carmonaweak}
Carmona, R. and D.~Lacker (2015).
\newblock A probabilistic weak formulation of mean field games and
  applications.
\newblock {\em The Annals of Applied Probability\/}~{\em 25\/}(3), 1189--1231.

\bibitem[\protect\citeauthoryear{Cellina}{Cellina}{1969}]{cellina}
Cellina, A. (1969).
\newblock Approximation of set valued functions and fixed point theorems.
\newblock {\em Annali di matem{\'a}tica pura ed applicata\/}~{\em 82\/}(1),
  17--24.

\bibitem[\protect\citeauthoryear{Chan and Sircar}{Chan and
  Sircar}{2015}]{chan2015bertrand}
Chan, P. and R.~Sircar (2015).
\newblock Bertrand and {C}ournot mean field games.
\newblock {\em Applied Mathematics \& Optimization\/}~{\em 71\/}(3), 533--569.

\bibitem[\protect\citeauthoryear{Chan and Sircar}{Chan and
  Sircar}{2017}]{chan2017fracking}
Chan, P. and R.~Sircar (2017).
\newblock Fracking, renewables, and mean field games.
\newblock {\em SIAM Review\/}~{\em 59\/}(3), 588--615.

\bibitem[\protect\citeauthoryear{Darling and Pardoux}{Darling and
  Pardoux}{1997}]{darling}
Darling, R. and E.~Pardoux (1997).
\newblock Backwards {SDE} with random terminal time and applications to
  semilinear elliptic {PDE}.
\newblock {\em The Annals of Probability\/}~{\em 25\/}(3), 1135--1159.

\bibitem[\protect\citeauthoryear{Delarue, Inglis, Rubenthaler, and
  Tanr{\'e}}{Delarue et~al.}{2015a}]{delarue2015global}
Delarue, F., J.~Inglis, S.~Rubenthaler, and E.~Tanr{\'e} (2015a).
\newblock Global solvability of a networked integrate-and-fire model of
  {M}c{K}ean--{V}lasov type.
\newblock {\em The Annals of Applied Probability\/}~{\em 25\/}(4), 2096--2133.

\bibitem[\protect\citeauthoryear{Delarue, Inglis, Rubenthaler, and
  Tanr{\'e}}{Delarue et~al.}{2015b}]{delarue2015particle}
Delarue, F., J.~Inglis, S.~Rubenthaler, and E.~Tanr{\'e} (2015b).
\newblock Particle systems with a singular mean-field self-excitation.
  application to neuronal networks.
\newblock {\em Stochastic Processes and their Applications\/}~{\em 125\/}(6),
  2451--2492.

\bibitem[\protect\citeauthoryear{Dembo and Zeitouni}{Dembo and
  Zeitouni}{2010}]{dembo}
Dembo, A. and O.~Zeitouni (2010).
\newblock {\em Large {D}eviations {T}echniques and {A}pplications. {S}tochastic
  {M}odelling and {A}pplied {P}robability, 38}.
\newblock Springer-Verlag, Berlin.

\bibitem[\protect\citeauthoryear{Dufour and Stockbridge}{Dufour and
  Stockbridge}{2012}]{dufour}
Dufour, F. and R.~Stockbridge (2012).
\newblock On the existence of strict optimal controls for constrained,
  controlled {M}arkov processes in continuous time.
\newblock {\em Stochastics: An International Journal of Probability and
  Stochastic Processes\/}~{\em 84\/}(1), 55--78.

\bibitem[\protect\citeauthoryear{El~Karoui, Nguyen, and
  Jeanblanc-Picqu{\'e}}{El~Karoui et~al.}{1987}]{karoui}
El~Karoui, N., D.~Nguyen, and M.~Jeanblanc-Picqu{\'e} (1987).
\newblock Compactification methods in the control of degenerate diffusions:
  existence of an optimal control.
\newblock {\em Stochastics: An International Journal of Probability and
  Stochastic Processes\/}~{\em 20\/}(3), 169--219.

\bibitem[\protect\citeauthoryear{Filippov}{Filippov}{1962}]{filippov}
Filippov, A. (1962).
\newblock On certain questions in the theory of optimal control.
\newblock {\em Journal of the Society for Industrial and Applied Mathematics,
  Series A: Control\/}~{\em 1\/}(1), 76--84.

\bibitem[\protect\citeauthoryear{Fleming and Rishel}{Fleming and
  Rishel}{2012}]{fleming}
Fleming, W. and R.~Rishel (2012).
\newblock {\em Deterministic and stochastic optimal control}, Volume~1.
\newblock Springer Science \& Business Media, New York.

\bibitem[\protect\citeauthoryear{Fouque and Sun}{Fouque and
  Sun}{2013}]{fouque2013systemic}
Fouque, J.-P. and L.-H. Sun (2013).
\newblock Systemic risk illustrated.
\newblock In {\em Handbook on Systemic Risk}, pp.\  444--452. Cambridge
  University Press, Cambridge.

\bibitem[\protect\citeauthoryear{Funaki}{Funaki}{1984}]{funaki1984certain}
Funaki, T. (1984).
\newblock A certain class of diffusion processes associated with nonlinear
  parabolic equations.
\newblock {\em Zeitschrift f{\"u}r Wahrscheinlichkeitstheorie und Verwandte
  Gebiete\/}~{\em 67\/}(3), 331--348.

\bibitem[\protect\citeauthoryear{G{\"a}rtner}{G{\"a}rtner}{1988}]{gartner}
G{\"a}rtner, J. (1988).
\newblock On the {M}c{K}ean--{V}lasov limit for interacting diffusions.
\newblock {\em Mathematische Nachrichten\/}~{\em 137\/}(1), 197--248.

\bibitem[\protect\citeauthoryear{Giesecke, Spiliopoulos, and Sowers}{Giesecke
  et~al.}{2013}]{giesecke2013default}
Giesecke, K., K.~Spiliopoulos, and R.~Sowers (2013).
\newblock Default clustering in large portfolios: Typical events.
\newblock {\em The Annals of Applied Probability\/}~{\em 23\/}(1), 348--385.

\bibitem[\protect\citeauthoryear{Giesecke, Spiliopoulos, Sowers, and
  Sirignano}{Giesecke et~al.}{2015}]{giesecke2015large}
Giesecke, K., K.~Spiliopoulos, R.~Sowers, and J.~Sirignano (2015).
\newblock Large portfolio asymptotics for loss from default.
\newblock {\em Mathematical Finance\/}~{\em 25\/}(1), 77--114.

\bibitem[\protect\citeauthoryear{Grigelionis and Mackevi{\v{c}}ius}{Grigelionis
  and Mackevi{\v{c}}ius}{2003}]{grigelionis}
Grigelionis, B. and V.~Mackevi{\v{c}}ius (2003).
\newblock The finiteness of moments of a stochastic exponential.
\newblock {\em Statistics \& probability letters\/}~{\em 64\/}(3), 243--248.

\bibitem[\protect\citeauthoryear{Hambly and Ledger}{Hambly and
  Ledger}{2017}]{hambly2017stochastic}
Hambly, B. and S.~Ledger (2017).
\newblock A stochastic {M}c{K}ean--{V}lasov equation for absorbing diffusions
  on the half-line.
\newblock {\em The Annals of Applied Probability\/}~{\em 27\/}(5), 2698--2752.

\bibitem[\protect\citeauthoryear{Hambly, Ledger, and S{\o}jmark}{Hambly
  et~al.}{2019}]{hambly2019mckean}
Hambly, B., S.~Ledger, and A.~S{\o}jmark (2019).
\newblock A {M}c{K}ean--{V}lasov equation with positive feedback and blow-ups.
\newblock {\em The Annals of Applied Probability\/}~{\em 29\/}(4), 2338--2373.

\bibitem[\protect\citeauthoryear{Hambly and S{\o}jmark}{Hambly and
  S{\o}jmark}{2019}]{hambly2019spde}
Hambly, B. and A.~S{\o}jmark (2019).
\newblock An {SPDE} model for systemic risk with endogenous contagion.
\newblock {\em Finance and Stochastics\/}~{\em 23\/}(3), 535--594.

\bibitem[\protect\citeauthoryear{Haussmann and Lepeltier}{Haussmann and
  Lepeltier}{1990}]{haussmann}
Haussmann, U. and J.~Lepeltier (1990).
\newblock On the existence of optimal controls.
\newblock {\em SIAM Journal on Control and Optimization\/}~{\em 28\/}(4),
  851--902.

\bibitem[\protect\citeauthoryear{Huang, Malham{\'e}, and Caines}{Huang
  et~al.}{2006}]{huang2006}
Huang, M., R.~Malham{\'e}, and P.~Caines (2006).
\newblock Large population stochastic dynamic games: closed-loop
  {M}c{K}ean--{V}lasov systems and the {N}ash certainty equivalence principle.
\newblock {\em Communications in Information \& Systems\/}~{\em 6\/}(3),
  221--252.

\bibitem[\protect\citeauthoryear{Jacod and Shiryaev}{Jacod and
  Shiryaev}{2013}]{jacod}
Jacod, J. and A.~Shiryaev (2013).
\newblock {\em Limit theorems for stochastic processes}, Volume 288.
\newblock Springer Science \& Business Media, Berlin.

\bibitem[\protect\citeauthoryear{Kallenberg}{Kallenberg}{2006}]{kallenberg}
Kallenberg, O. (2006).
\newblock {\em Foundations of modern probability}.
\newblock Springer Science \& Business Media, New York.

\bibitem[\protect\citeauthoryear{Karatzas and Shreve}{Karatzas and
  Shreve}{1987}]{karaztas}
Karatzas, I. and S.~Shreve (1987).
\newblock {\em Brownian motion and Stochastic calculus, {V}olume 113 of}.
\newblock Graduate Texts in Mathematics, New York.

\bibitem[\protect\citeauthoryear{Kushner and Dupuis}{Kushner and
  Dupuis}{2013}]{kushner}
Kushner, H. and P.~Dupuis (2013).
\newblock {\em Numerical methods for stochastic control problems in continuous
  time}, Volume~24.
\newblock Springer Science \& Business Media, New York.

\bibitem[\protect\citeauthoryear{Lacker}{Lacker}{2015}]{lacker}
Lacker, D. (2015).
\newblock Mean field games via controlled martingale problems: existence of
  {M}arkovian equilibria.
\newblock {\em Stochastic Processes and their Applications\/}~{\em 125\/}(7),
  2856--2894.

\bibitem[\protect\citeauthoryear{Lacker}{Lacker}{2018}]{lacker2018strong}
Lacker, D. (2018).
\newblock On a strong form of propagation of chaos for {M}c{K}ean--{V}lasov
  equations.
\newblock {\em {E}lectronic {C}ommunications in {P}robability\/}~{\em 23}.

\bibitem[\protect\citeauthoryear{Lacker}{Lacker}{2020}]{lacker2020convergence}
Lacker, D. (2020).
\newblock On the convergence of closed-loop nash equilibria to the mean field
  game limit.
\newblock {\em Annals of Applied Probability\/}~{\em 30\/}(4), 1693--1761.

\bibitem[\protect\citeauthoryear{Lang}{Lang}{1986}]{lang}
Lang, R. (1986).
\newblock A note on the measurability of convex sets.
\newblock {\em Archiv der Mathematik\/}~{\em 47\/}(1), 90--92.

\bibitem[\protect\citeauthoryear{Lasry and Lions}{Lasry and
  Lions}{2006a}]{lasry-lions2006a}
Lasry, J.-M. and P.-L. Lions (2006a).
\newblock Jeux {\`a} champ moyen. i--le cas stationnaire.
\newblock {\em Comptes Rendus Math{\'e}matique\/}~{\em 343\/}(9), 619--625.

\bibitem[\protect\citeauthoryear{Lasry and Lions}{Lasry and
  Lions}{2006b}]{lasry-lions2006b}
Lasry, J.-M. and P.-L. Lions (2006b).
\newblock Jeux {\`a} champ moyen. ii--horizon fini et contr{\^o}le optimal.
\newblock {\em Comptes Rendus Math{\'e}matique\/}~{\em 343\/}(10), 679--684.

\bibitem[\protect\citeauthoryear{Lasry and Lions}{Lasry and
  Lions}{2007}]{lasry-lions2007}
Lasry, J.-M. and P.-L. Lions (2007).
\newblock Mean field games.
\newblock {\em Japanese journal of mathematics\/}~{\em 2\/}(1), 229--260.

\bibitem[\protect\citeauthoryear{Ledger and S{\o}jmark}{Ledger and
  S{\o}jmark}{2020}]{ledger2020uniqueness}
Ledger, S. and A.~S{\o}jmark (2020).
\newblock Uniqueness for contagious {M}c{K}ean--{V}lasov systems in the weak
  feedback regime.
\newblock {\em Bulletin of the London Mathematical Society\/}~{\em 52\/}(3),
  448--463.

\bibitem[\protect\citeauthoryear{McKean}{McKean}{1966}]{mckean}
McKean, H. (1966).
\newblock A class of {M}arkov processes associated with nonlinear parabolic
  equations.
\newblock {\em Proceedings of the National Academy of Sciences\/}~{\em
  56\/}(6), 1907--1911.

\bibitem[\protect\citeauthoryear{M{\'e}l{\'e}ard}{M{\'e}l{\'e}ard}{1996}]{meleard1996asymptotic}
M{\'e}l{\'e}ard, S. (1996).
\newblock Asymptotic behaviour of some interacting particle systems;
  {M}c{K}ean-{V}lasov and {B}oltzmann models.
\newblock In {\em Probabilistic models for nonlinear partial differential
  equations}, pp.\  42--95. Springer.

\bibitem[\protect\citeauthoryear{Nadtochiy and Shkolnikov}{Nadtochiy and
  Shkolnikov}{2019}]{nadtochiy2019particle}
Nadtochiy, S. and M.~Shkolnikov (2019).
\newblock Particle systems with singular interaction through hitting times:
  application in systemic risk modeling.
\newblock {\em The Annals of Applied Probability\/}~{\em 29\/}(1), 89--129.

\bibitem[\protect\citeauthoryear{Nadtochiy and Shkolnikov}{Nadtochiy and
  Shkolnikov}{2020}]{nadtochiy2020mean}
Nadtochiy, S. and M.~Shkolnikov (2020).
\newblock Mean field systems on networks, with singular interaction through
  hitting times.
\newblock {\em Annals of Probability\/}~{\em 48\/}(3), 1520--1556.

\bibitem[\protect\citeauthoryear{Nutz}{Nutz}{2018}]{nutz2018mean}
Nutz, M. (2018).
\newblock A mean field game of optimal stopping.
\newblock {\em SIAM Journal on Control and Optimization\/}~{\em 56\/}(2),
  1206--1221.

\bibitem[\protect\citeauthoryear{Oelschl{\"a}ger}{Oelschl{\"a}ger}{1984}]{oelschlager1984martingale}
Oelschl{\"a}ger, K. (1984).
\newblock A martingale approach to the law of large numbers for weakly
  interacting stochastic processes.
\newblock {\em The Annals of Probability\/}~{\em 12}, 458--479.

\bibitem[\protect\citeauthoryear{Stroock and Varadhan}{Stroock and
  Varadhan}{1969}]{stroock1969}
Stroock, D. and S.~Varadhan (1969).
\newblock Diffusion processes with continuous coefficients, i-ii.
\newblock {\em Communications on Pure and Applied Mathematics\/}~{\em 22\/}(3),
  345--400.

\bibitem[\protect\citeauthoryear{Stroock and Varadhan}{Stroock and
  Varadhan}{2007}]{stroock}
Stroock, D. and S.~Varadhan (2007).
\newblock {\em Multidimensional diffusion processes}.
\newblock Springer, Berlin.

\bibitem[\protect\citeauthoryear{Sznitman}{Sznitman}{1991}]{sznitman}
Sznitman, A. (1991).
\newblock Topics in propagation of chaos.
\newblock In {\em Ecole d'{\'e}t{\'e} de probabilit{\'e}s de Saint-Flour XIX
  1989}, pp.\  165--251. Springer.

\bibitem[\protect\citeauthoryear{Villani}{Villani}{2003}]{villani}
Villani, C. (2003).
\newblock {\em Topics in optimal transportation}.
\newblock Number~58. American Mathematical Society, Providence.

\end{thebibliography}
}
\end{document}